\documentclass[10pt]{article}

\usepackage{amscd,amssymb,amsmath,latexsym,enumerate}
\usepackage{hyperref}
\usepackage{paralist}
\usepackage[mathscr]{euscript}
\usepackage{epsfig}
\usepackage{fancybox}
\usepackage{yhmath}
\usepackage{verbatim}
\usepackage{graphicx,bm,units,yfonts}
\usepackage{phonetic}
\usepackage{mathtools}
\usepackage[small,nohug]{diagrams}
\diagramstyle[labelstyle=\scriptstyle]
\usepackage{xcolor}
\usepackage{upgreek}
\usepackage{makeidx}
\usepackage{graphicx}

\newtheorem{theorem}{Theorem}[section]

\numberwithin{equation}{section}
\newtheorem{proposition}[theorem]{Proposition}
\newtheorem{lemma}[theorem]{Lemma}
\newtheorem{corollary}[theorem]{Corollary}
\newtheorem{definition}[theorem]{Definition}
\newtheorem{remark}[theorem]{Remark}
\newtheorem{example}[theorem]{Example}

\newcommand{\BM}{{\mathbb B}}
\newcommand{\CM}{{\mathbb C}}
\newcommand{\FM}{{\mathbb F}}
\newcommand{\NM}{{\mathbb N}}
\newcommand{\PM}{{\mathbb P}}
\newcommand{\RM}{{\mathbb R}}
\newcommand{\SM}{{\mathbb S}}
\newcommand{\TM}{{\mathbb T}}
\newcommand{\ZM}{{\mathbb Z}}
\newcommand{\KM}{{\mathbb K}}
\newcommand{\QM}{{\mathbb Q}}
\newcommand{\UM}{{\mathbb U}}
\newcommand{\EM}{{\mathbb E}}

\newcommand{\Aa}{{\mathcal A}}
\newcommand{\Bb}{{\mathcal B}}
\newcommand{\Cc}{{\mathcal C}}
\newcommand{\Dd}{{\mathcal D}}
\newcommand{\Hh}{{\mathcal H}}
\newcommand{\Kk}{{\mathcal K}}
\newcommand{\Tt}{{\mathcal T}}

\newcommand{\Uu}{{\mathcal U}}
\newcommand{\Mm}{{\mathcal M}}
\newcommand{\Nn}{{\mathcal N}}
\newcommand{\Pp}{{\mathcal P}}

\newcommand{\Ind}{\mbox{\rm Ind}}
\newcommand{\Exp}{\mbox{\rm Exp}}
\newcommand{\sgn}{{\rm sgn}}
\newcommand{\Ad}{{\rm Ad}}
\newcommand{\diag}{{\rm diag}}
\newcommand{\ev}{{\rm ev}}
\newcommand{\Tr}{\mbox{\rm Tr}}

\newcommand{\sa}{{\mbox{\rm\tiny sa}}}
\newcommand{\gu}{{\mbox{\rm\tiny u}}}
\newcommand{\sym}{{\mbox{\rm\tiny sym}}}

\newcommand{\proof}{\noindent {\bf Proof. }}
\newcommand{\qed}{\hfill $\Box$}

\newcommand{\DiracSA}{H}

\begin{document}

\title{Generalized Connes-Chern characters in $KK$-theory \\
with an application to weak invariants of topological insulators}

\author{Emil Prodan$^1$, Hermann Schulz-Baldes$^2$
\\
\\
{\small $^1$ Department of Physics and Department of Mathematical Sciences, Yeshiva University, USA}
\\
{\small $^2$ Department Mathematik, Friedrich-Alexander-Universit\"at Erlangen-N\"urnberg, Germany}
}
\date{ }

\maketitle

\begin{abstract}
We use constructive bounded Kasparov $K$-theory to investigate the numerical invariants stemming from the internal Kasparov products $K_i(\mathcal A) \times KK^i(\mathcal A, \mathcal B) \rightarrow K_0(\mathcal B) \rightarrow \RM$, $i=0,1$, where the last morphism is provided by a tracial state.  For the class of properly defined finitely-summable Kasparov $(\Aa,\Bb)$-cycles, the invariants are given by the pairing of $K$-theory of $\Bb$ with an element of the periodic cyclic cohomology of $\Bb$, which we call the generalized Connes-Chern character. When $\Aa$ is a twisted crossed product of $\Bb$ by $\ZM^k$, $\Aa = \Bb \rtimes_\xi^\theta \ZM^k$, we derive a local formula for the character corresponding to the fundamental class of a properly defined Dirac cycle. Furthermore, when $\Bb = C(\Omega) \rtimes_{\xi'}^{\phi} \ZM^j$, with $C(\Omega)$ the algebra of continuous functions over a disorder configuration space, we show that the numerical invariants are connected to the weak topological invariants of the complex classes of topological insulators, defined in the physics literature. The end products are generalized index theorems for these weak invariants, which enable us to predict the range of the invariants and to identify regimes of strong disorder in which the invariants remain stable. The latter will be reported in a subsequent publication.
\end{abstract}

\setcounter{tocdepth}{5}

\tableofcontents

\section{Introduction}

As Jensen and Thomsen state in the preface of their textbook ``Elements of $KK$-theory" \cite{JensenThomsenBook1991}, Kasparov's $KK$-theory \cite{KasparovIzv1975tr,KasparovJOT1980yt,KasparovMUI1981ui,KasparovJSM1987} is a powerful but forbiddingly difficult topic, encompassing several traditional branches of mathematics. With concepts such as $KK$-groups and Kasparov products, it provides the framework for a vast generalization of traditional $K$-theory of operator algebras. Unfortunately, many of the fundamental statements are not entirely constructive and this aspect has hampered a wider adoption of $KK$-theoretic methods and the applications remained limited for a while. As already emphasized in \cite{BKR}, this seems to be changing rapidly as recent years have seen a swift progress on what is now called the constructive Kasparov product \cite{BMS,FR,KL,KNR,Me} and on its associated numerical invariants \cite{CareyAMS2014jf,CPRS1,CPRS2}. In turn, this enabled several important applications in condensed matter \cite{Bou,BCR,BCR2,BKR,Pro3} as well as high-energy physics \cite{BMS,Su}. Our work fits into the applied side of this general program, with the caveat that our approach relies entirely on the bounded Kasparov theory, as opposed to most if not all the existing new developments which involve the un-bounded version of the theory \cite{BJ1,Ku}. Let us specify from the beginning, though, that the focus is rather on deriving local index formulas for numerical invariants derived from particular Kasparov products. More precisely, we seek a local index formula for the following numerical invariant
\begin{equation}\label{Diag-MainIntro1}
\begin{diagram}
& KK^\ast(\CM,\Aa) & \; \times \; & KK^\ast(\Aa,\Bb)  & \;\rTo{\rm Kasparov \ product} \;&  K_0(\Bb) & \;\rTo{\ \ {\rm Tr} \ \ } \;&  \RM \; .
\end{diagram}
\end{equation}
Here ${\rm Tr}$ denotes a tracial state on matrices with entries in $\Bb$. In principle one could also obtain numerical invariants in another way ({\it e.g.} pairings with higher cocycles), but this will not be studied here. As the reader perhaps immediately realizes (see also our discussion in Section~\ref{sec-finalrem}), \eqref{Diag-MainIntro1} is a straightforward generalization of the numerical invariant
\begin{equation}\label{Diag-MainIntro2}
\begin{diagram}
& KK^\ast(\CM,\Aa) & \; \times \; & KK^\ast(\Aa,\CM)  & \;\rTo{\rm Kasparov \ product} \;&  K_0(\CM) & \;\rTo{\ \ {\rm Tr} \ \ } \;&  \ZM \; ,
\end{diagram}
\end{equation}
which is at the heart of Connes' quantized calculus and defines the Connes-Chern character in the periodic cyclic homology and cohomology \cite{Connes:1994wk}. From a $KK$-theory viewpoint, the quantized calculus together with the local index formula of Connes-Moscovici \cite{CONNES:1995rv} provide a constructive framework for the Kasparov products in \eqref{Diag-MainIntro2}. That this is a fruitful approach to index theorems was already stressed in Kasparov's works and further advocated by Higson and Roe \cite{HR,HigsonPSPM1990} (see also \cite{Pus}). This point of view, for example, led to the index theory for locally compact noncommutative geometries \cite{CareyAMS2014jf}. 

\vspace{.2cm}

As in past works \cite{BELLISSARD:1994xj,ProdanJPA2013hg,PS,ProdanSpringer2016ds}, our interest in index theorems is driven by applications to solid state systems with topological properties. In this work we prove generalized index theorems for the so-called weak numerical invariants distinguishing weak phases of complex topological insulators. The algebra of the physical observables for these systems is generated by the non-commutative tori as well as algebras of continuous functions over disorder configuration spaces. The latter is assumed to be compact and metrizable. In this context, a condensed matter system is defined by a self-adjoint Hamiltonian $h$ from the stabilization of this algebra and its topology is encoded in either the Fermi projection or the Fermi unitary operator (see Section~\ref{Sec-ApplicationTI}). Both operators can be constructed by functional calculus from $h$, {\it e.g.} using $\sgn(h-\mu)$ where $\mu$ denotes the Fermi energy. In the regime of weak disorder, the Fermi energy is located inside a spectral gap of $h$ and these operators belong to a smooth sub-algebra of the $C^\ast$-algebra of physical observables. Hence the condensed matter phases are classified by its $K$-theory in this regime. When the disorder configuration space is contractible, a complete set of numerical invariants can be obtained by pairing the generators of the cyclic cohomology with the class of the topological phase in the $K$-theory.  As shown in \cite{ProdanSpringer2016ds}, the knowledge of the complete set of numerical invariants enables one to pin-point the $K$-theory class of a topological phase and each numerical invariant displays a bulk-boundary correspondence principle. The physical effect of the latter is that certain surfaces cut to a crystal induce essential spectrum at the Fermi level, and this spectrum cannot be removed by any local boundary condition. The strong invariant, {\it i.e.} the one obtained from the highest generating cyclic cocycle, was in \cite{BELLISSARD:1994xj,ProdanJPA2013hg,PS} shown to be stable in the regime of strong disorder where the Fermi level is located in the essential, but Anderson-localized spectrum of $h$. In this regime, the operators encoding the topology are no longer inside the $C^\ast$-algebra, but instead they belong to a non-commutative Sobolev space. Furthermore, it was shown in \cite{ProdanSpringer2016ds} that the boundary spectrum induced by the strong invariant is not Anderson-localized in the presence of boundary disorder. Key for these developments were the index theorems stemming from \eqref{Diag-MainIntro2}, together with derivations of local index formulas for bounded Fredholm modules. For the weak invariants, the questions of whether they are stable in the regime of strong disorder and if the topological boundary spectrum can be Anderson-delocalized remain completely open, despite of numerical evidence that there are situations where the answer is affirmative \cite{BF,KOI,MBM,RKS}. Although this work is restricted to the analysis of the weak disorder regime, we can announce here that the generalized index theorems stemming from \eqref{Diag-MainIntro2} enabled us to identify such strong disorder regimes and that these findings will be object of a subsequent paper.

\section{Outline and relation to other works}

Parts of this work are review. It was our aim to collect most of the notions and known results on $KK$-theoy needed for the study of topological invariants in solid state systems. Unfortunately, they cannot be found in a single reference, but are rather scattered over the literature of the last decades, with some only dating back a few years. We hope that the detailed treatment in Sections~\ref{chap-GradedHilbert}, \ref{chap-ComplexK} and \ref{chap-KK} will be useful for readers interested in applications of constructive $KK$-theory (which we expect to be far reaching). Not all, but many of the central arguments are included, in particular, those that we could not locate in the literature ({\it e.g.} Section~\ref{sec-K1KK0}). In any case, we took great care to provide detailed referencing. 

\vspace{0.2cm}

Sections \ref{chap-GCCC}, \ref{Sec-CrossedProd} and \ref{Sec-ApplicationTI} then mainly contain new material.  We develop a constructive theory of \eqref{Diag-MainIntro1} in the bounded version of Kasparov's theory without leaving the $C^\ast$-algebraic setting. This is a distinct feature of the present approach. For this purpose, the notion of finitely summable Kasparov cycles is used as an input to a generalization of the Calderon-Fedosov formula (Section~\ref{Sec-Fedosov}) for the index map in $K$-theory. We believe this formula is new in the $C^\ast$-algebraic setting. Then it is shown (Sections~\ref{Sec-EvenGChern} and \ref{Sec-OddGChern}) that \eqref{Diag-MainIntro1} is given by a direct lift of the classical Connes-Chern character to the present context. Furthermore (Section~\ref{Sec-CrossedProd}), when $\Aa$ in \eqref{Diag-MainIntro1} is a twisted crossed product of a $C^*$-algebra $\Bb$ by $\ZM^k$, we derive a local formula for the generalized Connes-Chern character corresponding to the fundamental class generated by the Dirac operator constructed by pairing of the derivations of an associated torus action with the generators of an appropriate complex Clifford algebra. The use of such auxiliary Clifford algebras already appeared in the work on the integer quantum Hall effect \cite{BELLISSARD:1994xj} and it was also successfully employed in the works of Pask, Rennie and Sims on graph algebras \cite{PR1,PR2}.  Then the final step of the calculation of the generalized index theorems for crossed product algebras is based on certain geometric identities which were already the key to \cite{BELLISSARD:1994xj,ProdanJPA2013hg,PS,ProdanSpringer2016ds}. This is an important feature of our work, as there are already strong indications that these identities pertain to extensions into the strong disorder regime. Lastly, in Section~\ref{Sec-ApplicationTI} we show that in the context of complex topological insulators, these local index formulas coincide with the bulk weak topological invariants, as defined in the physics literature. These local index formulas also cover the boundary invariants connected to the bulk weak invariants by the bulk-boundary principle \cite{ProdanSpringer2016ds}, but this is not explicitly discussed here.

\vspace{0.2cm}

The present work does not constitute the first application of $KK$-theory to solid state systems. The effectiveness of $KK$-theory in this context was first pointed out in \cite{BCR}, where the bulk-boundary correspondence  for Class A topological insulators in dimension 2 was reformulated in terms of the Kasparov product. A classification of the strong bulk-topological phases in the framework of $KK$-theory was developed in \cite{Bou, BCR2}, and in \cite{BKR} the results were generalized to the whole periodic table of topological insulators. Note however that these works cover only the strong topological phases and only the weak disorder regime. For this reason, no statement could be made about the central issue in the field of topological insulators, namely, the delocalized character of the boundary spectrum. Independently and at the same time, the unpublished manuscript \cite{Pro3} of one of the authors pointed out that $KK$-theory can be used to derive index theorems for the weak invariants. From a physical point of view, this development is important because, if these theorems can be pushed in regimes of strong disorder, then the delocalization of the boundary states induced by the weak invariants can be proven by following the arguments used in \cite{ProdanSpringer2016ds} for the strong invariants. 

\vspace{0.2cm}

Apart from the present approach based on the geometric identities alluded to above, it seems feasible to achieve index theorems for the weak invariants by evaluating the  Connes-Moscovici local index formula \cite{CONNES:1995rv,HigsonICTP}, more precisely its extensions to unital \cite{CPRS1,CPRS2,CPRS3} and non-unital \cite{CareyAMS2014jf} semifinite von Neumann algebras. We will refer to the later as the semifinite local index formula. In the unital odd case, for example, given a $\ast$-subalgebra $\Aa$ of a semifinite von Neumann algebra $(\Nn,\tau)$ and a semifinite spectral triple $(\Aa,\Hh,\Dd)$, smooth and summable, the von Neumann spectral flow Sf$(\Dd,u \Dd u^\ast)$ was interpreted in \cite{CPRS1} as the pairing $\langle [u]_1,[(\Aa,\Hh,\Dd)] \rangle$ between the $K$-theory and $K$-cohomology of the completion of $\Aa$ in the so-called $\delta$-$\phi$-topology. In this setting, the calculation of the semifinite local index formula is initiated from the integral formula for the spectral flow \cite{CP1,CP2}, which relies on the Breuer index for semifinite von Neumann algebras \cite{Bre1,Bre2}. This result is a general and yet constructive tool, which already found a variety of applications (see the last chapter of \cite{CareyAMS2014jf}). In particular, in Section~5.1 of \cite{CareyAMS2014jf}, it was pointed out that the local index formulas simplifies considerably in the broad context of possibly non-unital $C^\ast$-algebras equipped with a continuous torus action (the unital case is essentially contained in \cite{PR2}).

\vspace{0.2cm}

A delicate point of the semifinite index theory is the choice and the control of the various subalgebras involved, in particular, the algebra $\Aa$ above. This affects the range of the semifinite index pairing.  In general, it is possible to use a specific $C^\ast$-subalgebra of $\Kk_\Nn$, the algebra of $\tau$-compact operators in $\Nn$ (see \cite[p.~37]{CareyAMS2014jf}). An important observation is that this subalgebra is separable whenever $\Aa$ is separable, in which case the semifinite index lands in a discrete subgroup of the real numbers. Unfortunately this separable algebra is rigidly fixed by the algebra $\Aa$ and its $K$-theory is out of reach in most cases so that the range of the index is not easy to predict. In a further development, \cite{KNR} introduced the notion of a spectral flow relative to any normed-closed ideal $J$ of $\Nn$ and then a $J$-semifinite index. This extension appeals to the theory of $J$-Fredholm operators and $J$-indices \cite{Ols}, which are natural extensions of the Breuer's work. In this new form, the theory seems to give more freedom on choosing the algebra $\Bb$ in \eqref{Diag-MainIntro1} and, equally important, it was also shown that the von Neumann spectral flow can factor through a $C^\ast$-spectral flow (see Theorems~5.2 and 5.5 in \cite{KNR}). However, none of the existing results in the literature seem to cover the results proved here. Even if the technical problems linked to the choice of subalgebras can be overcome, a supplementary calculation is needed to evaluate the semifinite local index formula and connect it to the weak invariants. The present approach is considerably shorter and more direct. Furthermore, it can be extended to the regime of strong disorder. The difficulties encountered when trying to push the semifinite local index formula into this regime are discussed in \cite[p. 77]{Bou}. A further potential advantage of the present approach is that is easier to use a higher  cocycle in the last arrow of \eqref{Diag-MainIntro1} instead of a tracial state, though this is not worked out here.

\vspace{0.2cm}





\section{Graded Hilbert C$^\ast$-modules} 
\label{chap-GradedHilbert}

In an attempt of making the exposition self-contained, this section develops the minimal background on Hilbert C$^\ast$-modules with a particular focus on gradings, which are natural and central to Kasparov's $K$-theory but no so much in the other $K$-theories. Everything said here can be found in the original papers of Kasparov \cite{KasparovIzv1975tr,KasparovJOT1980yt,KasparovMUI1981ui,KasparovJSM1987} and in the standard textbooks on the subject, such as \cite{BlackadarBook1998,JensenThomsenBook1991,Ska,WeggeOlsenBook1993de,HR}. Other excellent sources are \cite{PaschkeTAMS1973rd, LanceBook1995vc}, especially for the first part of the exposition. Regarding the presentation, we chose to introduce the concepts using definitions immediately followed by examples. Most of these examples will play an important role later in the exposition, hence the they also serve us to set up our notation. Together with the accompanying remarks and statements, they gradually prepare the ground for the main developments in later chapters.

\subsection{Graded C$^\ast$-algebras} 

As a preparation for the definition of graded Hilbert C$^\ast$-modules, let us first recall the notion of grading on C$^\ast$-algebras. Throughout, the algebras will be denoted by calligraphic capital letters from the beginning of the alphabet, such as $\Aa$, $\Bb$, etc..
 
\begin{definition}[\cite{BlackadarBook1998} pp. 114, and \cite{HR}, Appendix A] Let $\Bb$ be a C$^\ast$-algebra.
\begin{enumerate}[\rm (i)]
\item The algebra is said to be graded if it is equipped with an order two $\ast$-automorphism $\gamma_\mathcal B$, that is, for all $b,b'\in\Bb$,
$$
\gamma_\mathcal B(bb') 
\;=\; 
\gamma_\mathcal B(b)\gamma_\mathcal B(b')
\;, 
\qquad 
\gamma_\mathcal B(b^\ast) 
\;=\; 
\gamma_\mathcal B(b)^\ast
\;, 
\qquad 
\gamma_\mathcal B \circ \gamma_\mathcal B\;=\;{\rm id}
\;.
$$
In this case, one says that $\mathcal B$ is graded by the grading automorphism $\gamma_\mathcal B$. 

\item An element $b\in\Bb$ of a graded C$^\ast$-algebra is called even if $\gamma_\Bb(b)=b$, and odd if $\gamma_\Bb(b)=-b$. Both even and odd elements are called homogeneous.

\item The grading is said to be even (sometimes also called inner) if there exists a self-adjoint unitary $g$ in the multiplier algebra $\Mm(\Bb)$ of $\Bb$ such that $\gamma_\Bb(b) = {\rm Ad}_g(b)= g^\ast b g$. Otherwise, the grading is said to be odd. 

\item A graded homomorphism $\phi : \mathcal A \rightarrow \mathcal B$ between graded C$^\ast$-algebras is a $\ast$-homomorphism commuting with the gradings,
$$
\phi \circ \gamma_\mathcal A \;=\; \gamma_\mathcal B \circ \phi
\;.
$$

\item A homomorphism is said to be non-degenerate if for any non-zero element $b$ of $\Bb$, there exists $a \in \Aa$ such that $\phi(a)b \neq 0$.
\end{enumerate}
\end{definition}

\begin{remark}{\rm In the case when the algebras are unital, which is mostly the case considered in these notes, non-degeneracy is equivalent with $\phi(1)=1$.
}
\hfill $\diamond$
\end{remark}

\begin{example} 
{\rm
An  ungraded C$^\ast$-algebra can and will be regarded as a graded one with the trivial grading given by the identity map.
}
\hfill $\diamond$
\end{example}

\begin{example} 
\label{ex-MatrixAlgGrading}
{\rm Throughout, we will denote the C$^\ast$-algebra of $N\times N$-matrices with complex entries by $\CM(N)$. In particular, $\CM(2)$ will appear quite often in our exposition. A useful set of (linear) generators is provided by the identity matrix together with Pauli's matrices
$$
\sigma_1\;=\;\begin{pmatrix} 0 & 1 \\ 1 & 0 \end{pmatrix}\;, 
\qquad 
\sigma_2\;=\;\begin{pmatrix} 0 & -\imath \\ \imath & 0 \end{pmatrix}
\;, 
\qquad 
\sigma_3\;=\;\begin{pmatrix} 1 & 0 \\ 0 & -1 \end{pmatrix}
\;.
$$
Most of the time, $\CM(2)$ will be considered with the inner grading ${\rm Ad}_{\sigma_3}$. In general, $\CM(2N)$ can be equipped with a canonical inner grading provided by the isomorphism $\CM(2N) \cong \CM(2) \otimes \CM(N)$ and by the grading ${\rm Ad}_{\sigma_3} \otimes 1$ on the latter algebra.   
}
\hfill $\diamond$
\end{example}

\begin{example} 
\label{Ex-CliffordAlg2}
{\rm
Let $\CM_k$ be the complex Clifford algebra generated by $k$ elements $\Gamma_1,\ldots,\Gamma_k$ satisfying
$$
\Gamma_i \Gamma_j + \Gamma_j \Gamma_j 
\;=\; 2 \,\delta_{i,j}\;, 
\qquad 
\Gamma_i^\ast \;=\; \Gamma_i\;, 
\qquad i,j = 1,\ldots, k
\;.
$$
Its standard grading is given by
$$
\gamma_{\CM_k}(\Gamma_j)\;=\;-\,\Gamma_j
\;,
\qquad
j=1,\ldots,k
\;.
$$
For $k$ even, the grading is even and  given by $\gamma_{\CM_k} = {\rm Ad}_{\Gamma_0}$, where the so-called chiral element is defined by
$$
\Gamma_0 \;= \;(-\imath)^\frac{k}{2} \Gamma_1 \cdots \Gamma_k\;.
$$
For $k$ odd, the grading is odd and no such element exists.
\hfill $\diamond$
}
\end{example}

\begin{example}\label{Ex-GradedAlg1} 
{\rm
Let $\mathcal B$ be a graded C$^\ast$-algebra and consider the direct sum $\mathcal B \oplus \mathcal B$. It is a C$^\ast$-algebra when equipped with the addition $(a\oplus b)+ (a' \oplus b')=(a+a') \oplus (b+b')$, multiplication $(a\oplus b)(a' \oplus b')=aa' \oplus bb'$, as well as the norm $\|a \oplus b \| = \max \{\|a\|,\|b\|\}$. This algebra accepts the following two natural gradings: 
\begin{enumerate}[\rm (i)] 

\item The standard even grading $\gamma_{\Bb\oplus \Bb}(b \oplus b') = \gamma_\Bb(b) \oplus (-\gamma_\Bb(b'))$.

\item The standard odd grading $\gamma_{\Bb\oplus \Bb}(b \oplus b') = \gamma_\Bb(b') \oplus \gamma_\Bb(b)$. 

\end{enumerate}
The C$^\ast$-algebra $\mathcal B \oplus \mathcal B$ equipped with the odd grading is usually denoted as $\mathcal B_{(1)}$. 
}
\hfill $\diamond$
\end{example}

\begin{example}
\label{ex-CliffMatIso}
{\rm There exists the isomorphisms of graded algebras 
$$
\CM_k \cong \left \{
\begin{array}{ll}
\CM(2^\frac{k}{2}) \quad  & \mbox{for $k$ even}\;, \medskip \\
\CM(2^\frac{k-1}{2}) \oplus \CM(2^\frac{k-1}{2}) \quad & \mbox{for $k$ odd}\;,
\end{array}
\right .
$$
where for $k$ even $\CM(2^\frac{k}{2})$ is graded as in Example~\ref{ex-MatrixAlgGrading}, and for $k$ odd $\CM(2^\frac{k-1}{2})$ is ungraded, but the direct sum is equipped with the standard odd grading introduced above. In particular, $\CM_1 \cong \CM \oplus \CM=\CM_{(1)}$ and, if $\epsilon$ is the generator of $\CM_1$, then $\epsilon$ is sent to $1 \oplus (-1)$ by this graded isomorphism. 
}
\hfill $\diamond$
\end{example}

\begin{example}
\label{Ex-GradedAlg4}
{\rm
Let $\Bb=\KM$, the algebra of compact operators on a separable Hilbert space. A natural inner grading $\gamma_\KM$ can be introduced by using the isomorphism $\KM \cong\CM(2)\otimes\KM$ together with the inner grading on $\CM(2)$ introduced above. While the isomorphism used in the construction is not canonical, the graded algebras obtained this way are all graded-isomorphic. This example is important because most of the time the stabilization $\KM \otimes \Bb$ of an ungraded C$^\ast$-algebra will be viewed as a graded algebra with the grading $\gamma_\KM \otimes {\rm id}$. 
\hfill $\diamond$
}
\end{example}

\begin{definition}[\cite{BlackadarBook1998} pp. 116, and \cite{HR}, Appendix A]
\label{Def-GradedProd} 
The (algebraic) graded tensor product $\Aa \hat \otimes \Bb$ of two graded C$^\ast$-algebras $\Aa$ and $\Bb$ are defined by the following set of rules:
\begin{enumerate}[\rm (i)]

\item $(a \hat \otimes b)(a' \hat \otimes b') = (-1)^{{\rm deg}(b){\rm deg}(a')} aa' \hat \otimes bb'$ (multiplication)

\vspace{0.1cm}

\item $(a \hat \otimes b)^\ast = (-1)^{{\rm deg}(a){\rm deg}(b)} \, a^\ast \hat \otimes b^\ast$ (conjugation)

\end{enumerate}
while the rule for addition remains unchanged. The grading is uniquely defined by the rule ${\rm deg}(a \hat \otimes b) = {\rm deg}(a) + {\rm deg}(b)$.
\end{definition}

\begin{remark}{\rm Since one of the algebras in the graded products encountered in our work will always be finite dimensional, the topology of the graded tensor product will never be an issue. } 
\hfill $\diamond$ \end{remark}

\begin{remark}{\rm Note that if one of the algebras is trivially graded, then the graded tensor product reduces to the ordinary one. However, we will continue to use the hat in such situations just to indicate that the result is a graded algebra.} 
\hfill $\diamond$ \end{remark}

\begin{definition}[Graded Commutator] 
\label{def-GradedComm}
The natural version of the commutator when dealing with graded algebras is the graded commutator defined as the unique bilinear map satisfying $[x,y] = xy - (-1)^{{\rm deg}(x){\rm deg}(y)} yx$ for homogeneous $x, \, y$ elements of the algebra.
\end{definition}

\begin{example}[\cite{BlackadarBook1998} pp. 119]\label{Re-B1Rep} 
{\rm For a graded C$^\ast$-algebra $\Bb$ with even grading $\gamma_\Bb$, the graded algebra $\Bb_{(1)}$  introduced in Example~\ref{Ex-GradedAlg1} is isomorphic to $\mathcal B \hat \otimes \CM_1$.}
\hfill $\diamond$
\end{example}

\begin{example}
\label{ex-Cliffo}
{\rm Another useful example of graded product is the graded isomorphism $\varphi:\CM_1 \hat \otimes \, \CM_1 \to \CM(2)$, with even grading on $\CM(2)$. If $1$ and $\epsilon$ are the generators of $\CM_1$, then
$$
\varphi(1\hat\otimes 1)\,=\, 1
\;,
\qquad
\varphi(\epsilon \hat \otimes 1)\,=\, \sigma_2
\;,
\qquad
\varphi(1 \hat \otimes \epsilon )\,=\, \sigma_1
\;,
\qquad
\varphi(\epsilon \hat \otimes \epsilon )\,=\, - \imath \sigma_3
\;
$$
provide this isomorphism. Indeed, then the following relations are satisfied: $(1 \hat \otimes \epsilon)^2 = 1 \hat \otimes 1$, $(\epsilon \hat \otimes 1)^2 = 1 \hat \otimes 1$, and
$$
(1 \hat \otimes \epsilon ) ( \epsilon \hat \otimes 1 ) 
\;=\; 
- \epsilon \hat \otimes \epsilon
\;, 
\qquad 
(\epsilon \hat \otimes 1 ) ( 1 \hat \otimes \epsilon ) 
\;=\; 
\epsilon \hat \otimes \epsilon
\;, 
\qquad (\epsilon \hat \otimes \epsilon )^2\; =\;  -1 \hat \otimes 1
\;,
$$
as required by the rules in Definition~\ref{Def-GradedProd}. Combining $\varphi$ with Example~\ref{ex-CliffMatIso} it also follows that $\CM_1 \hat \otimes \, \CM_1 \cong \CM_2$. This is true more generally, namely $\CM_p \hat \otimes \, \CM_q = \CM_{p+q}$.}
\hfill $\diamond$
\end{example}

\subsection{Right Hilbert modules over C$^\ast$-algebras} 

The Hilbert C$^\ast$-modules are natural extensions of the ordinary Hilbert spaces. In our applications, they play the role of representation spaces for the algebras and, as it is well known \cite{CalderonPNAS1967yt}, many of the classical tools can be generalized to this new setting. As a result, one can resolve the non-commutative geometry of a vast new class of examples. Throughout this section, $\Bb$ stands for a graded C$^\ast$-algebra, unless it is explicitly specified otherwise.

\begin{definition}[Hilbert C$^\ast$-module, \cite{JensenThomsenBook1991} p.~1,and \cite{BlackadarBook1998,WeggeOlsenBook1993de}] 
A right Hilbert C$^\ast$-module over $\mathcal B$, or simply a Hilbert $\Bb$-module, is a complex linear space $E_\mathcal B$ which is also a right $\mathcal B$-module equipped with an inner product
$$
\langle \,\cdot \,,\, \cdot\, \rangle \;:\; E_\mathcal B \times E_\mathcal B \rightarrow \mathcal B
\;,
$$
which is linear in the second variable and satisfies the following relations for all $\psi, \psi' \in E_\Bb$ and $b \in \Bb$:

\begin{enumerate}[\rm (i)]

\item $\langle \psi , \psi' b \rangle = \langle \psi , \psi' \rangle b$,

\item $\langle \psi ,\psi' \rangle^\ast = \langle \psi' ,\psi \rangle$,

\item $\langle \psi , \psi \rangle \geq 0$,

\item $\psi \neq 0$ implies $\langle \psi , \psi \rangle \neq 0$.

\end{enumerate}
In addition, the complex linear space $E_\mathcal B$ must be complete in the norm induced by the inner product
$$
\|\psi \| 
\;=\; 
\| \langle \psi,\psi \rangle \|^\frac{1}{2}
\;, 
\qquad \psi \in E_\Bb\;.
$$
\end{definition}

\begin{remark} {\rm The above rules imply $\langle xb,y \rangle = b^\ast \langle x,y \rangle$. Hence the inner product is linear in the second argument and anti-linear in the first one, just as in the physics literature.}
\hfill $\diamond$
\end{remark}

\begin{example}\label{Ex-AlgebraHMungraded} 
{\rm 
A C$^\ast$-algebra $\mathcal B$ can be viewed as a Hilbert $\mathcal B$-module equipped with the inner product $\langle b ,b' \rangle = b^\ast b'$. 
}
\hfill $\diamond$
\end{example}

\begin{example}[Ungraded standard Hilbert module]\label{Ex-StandardHM1} 
{\rm Let $\mathcal B$ be a trivially graded C$^\ast$-algebra. The standard Hilbert $\mathcal B$-module is defined as the space
$$
\mathcal H_\mathcal B 
\;=\; 
\left\{
(b_n)_{n\in\NM}\ : \  b_n \in \mathcal B, \ \sum_{n=1}^N b_n^\ast b_n \ \mbox{converges in} \ \mathcal B \ {\rm as} \ N\rightarrow \infty
\right\}
\;,
$$
endowed with the obvious addition and the inner product
$$
\langle (b_n)_{n\in\NM}  ,  (b'_n)_{n\in\NM}  \rangle 
\;=\; \sum_{n \in \mathbb N} b_n^\ast b'_n \in \Bb
\;.
$$
If $\Bb$ is stable ({\rm i.e.} $\Bb \cong \KM \otimes \Bb$), then $\Hh_\Bb\cong\Bb$ (see Lemma~1.3.2 in \cite{JensenThomsenBook1991}).
}
\hfill $\diamond$
\end{example}

\begin{definition}[Isomorphism of Hilbert modules, \cite{JensenThomsenBook1991} p.~9, and \cite{WeggeOlsenBook1993de}]
\label{def-HMIsomorphism}
Two Hilbert $\mathcal B$-modules $E_\mathcal B$ and $F_\mathcal B$ are isomorphic whenever there is a complex linear bijection $\Psi: E_\mathcal B \rightarrow F_\mathcal B$ such that
$$
\langle \Psi(\psi), \Psi(\psi') \rangle 
\;=\; 
\langle \psi,\psi' \rangle.
$$
One writes $E_\mathcal B \approx F_\mathcal B$ and calls $\Psi$ an isomorphism. 
\end{definition}

\begin{remark} {\rm Isomorphism defines an equivalence relation on the set of Hilbert $\Bb$-modules. It also follows immediately from the properties of the inner product that any isomorphism between Hilbert $\mathcal B$-modules is automatically a $\mathcal B$-module map.}\hfill $\diamond$
\end{remark}

\begin{definition}[Graded Hilbert C$^\ast$-module, \cite{JensenThomsenBook1991} p.~26]
\label{Def-GradedHM}
 A Hilbert $\mathcal B$-module $E_\mathcal B$ over a gra\-ded C$^\ast$-algebra $\Bb$ is called graded if it is equipped with a complex linear bijection $S: E_\mathcal B \rightarrow E_\mathcal B$, called the grading operator, satisfying $S^2=1$ as well as 
$$
S (\psi b) 
\;=\; 
S(\psi) \gamma_\mathcal B(b)\;, 
\qquad \psi \in E_\mathcal B\;, \;\; b \in \mathcal B
\;,
$$
and
$$
\langle S \psi, S \psi' \rangle 
\;=\; 
\gamma_\mathcal B( \langle \psi,\psi' \rangle )\;, 
\qquad 
\psi, \ \psi' \in E_\mathcal B
\;.
$$
\end{definition}

\begin{definition}[Graded isomorphisms, \cite{JensenThomsenBook1991} p.~27]\label{Def-GradedHom}
If $E_\mathcal B$ and $F_\mathcal B$ are graded Hilbert $\mathcal B$-mod\-ules with grading operators $S_E$ and $S_F$, respectively, then a graded isomorphism is an isomorphisms $\Psi: E_\mathcal B \rightarrow F_\mathcal B$ of Hilbert $\mathcal B$-modules satisfying $\Psi \circ S_E = S_F \circ \Psi$.
\end{definition}

\begin{example}
\label{Ex-GradedHS} 
{\rm
The ordinary separable Hilbert space $\Hh$ can be regarded as a Hilbert $\CM$-module. If there exists a symmetry on $\Hh$, that is, a self-adjoint operator $J$  such that $J^2=1$, then $\Hh$ can be regarded as a graded Hilbert $\CM$-module with the grading operator $J$. Useful gradings are obtained when the invariant spectral subspaces of $S$ are infinite. Note that in these cases the spectral subspaces are isomorphic and in fact the graded Hilbert spaces are all isomorphic. The universal grading can be defined using the isomorphism $\Hh \cong \CM^2 \otimes \Hh$ and the grading $\sigma_3 \otimes 1$. Throughout, we will use the notation $\Sigma_3$ for this grading on $\Hh$, and when $\Hh$ is considered with the grading we use the notation $\widehat \Hh$. It is also useful to introduce a notation for the other two possible gradings $\Sigma_1$ and $\Sigma_2$ induced by the other generators $\sigma_1$ and $\sigma_2$ of $\CM(2)$.
}
\hfill $\diamond$
\end{example}

\begin{example}\label{Ex-AlgebraHM} 
{\rm 
A graded C$^\ast$-algebra $\mathcal B$ is a Hilbert $\Bb$-module (see Example~\ref{Ex-AlgebraHMungraded}), which becomes graded the grading operator $\gamma_\mathcal B$.
}
\hfill $\diamond$
\end{example}

\begin{example}\label{Ex-ProjectiveHM} 
{\rm Let $\mathcal B$ be a trivially graded C$^\ast$-algebra and $p\in \mathcal B$ a projection, that is, $p^2=p=p^\ast$. Let 
$$
p\mathcal B \;=\;\{pb \ : \ b \in \mathcal B\}
$$ 
be equipped with the right action $(pb)b' = pbb'$ and the inner product inherited from Example~\ref{Ex-AlgebraHM}, $\langle pb, pb' \rangle = b^\ast p b'$. Then $p\mathcal B$ is a trivially graded Hilbert $\mathcal B$-module. 

\vspace{.1cm}

Suppose next that $\mathcal B$ is graded and that $p$ is of graded degree $0$, that is, $\gamma_\Bb(p) = p$. Actually, a projection can never be of degree $1$ because this would lead to the contradiction $\gamma_\Bb(p) = \gamma_\Bb(p^2)=\gamma_\Bb(p)\gamma_\Bb(p)=(-p)(-p)=p$. Furthermore, every non-homogeneous projection in a matrix algebra over $\Bb$ can always be homotopically deformed in a projection of degree $0$, see Remark~\ref{rem-ProjDeg}. For a projection $p$ of degree $0$, $p\mathcal B$ is a graded Hilbert $\mathcal B$-module by the grading operator $S(pb) = \gamma_\mathcal B (pb) = p\big (\gamma_\mathcal B (b) \big )$. Furthermore, if $\Bb$ is evenly graded by the (self-adjoint and unitary) grading element $g\in\Mm(\Bb)$ ({\it i.e.} $\gamma_\Bb(b) = g^\ast b g$) and $p$ is such that $pg =p$, then $p\Bb$ is can be graded by $S(pb) = pbg$, which clearly satisfies $S(pbb') =S(pb) \gamma_\Bb (b')$. 
}
\hfill $\diamond$
\end{example}

\begin{example}[Graded standard Hilbert module]\label{Ex-GradedStandardHM1} 
{\rm If $\mathcal B$ is graded in the previous example, then $\mathcal H_\mathcal B$ accepts two natural gradings, one by 
$$
S(b_1,b_2, \ldots) 
\;=\; 
(\gamma_\mathcal B (b_1), \gamma_\mathcal B (b_2), \ldots)
$$
and one by $-S$. The graded standard module is defined as the graded Hilbert $\mathcal B$-module $\widehat{\mathcal H}_\mathcal B = \mathcal H_\mathcal B \oplus \mathcal H_\mathcal B$ equipped with the grading $\widehat{S}=S \oplus (-S)$. Note that, if $\Bb$ is ungraded, then simply $\widehat{S}=1\oplus(-1)$.}
\hfill $\diamond$
\end{example}

\begin{remark}{\rm Further examples of graded Hilbert C$^\ast$-modules are provided by the inner tensor product of existing modules but this operation requires the notion of operators on Hilbert C$^\ast$-modules, which is introduced next.}
\hfill $\diamond$
\end{remark}

\subsection{Classes of operators on Hilbert C$^\ast$-modules} 

Here again $\Bb$ stands for a graded C$^\ast$-algebra.

\subsubsection{Adjointable operators} 

The class of adjointable operators is the generalization of the bounded operators on Hilbert spaces. They will be denoted by large roman letters, while the elements in the underlying algebra $\Bb$ will still be denoted by small roman letters.

\begin{definition}[\cite{JensenThomsenBook1991} p.~4, and \cite{BlackadarBook1998,WeggeOlsenBook1993de}]  The space of adjointable operators $\mathbb B(E_\mathcal B)$ over a Hilbert $\mathcal B$-module consists of the maps $T: E_\mathcal B \rightarrow E_\mathcal B$ for which there exists a linear map $T^\ast : E_\mathcal B \rightarrow E_\mathcal B$ such that
$$
\langle T \psi,\psi' \rangle 
\;=\; 
\langle \psi, T^\ast \psi' \rangle
\;, 
\qquad \psi,\psi' \in E_\mathcal B
\;.
$$
\end{definition}

\begin{remark} 
{\rm 
An adjointable operator is automatically a bounded $\mathcal B$-module map, but the reversed implication is not always true. This is one of the major differences encountered when passing from Hilbert spaces to Hilbert C$^\ast$-modules. If the operator $T^\ast$ exists, it is unique and is called the adjoint of $T$. 
\hfill $\diamond$
}
\end{remark}

\begin{proposition}[\cite{JensenThomsenBook1991} p.~4, and \cite{BlackadarBook1998,WeggeOlsenBook1993de}] When endowed with the operator norm,
$$
\|T \| \;=\; \sup \big \{\|T \psi\| \ : \ \psi \in E_\mathcal B, \ \|\psi\| \leq 1 \big \}
\;,
$$
the space $\BM(E_\Bb)$ becomes a unital C$^\ast$-algebra.
\end{proposition}

\begin{definition}[Strict topology, \cite{JensenThomsenBook1991} p.~6] The strict topology of $\BM(E_\Bb)$ is defined as the locally convex topology  induced by the semi-norms
$$
\| T\|_\psi 
\;=\; 
\|T\psi\| \,+\, \|T^\ast \psi\| 
\;, 
\qquad \psi \in E_\Bb
\;.
$$
\end{definition}

\begin{remark}{\rm As in the case of Hilbert spaces, the strict topology is useful for expanding the operators as infinite sums of more elementary operators, such as the finite-rank operators introduced below. This will often occur later in our exposition.}
\hfill $\diamond$
\end{remark}

\begin{definition}[Induced grading, \cite{JensenThomsenBook1991} p.~6]
If $E_\mathcal B$ is a graded Hilbert $\mathcal B$-module with grading operator $S$, then $\mathbb B(E_\mathcal B)$ inherits the canonical grading
$$
\BM(E_\mathcal B) \,\ni\, T\;  \mapsto\; \mbox{\rm Ad}_S(T)\;=\;S T S^{-1}
\;.
$$
In such situations, $\mathbb B(E_\mathcal B)$ will always be viewed as a graded C$^\ast$-algebra with this particular grading.
\end{definition}

\begin{example}
\label{ex-AlgOperIso}
{\rm Let $\Bb$ be a graded unital C$^\ast$-algebra. As in Example~\ref{Ex-AlgebraHM}, $\Bb$ is also viewed as a graded Hilbert $\Bb$-module. Then $\BM(\Bb) \cong \Bb$, with the graded isomorphism provided by
$$
\BM(\Bb) \;\ni\; T \;\mapsto\;W(T) \;=\; T(1) \;\in \;\Bb\;.
$$
Indeed, one can easily verify that $W$ respects the multiplication,
$$
W(TT')\; =\; T\big ( T'(1) \big ) \;= \;T\big ( 1 \cdot T'(1) \big ) \;=\; T(1) \, T'(1)\;,
$$
and that $W^{-1}(b) = L_b$, $L_b(b') = bb'$, is the inverse of $W$. As for the grading, we have
$$
(W\circ \gamma_{\BM(\Bb)})(T)
\;=\;
W(\gamma_\Bb \circ T \circ \gamma_\Bb^{-1}) 
\;=\; 
\gamma_\Bb \big ( T(\gamma_\Bb^{-1}(1)) \big ) 
\;=\; \gamma_\Bb\big ( T(1) \big ) 
\;=\; 
(\gamma_\Bb \circ W)(T)
\;,
$$
or simply $W\circ \gamma_{\BM(\Bb)} = \gamma_\Bb \circ W$, hence the isomorphism is graded.
}
\hfill $\diamond$
\end{example}

\subsubsection{Finite rank operators} 

Again this is an immediate generalization of the corresponding concept for ordinary Hilbert spaces.

\begin{definition}[\cite{JensenThomsenBook1991} p.~5, and \cite{BlackadarBook1998,WeggeOlsenBook1993de}] The rank-one operators are defined as
\begin{equation}
\Theta_{\phi,\phi'} (\psi) \;=\; \phi \,\langle \phi', \psi \rangle\;, \qquad \psi,\, \phi, \, \phi' \in E_\mathcal B\;.
\end{equation} 
\end{definition}

\begin{proposition}[\cite{JensenThomsenBook1991,BlackadarBook1998,WeggeOlsenBook1993de}] The following are some elementary properties of the rank-one operators:
\begin{enumerate}[\rm (i)]

\item They are adjointable, with $\Theta_{\phi,\phi'}^\ast = \Theta_{\phi',\phi}$. 

\item Their product is again a finite-rank operator
$$
(\Theta_{\varphi,\varphi'} \, \Theta_{\phi,\phi'}) (\psi) 
\;=\; 
\Theta_{\varphi,\varphi'} \big ( \phi \,\langle \phi', \psi \rangle \big ) 
\;=\; 
\varphi \langle \varphi', \phi \rangle \langle \phi', \psi \rangle 
\; .
$$

\item For any adjointable operator $T$ over $E_\Bb$, 
$$
T \,\Theta_{\phi,\phi'} 
\;=\; 
\Theta_{T\phi,\phi'} \; , 
\qquad \Theta_{\phi,\phi'} (\psi)\, T 
\;=\; 
\Theta_{\phi,T^\ast\phi'} \; .
$$
Hence these products are also rank-one operators.
\end{enumerate}
\end{proposition}

\begin{proposition}[Ideal of finite-rank operators, \cite{JensenThomsenBook1991} p.~5, and \cite{BlackadarBook1998,WeggeOlsenBook1993de}] The algebraic span of the rank-one operators generates a double-sided ideal in $\BM(E_\Bb)$. This ideal is sometimes denoted by $\Theta(E_\Bb)$, a notation which will be adopted here.
\end{proposition}

\begin{example} {\rm According to Example~\ref{ex-AlgOperIso}, any operator from $\BM(\Bb)$ is of the form
$$
L_b(b') = bb' = b (1 \cdot b') = \Theta_{b, 1}(b').
$$
Hence, $\Theta(\Bb) = \BM(\Bb)$.
}
\hfill $\diamond$
\end{example}

\subsubsection{Compact operators}

\begin{definition}[\cite{JensenThomsenBook1991} p.~5, and \cite{BlackadarBook1998,WeggeOlsenBook1993de}] The closure of the set of finite rank operators on $E_\Bb$ in the norm topology of $\mathbb B(E_\mathcal B)$ is a closed double-sided essential ideal of $\BM(E_\Bb)$, called the ideal of compact operators and denoted by $\KM (E_\Bb)$. 
\end{definition}

\begin{remark}\label{Re-OpSM} {\rm For the standard ungraded Hilbert module $\mathcal H_\mathcal B$ from Example~\ref{Ex-StandardHM1}, the classes of operators introduced so far can be characterized as follows \cite{BlackadarBook1998,WeggeOlsenBook1993de}. First, there exists an isomorphism of C$^\ast$-algebras
$$
\mathbb K(\mathcal H_\mathcal B)\; \simeq\; \mathbb K \otimes \mathcal B,
$$
where $\mathbb K$ is the ungraded algebra of compact operators over a separable Hilbert space. Furthermore,
$$
\mathbb B(\mathcal H_\mathcal B)\; \simeq \;\mathcal M(\mathbb K \otimes \mathcal B)\;,
$$
where $\mathcal M$ indicates the multiplier algebra or the double centralizer \cite{BusbyTAMS1968jf}. From a standard property of the tensor products and their multiplier algebras (\cite{WeggeOlsenBook1993de}, Corollary~T.6.3), we have the chain of C$^\ast$-algebra inclusions
\begin{equation}\label{Eq-Chain}
\KM(\Hh_\Bb)
\;=\;
\mathbb K \otimes \mathcal B \;\subset\; 
\mathbb B \otimes \mathcal B \;\subset\;
\mathcal M(\mathbb K \otimes\mathcal B) 
\;=\; 
\BM(\Hh_\Bb)
\;,
\end{equation}
where $\mathbb B$ is the algebra of bounded operators over a separable Hilbert space and the tensor product  in $\BM \otimes \Bb$ is with the spatial C$^\ast$-norm, see \cite{WeggeOlsenBook1993de}.
\hfill $\diamond$
}
\end{remark}

\begin{example} 
\label{ex-StandardOps}
{\rm
Let us fix an orthonormal basis $(\psi_n)_{n\in \mathbb N}$ in the separable Hilbert space $\mathcal H$. Then any element of $\mathcal H_\mathcal B$ can be represented uniquely as  
$$
\psi\; =\; \sum_{n\in \mathbb N} \psi_n \otimes b_n
\;, 
\qquad b_n\in\Bb\; \mbox{ such that }\;\lim_{N \rightarrow \infty} \sum_{n=1}^N b_n^\ast b_n\; \in\; \mathcal B
\;.
$$
Furthermore, the elements from $\KM(\Hh_\Bb)$  can be uniquely represented as norm convergent sums in $\BM(\mathcal H_\mathcal B)$ of the form
$$
T \;=\; \sum_{n,m \in \mathbb N} \Theta_{\psi_n,\psi_m} \otimes b_{n,m}
\;,
\qquad
b_{n,m}\in\Bb
\;.
$$
The elements from $\BM(\Hh_\Bb)$ can be represented similarly, but the convergence of the sum needs to be considered in the strict topology and not the norm topology. Lastly, the rank-one operators take the form
$$
\Theta_{\phi',\phi} \;=\; \sum_{n,m \in \mathbb N} \Theta_{\psi_n,\psi_m} \otimes b'_n b_m^\ast
\;,
$$
if $\phi = \sum_{n \in \mathbb N} \psi_n \otimes b_n$ and $\phi' = \sum_{n \in \mathbb N} \psi_n \otimes b'_n$.
\hfill $\diamond$
}
\end{example}

\subsubsection{Generalized Fredholm operators} 

Given a Hilbert $\mathcal B$-module $E_\Bb$, let us consider the short exact sequence of C$^\ast$-algebras,
\begin{equation}
\begin{diagram}
0 & \rTo & \mathbb K(E_\Bb) & \rTo{\mathfrak i}  & \mathbb B(E_\Bb) & \rTo{ \mathrm{ev} } & \mathbb Q(E_\Bb) & \rTo & 0 \; ,
\end{diagram}
\label{eq-Calkin}
\end{equation}
where  $\mathbb Q(E_\Bb)=\mathbb B(E_\Bb)/ \mathbb K(E_\Bb)$ is the corona or the generalized Calkin algebra and ${\rm ev}$ is the quotient map. The following is a minor modification of 17.5 in \cite{BlackadarBook1998} and \cite{MingoTAMS1987fg}.

\begin{definition} The class of generalized Fredholm operators over $E_\Bb$ is defined as
$$
\mathbb F(E_\Bb) 
\;=\; 
\big \{F \in \mathbb B(E_\Bb) \ : \ {\rm ev}(F) \ \mbox{invertible in} \ \mathbb Q(E_\Bb) \big \}
\;.
$$
This class can be further sub-divided in the sub-classes of unitary and self-adjoint generalized Fredholm operators, as well as the generalized Fredholm symmetries:
\begin{align*}
& \mathbb F_\gu(E_\Bb)\; = \;\big \{T \in \mathbb F(E_\Bb) \ : \ {\rm ev}(TT^*)={\rm ev}(T^*T)=1 \big \}
\;,
\\
& \mathbb F_\sa(E_\Bb)\; = \;\big \{T \in \mathbb F(E_\Bb) \ : \ {\rm ev}(T)^*={\rm ev}(T) \big \}
\;,
\\
& \mathbb F_\sym(E_\Bb)\; = \;\mathbb F_\gu(E_\Bb)\,\cap\,\mathbb F_\sa(E_\Bb)
\;.
\end{align*}
In case $E_\Bb$ has a grading $S_E$, each of these classes has even and odd sectors, 
$$
\mathbb F_{\pm}(E_\Bb)\; = \;\big \{T \in \mathbb F(E_\Bb) \ : \ S_E T S_E=\pm T \big \}
\;.
$$
The odd Fredholm symmetries $\mathbb F_{\sym,-}(E_\Bb)=\mathbb F_\sym(E_\Bb)\cap \mathbb F_{-}(E_\Bb)$ are also called Fredholm supersymmetries.
\end{definition} 

\begin{remark} {\rm The product of two generalized Fredholm operators is also a Fredholm operator. As such, $\mathbb F(E_\Bb)$  forms a group. Furthermore, when the quotient w.r.t. homotopy equivalence is taken, the group becomes abelian. Later it will become clear that also $\mathbb F_\sa(E_\Bb)$ has a group structure when $E_\Bb$ is the standard Hilbert module. 
} 
\hfill $\diamond$  
\end{remark}

\subsection{Kasparov's stabilization theorem} 

This is the central basic result in Kasparov's $K$-theory. It is mentioned here explicitly because it will be used often, primarily to characterize the so-called countably generated Hilbert C$^\ast$-modules (see Corollary~\ref{Co-FiniteHM}). The latter are the Hilbert C$^\ast$-modules which enter in the definition of the $KK$-groups, and hence will play an important role in the following. Throughout, $\Bb$ stands again for a graded C$^\ast$-algebra.

\begin{definition}[Countably generated Hilbert C$^\ast$-modules, \cite{JensenThomsenBook1991} pp.~11, 26, and \cite{BlackadarBook1998,WeggeOlsenBook1993de}] \label{Def-CountableHM}
A graded Hil\-bert $E_\mathcal B$-module is said to be countably generated whenever there is a countable subset $\{\psi_n\}_{n \in \NM}$ of $E_\mathcal B$ such that the linear span of the set $\{\psi_n b \,: \, n \in \NM, \, b \in \mathcal B\}$ is dense in $E_\mathcal B$. 
\end{definition}

\begin{theorem}[Kasparov's stabilization theorem, \cite{KasparovMUI1981ui}, and \cite{BlackadarBook1998,HR}] 
\label{theo-Stabilize}
If $E_\mathcal B$ is a countably generated graded Hilbert $\mathcal B$-module, then there always exists an isomorphism of graded Hilbert $\mathcal B$-modules
$$
E_\mathcal B \oplus \widehat{\mathcal H}_\mathcal B 
\;\simeq\; 
\widehat{\mathcal H}_\mathcal B
\;.
$$
Recall that the grading operator for the direct sum is given by $S = S_E \oplus \widehat S$.
\end{theorem}

\begin{corollary}
\label{Co-FiniteHM} 
Every countably generated graded Hilbert $\mathcal B$-module is isomorphic to $P\,\widehat{\mathcal H}_\mathcal B$, where $P$ is a  projection from $\mathbb B(\widehat{\mathcal H}_\mathcal B)$ commuting with the grading $\widehat S$.
\end{corollary}

\proof Let $R \in \mathbb B(E_\mathcal B \oplus \widehat{\mathcal H}_\mathcal B)$ be defined by $R(\psi \oplus \phi) = \psi \oplus 0$. Let $S=S_E \oplus \widehat S$ be the grading of $E_\mathcal B \oplus \widehat{\mathcal H}_\mathcal B$. Then $S R = R S$ and $R^2 =R$, as well as 
$$
\langle \psi \oplus \phi, R(\psi' \oplus \phi') \rangle 
\;=\; 
\langle R(\psi \oplus \phi), \psi' \oplus \phi' \rangle
\;.
$$
Hence, $R$ is a projection which commutes with $S$. Let $\Phi:E_\mathcal B \oplus \widehat{\mathcal H}_\mathcal B \to\widehat{\mathcal H}_\mathcal B$ be the graded isomorphism in the Kasparov's stabilization theorem. Then
$$
(\Phi R \Phi^{-1}) (\widehat{\mathcal H}_\mathcal B ) 
\;=\; 
(\Phi R)(E_\mathcal B \oplus \widehat{\mathcal H}_\mathcal B) 
\;=\; 
\Phi(E_\mathcal B \oplus 0)
\;\simeq \; 
E_\mathcal B
\;.
$$
Hence $P = \Phi R \Phi^{-1}$ gives the desired projection from $\mathbb B(\widehat{\mathcal H}_\mathcal B)$. Since $\Phi \circ S = \widehat S \circ \Phi$, we can also conclude that $\widehat S P = P \widehat S$.
\qed 

\subsection{Products of Hilbert C$^\ast$-modules}
\label{Sec-ProdHM}

\subsubsection{The internal tensor product}
\label{Sec-IntProdHM}

The internal tensor product is at the heart of Kasparov's cycles, $KK$-groups and Kasparov's product. We will review here the basic concept and provide the relevant examples.

\begin{definition}[Ungraded internal tensor product, \cite{JensenThomsenBook1991} p. 19, and \cite{BlackadarBook1998} 13.5]\label{Def-UnGradedIntProd} Let be given a  Hilbert $\mathcal A$-module $E_\mathcal A$  and a Hilbert $\mathcal B$-module $E_\mathcal B$  which is, moreover, equipped with a $\ast$-representation $\pi:\mathcal A\to \BM(E_\Bb)$ (a left action of $\Aa$ as linear operators on $E_\Bb$) which makes $E_\Bb$ into a $(\Aa,\Bb)$-bimodule. Consider the algebraic tensor product $E_\mathcal A \odot E_\mathcal B$ equipped with the $\Bb$-valued bilinear form
\begin{equation}\label{Eq-ScalarProd1}
\big \langle \psi_\mathcal A \odot \psi_\mathcal B, \psi'_\mathcal A \odot \psi'_\mathcal B \big \rangle 
\;=\; 
\big \langle \psi_\mathcal B, \pi(\langle \psi_\mathcal A,\psi'_\mathcal A \rangle) \psi'_\mathcal B \big \rangle
\;.
\end{equation}
Then the internal tensor product $E_\mathcal A \otimes_\mathcal A E_\mathcal B$ is defined as the Hilbert $\mathcal B$-module  obtained by the completion w.r.t. \eqref{Eq-ScalarProd1} of the quotient of $E_\mathcal A \odot E_\mathcal B$ by the sub-module of elements $\psi \in  E_\mathcal A \odot E_\mathcal B$ for which $\langle \psi,\psi \rangle =0$.
\end{definition}

\begin{remark} 
{\rm The submodule which is divided out is generated by the elements of the form 
$$
\psi_\mathcal A \, a \otimes \psi_\mathcal B \,-\, \psi_\mathcal A \otimes \pi(a) \psi_\mathcal B
\;.
$$
In other words, the elements $\psi_\mathcal A \, a \otimes \psi_\mathcal B$ and $\psi_\mathcal A \otimes \pi(a) \psi_\mathcal B$ are identical in the internal product $E_\mathcal A \otimes_\mathcal A E_\mathcal B$. The internal tensor product is associative, namely
$$
(E_\mathcal A \otimes_\mathcal A E_\mathcal B) \otimes_\mathcal B E_\mathcal C
\;=\;
E_\mathcal A \otimes_\mathcal A (E_\mathcal B \otimes_\mathcal B E_\mathcal C)
\;.
$$
This can be checked from the definition.
\hfill $\diamond$
}
\end{remark}

\begin{example}[Ungraded standard Hilbert module]\label{Ex-StandardHM2} 
{\rm
For a trivially graded C$^\ast$-algebra $\Bb$, the standard ungraded Hilbert module $\Hh_\Bb$ from Example~\ref{Ex-StandardHM1} accepts the following equivalent presentation $\Hh_\Bb = \Hh \otimes_\CM \Bb$, with
$$
\Psi\;:\;\mathcal H \otimes_\CM \mathcal B\;\to\;\Hh_\Bb \;, 
\qquad \Psi(\psi_n\otimes_\CM b)\;=\; (\ldots,0,b,0,\ldots)
\;,
$$
providing the isomorphism between the presentations. Above, $(\psi_n)_{n\in \NM}$ is an orthogonal basis of $\Hh$ and, on the r.h.s., $b$ sits precisely at the $n$-th entry.}
\hfill $\diamond$
\end{example}

\begin{proposition}
\label{prop-UnGradedIntProd} 
Let $E_\mathcal A = p\mathcal A$ be the rank-one Hilbert $\mathcal A$-module from Example~\ref{Ex-ProjectiveHM} associated to an ungraded C$^*$-algebra $\Aa$ and a projection $p \in \mathcal A$, and let $\pi$ be a representation of $\Aa$ in $\BM(E_\mathcal B)$. Then
$$
\Psi\,:\, p \mathcal A \otimes_\mathcal A E_\Bb \;\rightarrow\; \pi(p) E_\mathcal B\;, 
\qquad 
\Psi(pa \otimes_\mathcal A \psi )\;=\; \pi(pa)\psi
\;,$$
is an isomorphism of Hilbert $\mathcal B$-modules. Here, $\pi(p) E_\mathcal B$ is the Hilbert $\mathcal B$-module generated by the range of the projection $\pi(p)$.
\end{proposition}

\proof The map is well-defined since $\Psi$ is identically zero on the submodule which is quoted out, generated by $\psi_\mathcal A a \otimes \psi_\mathcal B - \psi_\mathcal A \otimes \pi(a) \psi_\mathcal B$, $\psi_\mathcal A \in E_\mathcal A$ and $\psi_\mathcal B \in E_\mathcal B$. It is bijective because $\Psi(p \otimes_\mathcal A \psi) = \pi(p) \psi$, $\psi \in E_\mathcal B$, span the entire $\pi(p) E_\mathcal B$ and $\Psi(p a \otimes_\mathcal A \psi) = 0$ implies $\pi(pa)\psi=0$, hence 
$$
pa \otimes_\mathcal A \psi \;=\; p(pa) \otimes_\mathcal A \psi 
\;=\; p \otimes_\mathcal A \pi(pa) \psi 
\;=\;0
\;,
$$
so that the map is injective. Lastly
$$
\langle \Psi(pa \otimes_\mathcal A \psi) , \Psi(pa' \otimes_\mathcal A \psi' )\rangle 
\;=\; 
\langle \pi(pa)\psi,\pi(pa') \psi' \rangle 
\;=\; 
\langle \psi, \pi(a^\ast p a)\psi' \rangle
\;.
$$
Hence the map preserves the inner product. \qed

\begin{example}\label{Ex-Tensorizing} 
{\rm
Assume the setting of Definition~\ref{Def-UnGradedIntProd}. To any given $\psi_\Aa\in E_\Aa$, one can associate the operator
$$
T_{\psi_\Aa}:E_\Bb\to E_\mathcal A \otimes_\mathcal A E_\mathcal B \; , 
\qquad 
T_{\psi_\Aa}\psi'_\Bb\;=\;\psi_\Aa\otimes_\mathcal A \psi'_\Bb
\;.
$$
This operator is adjointable, with adjoint
$$
T_{\psi_\Aa}^*:E_\mathcal A \otimes_\mathcal A E_\Bb\to  E_\mathcal B\;, 
\qquad 
T_{\psi_\Aa}^*\psi'_\Aa\otimes_\mathcal A \psi'_\Bb\;=\;\pi(\langle\psi'_\Aa, \psi_\Aa\rangle)^*\psi'_\Bb
\;,
$$
as can be readily checked. It will play a central role in the definition of the Kasparov product.
}
\hfill $\diamond$
\end{example}

\begin{definition}[Graded internal product, \cite{JensenThomsenBook1991} p. 65, and \cite{BlackadarBook1998} 14.4]\label{Def-GradedIntProd} 
If both the C$^\ast$-algebras and the Hilbert modules  in Definition~\ref{Def-UnGradedIntProd} are graded and the representation $\pi$ is a graded map from $\mathcal A$ to $\mathbb B(E_\Bb)$, then the internal product $E_\mathcal A \otimes_\mathcal A E_\mathcal B$ is a graded Hilbert $\mathcal B$-module w.r.t. the grading operator
$$
S (\psi_\mathcal A \otimes_\mathcal A \psi_\mathcal B) 
\;=\; 
S_{E_\mathcal A}(\psi_\mathcal A) \otimes_\mathcal A S_{E_\mathcal B}(\psi_\mathcal B)
\;.
$$
We then use the notation $E_\mathcal A \, \hat \otimes_\mathcal A \, E_\mathcal B$ to indicate the presence of this extra structure, even though  the internal product has not been modified.
\end{definition}

\begin{example}[Graded standard Hilbert module]\label{Ex-GradedStandardHM2} 
{\rm
For a graded C$^\ast$-algebra $\Bb$, the standard graded Hilbert module $\widehat \Hh_\Bb$ of Example~\ref{Ex-GradedStandardHM1} accepts the following equivalent presentation $\widehat \Hh_\Bb = \widehat \Hh \hat \otimes_\mathbb C \mathcal B$, with 
$$
\Psi\;:\;\widehat \Hh \hat \otimes_\CM \Bb \;\to\; \widehat \Hh_\Bb \; , 
\qquad \Psi \big ( (\psi^+_n \oplus \psi^-_m) \hat \otimes_\CM b \big )\;=\; ( \ldots, 0,b,0,\ldots) \oplus (\ldots , 0,b,0, \ldots)
\;,
$$
providing the isomorphism between the presentations. Above, $(\psi^\pm_n)_{n \in \NM}$ are orthonormal bases for the graded sectors $\Hh_\pm$ of $\widehat \Hh$ and on the r.h.s. $b$ sits at the $n$-th and $m$-th entries, respectively. Indeed, the grading of  $\widehat{\Hh} \, \hat \otimes_\mathbb C \, \mathcal B$ is $S=(1\oplus(-1)) \otimes_\mathbb C \gamma_\mathcal B$ and we can check that
$$
\Psi \big ( S ((\psi^+_n \oplus \psi^-_m) \hat \otimes_\CM b) \big ) \,=\, \Psi \big ( (\psi^+_n \oplus (-\psi^-_m)) \hat \otimes_\CM \gamma_\Bb(b) \big ) \,=\, ( \ldots, 0,\gamma_\Bb(b),0,\ldots) \oplus (\ldots , 0,-\gamma_\Bb(b),0, \ldots) \,,
$$
hence $\Psi$ respects the gradings ({\it cf.} Definition~\ref{Ex-GradedStandardHM1}).}
\hfill $\diamond$
\end{example}

\begin{remark} {\rm In general, there is no canonical characterization of the algebra $\BM(E_\Aa \hat \otimes_\Aa E_\Bb)$ of adjointable operators over the internal product in terms of the algebras $\BM(E_\Aa)$ and $\BM(E_\Bb)$. However, when $[\pi(a),T]=0$ (graded commutator) for all $a\in \Aa$ and $T \in \BM(E_\Bb)$, such a characterization exists and will prove extremely useful in several instances, notably when $\Aa = \CM$.}
\hfill $\diamond$
\end{remark}

\begin{proposition} 
\label{prop-internalProdOp}
Let $E_\Aa$ and $E_\Bb$ be graded Hilbert C$^\ast$-modules over the graded C$^\ast$-algebras $\Aa$ and $\Bb$, respectively, and let $\pi:\Aa \rightarrow \BM(E_\Bb)$ be a graded homomorphism such that $[\pi(a),T]=0$ for all $a\in \Aa$ and $T \in \BM(E_\Bb)$. Then there exists the homomorphism of graded algebras
$$
\BM(E_\Aa)\, \hat \otimes \,\BM(E_\Bb)\; \rightarrow \;\BM(E_\Aa \hat \otimes_\Aa E_\Bb)\;, 
\qquad
T_1 \hat \otimes T_2\;\mapsto\; T_1 \hat \otimes_\Aa T_2\;,
$$
where the latter operator is defined by
$$
(T_1 \hat \otimes_\Aa T_2) \, \psi_\Aa \hat \otimes_\Aa \psi_\Bb 
\;=\; 
(-1)^{{\rm deg}(T_2){\rm deg}(\psi_\Aa)} \,  T_1 \psi_\Aa \hat \otimes_\Aa T_2 \psi_\Bb
\;.
$$
\end{proposition}

\proof The above definitions are very similar to the ones for the external graded tensor product (see \cite{BlackadarBook1998} p.~139 and Section~\ref{Sec-ExtProdHM} below). It is instructive to see why for the internal graded tensor product the range of $\pi$ must be contained in the center of $\BM(E_\Bb)$. For this purpose, let us consider two equivalent elements $\psi_\Aa a \otimes_\Aa \psi_\Bb \sim \psi_\Aa \otimes_\Aa \pi(a)\psi_\Bb$. Then
\begin{align*}
(T_1 \hat \otimes_\Aa T_2) \psi_\Aa a \otimes_\Aa \psi_\Bb &  
\;=\; (-1)^{{\rm deg}(T_2){\rm deg}(\psi_\Aa a)} \, T_1 \psi_\Aa a \hat \otimes_\Aa T_2 \psi_\Bb \\
& \;= \;(-1)^{{\rm deg}(T_2)({\rm deg}(\psi_\Aa) + {\rm deg}(a))} \, T_1 \psi_\Aa \hat \otimes_\Aa \pi(a) T_2 \psi_\Bb
\;,
\end{align*}
while
\begin{align*}
(T_1 \hat \otimes_\Aa T_2) \psi_\Aa \otimes_\Aa \pi(a)\psi_\Bb
\; =\; (-1)^{{\rm deg}(T_2){\rm deg}(\psi_\Aa)} \, T_1 \psi_\Aa \hat \otimes_\Aa T_2 \pi(a) \psi_\Bb
\;.
\end{align*}
The difference of the two is
$$
(-1)^{{\rm deg}(T_2){\rm deg}(\psi_\Aa)}  \, T_1 \psi_\Aa \hat \otimes_\Aa \big (T_2 \pi(a) - (-1)^{{\rm deg}(T_2){\rm deg}(\pi(a))}\pi(a) T_2 \big ) \psi_\Bb
$$
and now one sees that $T_1 \hat \otimes_\Aa T_2$ is ill-defined unless $[\pi(a),T_2]=0$ for all $a \in \Aa$. Next, let us verify that the homomorphism respects multiplications. According to Definition~\ref{Def-GradedProd} applied to $\BM(E_\Aa)\, \hat \otimes \,\BM(E_\Bb)$, 
\begin{align*}
(T_1 \hat \otimes T_2)(T'_1 \hat \otimes T'_2) 
\;=\; (-1)^{{\rm deg}(T_2) {\rm deg}(T_1')} \, T_1 T_1' \hat \otimes T_2 T_2'
\end{align*}
and 
$$
(T_1 T_1' \hat \otimes_\Aa T_2 T_2') \, \psi_\Aa \hat \otimes_\Aa \psi_\Bb 
\;=\; 
(-1)^{({\rm deg}(T'_2)+{\rm deg}(T_2)){\rm deg}(\psi_\Aa)} \, T_1 T_1'\psi_\Aa \hat \otimes_\Aa T_2 T_2'\psi_\Bb \; ,
$$
while
\begin{align*}
 \Big ( (T_1 \hat \otimes_\Aa T_2)(T'_1 \hat \otimes_\Aa T'_2)\Big )   & \, \psi_\Aa \hat \otimes_\Aa \psi_\Bb 
\;=\; (-1)^{{\rm deg}(T'_2){\rm deg}(\psi_\Aa)} (T_1 \hat \otimes_\Aa T_2) \,  T'_1 \psi_\Aa \hat \otimes_\Aa T'_2 \psi_\Bb \\
& \;=\; (-1)^{{\rm deg}(T'_2){\rm deg}(\psi_\Aa)} (-1)^{{\rm deg}(T_2){\rm deg}(T_1'\psi_\Aa)} \, T_1 \psi_\Aa \hat \otimes_\Aa T_2 \psi_\Bb \\
& \;=\; (-1)^{({\rm deg}(T'_2)+{\rm deg}(T_2)){\rm deg}(\psi_\Aa)} (-1)^{{\rm deg}(T_2){\rm deg}(T_1')} \, T_1T_1' \psi_\Aa \hat \otimes_\Aa T_2 T_2' \psi_\Bb \; ,
\end{align*}
hence the answer is affirmative. 
\qed

\begin{remark}
\label{rem-picommute}
{\rm The calculation in the proof of Proposition~\ref{prop-internalProdOp} shows that there are other instances for which $T_1 \hat \otimes_\Aa T_2$ provides a well-defined operator. One example is that $T_2$ is of degree $0$ and commutes with $\pi$. A special case is $T_2=1$ and then the action is simply
$$
(T_1 \hat \otimes_\Aa 1) \, \psi_\Aa \hat \otimes_\Aa \psi_\Bb 
\;=\; T_1 \psi_\Aa \hat \otimes_\Aa \psi_\Bb
\;.
$$
This is of particular importance for the construction of the Kasparov product, see Section~\ref{sec-KKProd}.
}
\hfill $\diamond$
\end{remark}

\begin{example}
\label{ex-GradedEquiv1}
{\rm Let $\Aa$ and $\Bb$ be unital graded C$^\ast$-algebras, $E_\Bb$ a graded Hilbert $\Bb$-module and $\pi:\Aa \rightarrow \BM(E_\Bb)$ a graded non-degenerate homomorphism. Then $\Aa \hat \otimes_\Aa E_\Bb \cong E_\Bb$, with
$$
\Psi_1(a \hat \otimes_\Aa \psi_\Bb) 
\;=\; 
\pi(a) \psi_\Bb\;, 
\qquad 
\Psi_1^{-1}(\psi_\Bb) \;=\; 1 \hat \otimes_\Aa \psi_\Bb
$$
providing the isomorphism. One can check explicitly that 
$$
\Psi_1\big ( \gamma_\Aa(a) \hat \otimes_\Aa S_E \psi_\Bb \big ) 
\;=\; 
\pi\big (\gamma_\Aa(a) \big )  S_E \psi_\Bb 
\;=\; 
S_E \big ( \pi(a) \psi_\Bb \big )
\;,
$$
hence $\Psi_1$ is indeed a graded isomorphism.
}
\hfill $\diamond$
\end{example}

\begin{example} 
\label{ex-GradedEquiv2}
{\rm The above example can be pushed one step further to $\widehat \Hh_\Aa \hat \otimes_\Aa E_\Bb \cong \widehat\Hh\, \hat \otimes_\CM E_\Bb$. Indeed,
$$
\widehat \Hh_\Aa \hat \otimes_\Aa E_\Bb 
\;\cong\; 
(\widehat \Hh \,\hat \otimes_\CM \Aa) \hat \otimes_\Aa E_\Bb 
\;\cong \;
\widehat \Hh \,\hat \otimes_\CM (\Aa \hat \otimes_\Aa E_\Bb) 
\;\cong\; 
\widehat \Hh \,\hat \otimes_\CM E_\Bb\;,
$$
where the associativity of the internal tensor product was put to work. The isomorphism $\Psi_2 : \widehat \Hh_\Aa \hat \otimes_\Aa E_\Bb \to   \widehat \Hh \,\hat \otimes_\CM E_\Bb$ is given explicitly by
$$
\Psi_2 \big ( (\psi \hat \otimes_\CM a) \hat \otimes_\Aa \psi_\Bb \big ) \;=\; \psi \hat \otimes_\CM \pi(a) \psi_\Bb
\;, 
\qquad 
\Psi_2^{-1}(\psi \hat \otimes_\CM \psi_\Bb) \;=\; (\psi \hat \otimes_\CM 1 ) \hat \otimes_\Aa \psi_\Bb
\;. 
$$
If, moreover, $E_\Bb=\widehat{\Hh}_\Bb$, then
$$
\widehat \Hh_\Aa \hat \otimes_\Aa \widehat{\Hh}_\Bb 
\;\cong\; 
\widehat \Hh \,\hat \otimes_\CM \widehat{\Hh}_\Bb\;
\;\cong\; 
\widehat \Hh \,\hat \otimes_\CM \widehat{\Hh} \,\hat \otimes_\CM \Bb\;
\;\cong\; 
\widehat \Hh  \,\hat \otimes_\CM \Bb\;
\;\cong\; 
\widehat \Hh_\Bb\;.
$$
Equivalences of these types are often used below, and next follows yet another such example.}
\hfill $\diamond$
\end{example} 

\begin{example}
\label{ex-GradedEquiv3}
{\rm Consider the conditions from Example~\ref{ex-GradedEquiv2} and let $\widetilde \Pi : \BM(\widehat \Hh_\Aa \hat \otimes_\Aa E_\Bb) \rightarrow \BM(\widehat \Hh \hat \otimes_\CM E_\Bb)$ be the isomorphism induced by $\Psi_2$. Then there is a canonical homomorphism \cite[p.~138]{BlackadarBook1998}
$$
\Pi\, :\, \BM(\widehat \Hh_\Aa) 
\;\rightarrow\; \BM(\widehat \Hh \hat \otimes_\CM E_\Bb)
\;, 
\qquad 
\Pi(T) 
\;=\; 
\widetilde \Pi(T \hat \otimes_\Aa 1)
\;=\;
\Psi_2\,(T \hat \otimes_\Aa 1) \,\Psi_2^{-1}
\;.
$$
Let now $E_\Aa$ be a countably generated graded Hilbert $\Aa$-module. Then $E_\Aa \cong P \widehat \Hh_\Aa$, with $P\in \BM(\widehat \Hh_\Aa)$ the projection of graded degree 0 from Corollary~\ref{Co-FiniteHM}, and
$$
E_\Aa \hat \otimes_\Aa E_\Bb 
\;\cong\; 
P\widehat \Hh_\Aa \hat \otimes_\Aa E_\Bb 
\;=\; 
(P \hat \otimes_\Aa 1) \big (\widehat \Hh_\Aa \hat \otimes_\Aa E_\Bb \big ) 
\;\cong\; 
\widetilde \Pi(P \hat \otimes_\Aa 1) \big (\widehat \Hh \hat \otimes_\CM E_\Bb \big ) 
\;=\; 
\Pi(P) \big (\widehat \Hh \hat \otimes_\CM E_\Bb \big )
\;,
$$
or in short, $E_\Aa \hat \otimes_\Aa E_\Bb \cong \Pi(P) \big (\widehat \Hh \hat \otimes_\CM E_\Bb \big )$, with
$$
\Psi_3 \big ( P\psi_\Aa \hat \otimes_\Aa \psi_\Bb ) 
\;=\; 
\Psi_2 \big ( P\psi_\Aa \hat \otimes_\Aa \psi_\Bb ) 
\;=\; 
\Pi(P) \Psi_2(\psi_\Aa \hat \otimes_\Aa \psi_B)
$$
providing the isomorphism.
}
\hfill $\diamond$
\end{example}

\subsubsection{The external tensor product}
\label{Sec-ExtProdHM}

The external tensor product of Hilbert C$^\ast$-modules plays a role in the representation of higher Kasparov groups, where one needs to deal with modules over external tensor products such as $\Bb \hat \otimes \CM_k$ (see Section~\ref{sec-KKGroups}).  As we have seen above, the internal product $E_\Aa \hat \otimes_\Aa E_\Bb$ is just a Hilbert $\Bb$-module. In contradistinction, the result of the external tensor product $E_\Aa \hat \otimes E_\Bb$ will be a module over $\Aa \hat \otimes \Bb$.

\begin{definition}[Ungraded external tensor product, \cite{BlackadarBook1998} p.~111] Let $\Aa$ and $\Bb$ be ungraded unital C$^\ast$-algebras and $E_\Aa$, $E_\Bb$ be two ungraded Hilbert C$^\ast$-modules. The external tensor product $E_\Aa \otimes E_\Bb$ is the completion of the algebraic tensor product $E_\mathcal A \odot E_\mathcal B$  w.r.t. the inner product
$$
\langle \psi_\Aa \odot \psi_\Bb,  \psi'_\Aa \odot \psi'_\Bb \rangle 
\;=\; \langle \psi_\Aa,\psi'_\Aa \rangle \otimes \langle \psi_\Bb, \psi'_\Bb \rangle
\;.
$$
The result is a Hilbert $\Aa \otimes \Bb$-module, with
$$
(\psi_\Aa \otimes \psi_\Bb)(a \otimes b) 
\;=\; \psi_\Aa a \otimes \psi_\Bb b
$$
providing the right multiplication.
\end{definition}

\begin{remark} {\rm In general, the completion of the algebraic tensor product w.r.t. above scalar product is done in two steps ({\it cf.} \cite{JensenThomsenBook1991} p.~21), the complications being related with the topology of $\Aa \otimes \Bb$.
}
\hfill $\diamond$
\end{remark}

\begin{definition}[Graded external tensor product, \cite{BlackadarBook1998} p.~118] Let $\Aa$ and $\Bb$ be graded unital C$^\ast$-algebras and $E_\Aa$, $E_\Bb$ be two graded Hilbert C$^\ast$-modules. Then the algebraic tensor product $E_\mathcal A \odot E_\mathcal B$ can be transformed in a graded Hilbert $\Aa \hat \otimes \Bb$-module with multiplication
$$
(\psi_\Aa \hat \otimes \psi_\Bb)(a \hat \otimes b) 
\;=\; 
(-1)^{{\rm deg}(\psi_\Bb){\rm deg}(a)} \, \psi_\Aa a\, \hat \otimes \psi_\Bb b
\;,
$$
inner product
$$
\langle \psi_\Aa \hat \otimes \psi_\Bb, \psi'_\Aa \hat \otimes \psi'_\Bb \rangle 
\;=\; 
(-1)^{{\rm deg}(\psi_\Bb)({\rm deg}(\psi_\Aa) + {\rm deg}(\psi'_\Aa))}\; 
\langle \psi_\Aa,\psi'_\Aa \rangle \hat \otimes \langle \psi_\Bb, \psi'_\Bb \rangle 
\;\in \;\Aa \hat \otimes \Bb
$$
and grading given by ${\rm deg}(\psi_\Aa \hat \otimes \psi_\Bb) = {\rm deg}(\psi_\Aa) + {\rm deg}(\psi_\Bb)$.
\end{definition}

\begin{remark} {\rm Let us point out  that  the sign $(-1)^{{\rm deg}(\psi_\Aa){\rm deg}(\psi_\Bb)}$ results from the fact that $\psi_\Aa$ is on the right of $\psi_\Bb$ on the l.h.s. (due to the adjoint in the inner product), while on the r.h.s. it is to the left.
}
\hfill $\diamond$
\end{remark}

\begin{remark}[\cite{BlackadarBook1998} p.~118]
\label{re-ExtProdOp}
{\rm There exists a natural imbedding of $\BM(E_\Aa) \hat \otimes \BM(E_\Bb)$ into $\BM(E_\Aa \hat \otimes E_\Bb)$ given by
$$
(T_\Aa \hat \otimes T_\Bb)(\psi_\Aa \hat \otimes \psi_\Bb) 
\;=\; 
(-1)^{{\rm deg}(T_\Bb){\rm deg}(\psi_\Aa)} T_\Aa \psi_\Aa \hat \otimes T_\Bb \psi_\Bb
\;.
$$
This imbedding generates an isomorphism $\KM(E_\Aa) \hat \otimes \KM(E_\Bb) \cong \KM(E_\Aa \hat \otimes E_\Bb)$.
}
\hfill $\diamond$
\end{remark}

\begin{remark}
{\rm Note that if one of the algebras and its corresponding Hilbert module are ungraded, then the graded external product reduces to the ordinary one. However, if the remaining ones are graded, we will continue to use the hat just to indicate that the result is a graded Hilbert module over a graded algebra.
}
\hfill $\diamond$
\end{remark}

\section{The complex $K$-groups of C$^\ast$-algebras}
\label{chap-ComplexK}
 
In this chapter we briefly recall various characterizations of the $K$-groups of unital C$^\ast$-algebras. Of particular importance are the characterizations in terms of Hilbert C$^\ast$-modules and their various classes of operators, which enables one to imbed the $K$-groups in the broader family of Kasparov's $KK$-groups. The latter will provide yet another characterization of the $K$-groups and this will be elaborated in the following chapter.

\subsection{Classic characterization}

The $K_{0,1}$-groups of C$^\ast$-algebras classify the projections and the unitary elements w.r.t. the equivalence provided by the stable homotopy. Their standard definition can be found, {\it e.g.}, in \cite{RordamBook2000bv,BlackadarBook1998,WeggeOlsenBook1993de,HR}. Let $\Bb$ be a unital ungraded C$^\ast$-algebra and $\Pp_N(\Bb)$ the set of projections from $\CM(N) \otimes \Bb$. Consider the union 
\begin{equation}
\Pp_\infty (\Bb) 
\;=\; 
\bigcup_{N=1}^\infty \Pp_N(\Bb)
\;,
\end{equation}
where $\Pp_N(\Bb)$'s are considered pairwise disjoint. On $\Pp_\infty (\Bb)$, one introduces the addition operation
$$
p \oplus p'  
\;=\; 
\mathrm{diag}(p,p')
\;,
$$
which takes value in $\Pp_{N+M}(\Bb)$ when $p \in \Pp_N(\Bb)$ and $p' \in \Pp_M(\Bb)$. The space $\Pp_\infty (\Bb)$ accepts the following equivalence relation which is compatible with the addition,
\begin{equation}\label{Equivalence1}
\Pp_N(\Bb) \;\ni\; p \;\sim_0 \; p' \;\in \;\Pp_M(\Bb) 
\quad \Longleftrightarrow \quad 
\left \{
\begin{array}{l}
p \,=\, vv^* \\
 p'\,=\,v^*v
 \end{array}
 \right.
\end{equation}
for some $v$ from $\CM(N\times M) \otimes \Bb$.   

\begin{definition} The group $K_{0}(\Bb)$ is defined as the enveloping abelian Grothen\-dieck group of the semigroup $(\Pp_\infty (\mathcal A)/\sim_0, \oplus)$.
\end{definition}

As for the $K_1$-group, let $\mathcal U_N (\Bb)$ denote the group of unitary elements of $\CM(N) \otimes \Bb$ and let
\begin{equation}
\mathcal U_\infty (\Bb) 
\;=\; 
\cup_{N=1}^\infty \mathcal U_N(\Bb)
\;,
\end{equation}
where $\mathcal U_N(\Bb)$'s are again considered pairwise disjoint. Define the following binary operation
$$
u \oplus v 
\;=\; 
\mathrm{diag}(u,v)
\;,
$$
so that $u\oplus v \in \mathcal U_{N+M}(\Bb)$ when $u \in \mathcal U_N(\Bb)$ and $v \in \mathcal U_M(\Bb)$. The space $\mathcal U_\infty (\Bb)$ accepts an equivalence relation which is compatible with the binary operation,
$$ 
\mathcal U_N(\Bb) \;\ni\; u \sim_1 v \in \mathcal U_M(\Bb) 
\;\;\Longleftrightarrow \;\;
u\oplus 1_{K-N} \sim_h v\oplus 1_{K-M}
$$
for some $K\geq \max(N,M)$. Here, $\sim_h$ denotes the homotopy equivalence inside $\mathcal U_K(\Bb)$, namely, $u \sim_h v$ if $u$ can be continuously deformed into $v$, with respect to the topology of $M_K(\mathbb C) \otimes \Bb$, without ever leaving the unitary group $\mathcal U_K(\Bb)$. 

\begin{definition} The group $K_1(\Bb)$ is defined as the abelian group $(\mathcal U_\infty (\Bb)/\sim_1, \oplus)$.
\end{definition}

\begin{remark}{\rm If $\Bb$ is non-unital, the definitions of the $K$-groups are more involved (see \cite{RordamBook2000bv} for example), but they will not be needed here because we will deal exclusively with unital algebras.}
\hfill $\diamond$
\end{remark}

\subsection{Characterization in terms of Hilbert modules}

In this section, the algebra $\Bb$ is unital and ungraded.

\begin{proposition}[\cite{WeggeOlsenBook1993de}, pg.~264]
\label{Pro-K0G} 
The group $K_0(\mathcal B)$ can be equivalently characterized as
$$
K_0(\mathcal B) \;=\;
\big \{[P]_0-[Q] _0\ : \ P, Q  \in \KM(\mathcal H_\mathcal B) \big \}
\;,
$$
where $[ \,.\, ]_0$ indicates the classes under the Murray-von Neumann equivalence relation in $\mathbb B(\mathcal H_\mathcal B)$.
\end{proposition} 

\begin{proposition}
\label{Pr-FiniteRank} 
For a unital C$^\ast$-algebra $\Bb$, any projection from $\mathbb K(\mathcal H_\mathcal B)$ is finite-rank.
\end{proposition}

\proof We will take advantage of the fact that the finite range operators on $\Hh_\Bb$ form a double-sided ideal in $\Bb(\Hh_\Bb)$ and that this ideal is dense in $\KM(\Hh_\Bb)$. Let
\begin{equation}
\label{Eq-ApproxUnit}
P_N 
\;=\; 
\sum_{n=1}^N\Theta_{\psi_n,\psi_n}
\;, 
\qquad \psi_n =(0, \ldots,0,1,0 \ldots )
\;,
\end{equation}
be the standard approximate identity for $\KM(\mathcal H_\mathcal B)$, where $1$ in \eqref{Eq-ApproxUnit} sits at the $n$-th entry. Then, for any projection $P \in \mathbb K(\mathcal H_\mathcal B)$, $P P_N$ converges to $P$ in the norm topology of $\mathbb B(\mathcal H_\mathcal B)$ as $N \rightarrow \infty$. Consequently, for $N$ large enough,
\begin{equation}
R \;=\;
1\,-\,P (1-P_N)
\end{equation}
is invertible in $\mathbb B(\mathcal H_\mathcal B)$. Now, observe that $P R = PP_N$, hence $PR$ is a finite rank operator. Therefore $P=(P R) R^{-1}$ is finite-rank.
\qed

\vspace{.2cm}

By combining the above statement with Kasparov's stabilization theorem, more precisely with Corollary~\ref{Co-FiniteHM}, one arrives at the following equivalent characterization of the $K_0$-group:

\begin{corollary}\label{Co-K0Modules} 
The group $K_0(\mathcal B)$ can be represented in the form
$$
K_0(\mathcal B) 
\;=\;
\big \{[E_\mathcal B]_0-[E'_\mathcal B]_0 \ : \ E_\mathcal B, E'_\mathcal B \ \mbox{finitely generated projective $\mathcal B$-modules} \big \}
\;.
$$
This time, $[\,.\,]_0$ indicates the equivalence classes w.r.t. the isomorphism between the Hilbert modules and the abelian group structure on $K_0(\Bb)$ is induced by the direct sum of Hilbert modules.
\end{corollary}

\begin{remark}{\rm In the light of the above statements, it becomes clear that the elements of the $K_0(\Bb)$ group can be represented by compact projections on countably generated Hilbert $\Bb$-modules and one does not have to restrict to the standard Hilbert C$^\ast$-module as in Proposition~\ref{Pro-K0G}.} 
\hfill $\diamond$
\end{remark}

\begin{remark}
\label{rem-ProjDeg}
{\rm Even though $\Bb$ is ungraded, there is a natural grading on $\KM(\mathcal H_\mathcal B)=\KM\otimes \Bb$, inherited from the natural grading $\gamma_\KM$ on $\KM$ described in Example~\ref{Ex-GradedAlg4}. It was already pointed out in Example~\ref{Ex-ProjectiveHM} that every homogeneous projection in  $\KM\otimes \Bb$ has to be of degree $0$ (this holds for every graded algebra). Here we go one step further and show that every projection in  $\KM\otimes \Bb$ is homotopic to such a homogeneous projection, namely it is always possible to choose representatives of $K_0$-classes that are of degree $0$. For that, let $P$ be a finite rank projection which is decomposed over the eigenspaces of $\gamma_\KM=\Ad_{\sigma_3}$ as follows:
$$
P\;=\;\begin{pmatrix} a & b\\ c & d \end{pmatrix}
\;,
\qquad
\gamma_\KM (P)\;=\;\begin{pmatrix} a & -b\\ -c & d \end{pmatrix}
\;.
$$
This is identified with a $4\times 4$ block matrix is homotopy equivalent to yet another $4\times 4$ projection:
$$
P\;=\;\begin{pmatrix} a & b & 0 & 0\\ c & d & 0 & 0 \\ 0 & 0 & 0 & 0 \\ 0 & 0 & 0 & 0 \end{pmatrix}
\;\sim\;
\begin{pmatrix} a & 0 & b & 0\\ 0 & 0 & 0 & 0 \\ c & 0 & d & 0 \\ 0 & 0 & 0 & 0 \end{pmatrix}
\;,
$$
where the homotopy is given by a rotation in the $2$ and $3$ components by $\frac{\pi}{2}$. Now, indeed the r.h.s. is homogeneous of degree $0$ w.r.t. $\gamma_\KM=\Ad_{\sigma_3\oplus\sigma_3}$. 
} 
\hfill $\diamond$
\end{remark}

\begin{proposition}[17.5.4 in \cite{BlackadarBook1998}]
\label{Pro-K1G} 
The group $K_1(\mathcal B)$ of the complex $K$-theory of a unital C$^\ast$-algebra $\Bb$ can be characterized as 
$$
K_1(\Bb) \;=\;
\{[U]_1 \ : \ 
U\in\UM(\Hh_\Bb), \ U-1 \in  \KM(\Hh_\Bb) \}
\;.
$$
where $\UM(\Hh_\Bb)$ is the group of unitary operators over the standard Hilbert $\Bb$-module
$$
\UM(\Hh_\Bb) 
\;=\; 
\{ U \in \BM(\Hh_\Bb) \ : \ UU^\ast = U^\ast U = 1\}
\;,
$$
and $[ \,.\, ]_1$ indicates the homotopy classes in $\UM(\Hh_\Bb)$. The group structure is induced by the product of unitary operators. This is actually an abelian binary operation, as the non-trivial part of two finite-dimensional unitaries can be homotopically moved into orthogonal subspaces of $\Hh_\Bb$.
\end{proposition}

\begin{remark}
{\rm The condition $U-1 \in \KM(\Hh_\Bb)$ implies $U^*-1 \in \KM(\Hh_\Bb)$ and thus also $U^*-U \in \KM(\Hh_\Bb)$. Furthermore, $U^2 -1=U(U-U^\ast) \in \KM(\Hh_\Bb)$. In other words, $U$ in Proposition~\ref{Pro-K1G} is a symmetry modulo compact operators.
}
\hfill $\diamond$
\end{remark}

\begin{remark}{\rm The elements of $K_1(\Bb)$ can be characterized by unitaries over generic countably generated Hilbert $\Bb$-modules and not necessarily over $\Hh_\Bb$ as in Proposition~\ref{Pro-K1G}. Indeed, let $U \in \UM(E_\Bb)$ be such that $U-1 \in \KM(E_\Bb)$ for some countably generated Hilbert $\Bb$-module $E_\Bb$. Due to Kasparov's stabilization theorem, more precisely Corollary~\ref{Co-FiniteHM}, $E_\Bb \simeq P \,\Hh_\Bb$ with a projection $P\in\BM(\Hh_\Bb)$. The image $V$ of $U$ under this isomorphism can be extend to a unitary $\tilde U = V\oplus (1-P)$ on the whole $\Hh_\Bb$. This extension clearly satisfies $\tilde U - 1 \in \KM(\Hh_\Bb)$ and hence defines a class in the $K_1(\Bb)$ group as in Proposition~\ref{Pro-K1G}.
\hfill $\diamond$
}\end{remark}

\subsection{Characterization in terms of generalized Fredholm operators}

Choosing $E_\Bb$ to be the standard module $\Hh_\Bb$ in the exact sequence \eqref{eq-Calkin} and using Remark~\ref{Re-OpSM}, the exact sequence becomes
\begin{equation}
\begin{diagram}
0 & \rTo & \KM \otimes \Bb & \rTo{\mathfrak i}  & \Mm(\KM \otimes \Bb) & \rTo{ \mathrm{ev} } & \QM(\Hh_\Bb) & \rTo & 0 \; .
\end{diagram}
\label{eq-CalkinStandard}
\end{equation}
As is well known \cite{BlackadarBook1998,WeggeOlsenBook1993de,RordamBook2000bv}, any short exact sequence of C$^\ast$-algebras generate a six-term exact sequence at the level of $K$-theory:
\begin{equation}\label{SixTermDiagram}
\begin{diagram}
& K_0(\KM \otimes \Bb) & \rTo{\mathfrak i_\ast } & K_0(\Mm(\KM \otimes \Bb)) & \rTo{\ \ {\rm ev}_\ast \ \ } & K_0(\QM(\Hh_\Bb)) &\\
& \uTo{\rm Ind} & \  &  \ & \ & \dTo{\rm Exp} & \\
& K_1(\QM(\Hh_\Bb))  & \lTo{\ \ {\rm ev}_\ast} & K_1(\Mm(\KM \otimes \Bb)) & \lTo{\mathfrak i_\ast} & K_1(\KM \otimes \Bb) &
\end{diagram}
\end{equation}
The aim in the following is to characterize the $K$-groups in terms of the homotopy classes $[\mathbb F(\Hh_\Bb)]$ and $[\mathbb F_\sa(\Hh_\Bb)]$ of generalized Fredholm operators, and to provide explicit formulas for the connecting maps.

\begin{proposition}[\cite{WeggeOlsenBook1993de} p.~180] 
Any multiplier algebra has vanishing $K$-groups.
\end{proposition}

\begin{corollary}\label{Co-TwoIso} The six-term exact sequence of $K$-theory reduces to two isomorphisms:
\begin{align}
& \mbox{\rm Ind}\;:K_1(\mathbb Q(\mathcal H_\mathcal B))\;\to\; K_0(\KM \otimes \Bb)\;\cong \;K_0(\mathcal B)\;,
\label{eq-GenInd}
\\ 
& \mbox{\rm Exp}:K_0(\mathbb Q(\mathcal H_\mathcal B))\;\to\; K_1(\KM \otimes \Bb)\;\cong \;K_1(\mathcal B)\;,
\label{eq-GenExp}
\end{align}
where the two identifications on the r.h.s. follow from the stability of the $K$-groups.
\end{corollary}

\begin{remark}{\rm Unlike $\BM$, the algebra $\BM(\Hh_\Bb)$ is generally not a von Neumann algebra. As such, $\BM(\Hh_\Bb)$ is stable only under continuous functional calculus and the polar decomposition is not always available. This is one of the additional difficulties encountered when passing from ordinary Hilbert spaces to Hilbert C$^\ast$-modules which has to be dealt with next.}
\hfill $\Box$
\end{remark}


\begin{lemma}[\cite{PimsnerPopaVoiculescuJOT1980fj,MingoTAMS1987fg} and \cite{WeggeOlsenBook1993de}] 
\label{lem-polar}
For any $F \in \FM(\Hh_\Bb)$, there exists $G\in \FM(\Hh_\Bb)$ with  $F-G\in \KM(\Hh_\Bb)$ having a polar decomposition $G=W |G|$, where $|G| = \sqrt{G^\ast G}$ and $W\in\FM(\Hh_\Bb) $ is a partial isometry $W=W W^\ast W $ such that $W^\ast W |G| = |G| W^\ast W = |G|$. 
\end{lemma}

\begin{proposition} 
\label{prop-polar}
$[\mathbb F(\Hh_\Bb)]=[\FM_\gu(\Hh_\Bb)]=[\mathbb F_{\sym,-}(\widehat{\Hh}_\Bb)]\cong K_1(\QM(\Hh_\Bb))$.
\end{proposition}

\proof Let us first show that $\mathbb F_\gu(\Hh_\Bb)$ is a homotopy retract of $\mathbb F(\Hh_\Bb)$. Given $F,G,W \in \mathbb F({\mathcal H_\mathcal B})$ as in Lemma~\ref{lem-polar}, one defines a path in $\FM({\mathcal H_\mathcal B})$ by $t\in[-1,0]\mapsto F_t=F+(1+t)(G-F)$ and $t\in[0,1]\mapsto F_t=W|G|^{1-t}$. Note that here $|G|^{1-t}$ is defined by continuous functional calculus over the C$^\ast$-algebra $\BM(\Hh_\Bb)$, and that $F_1=W\in  F_\gu(\Hh_\Bb)$. This shows the first equality. For the second, one uses that odd operators on $\widehat{\Hh}_\Bb=\Hh_\Bb\oplus\Hh_\Bb$ are off-diagonal in the standard grading given in Example~\ref{Ex-GradedStandardHM1}. For a supersymmetry these off-diagonal operators have to be unitary, showing $\FM_\gu(\Hh_\Bb)$ is naturally identified with $\mathbb F_{\sym,-}(\widehat{\Hh}_\Bb)$. As to the third equality, by definition, homotopy classes of unitaries in $\QM(\Hh_\Bb)$ correspond to homotopy classes in $\FM_\gu(\Hh_\Bb)$.
\qed

\begin{definition}[Generalized Fredholm index]\label{Def-GenFredInd} Let $F,G,W \in \mathbb F(\mathcal H_\mathcal B)$ be as in Lemma~\ref{lem-polar}. Then 
$$
\mathrm{Ker} (G)\; = \;{\bf 1}-W W^\ast\;, 
\qquad 
\mathrm{Ker}(G^\ast)\; =\; {\bf 1} - W^\ast W 
\;,
$$
are two compact projections in $\mathbb K \otimes \mathcal B$ which define an element of the group $K_0(\mathcal B)$
$$
[\mathrm{Ker} (G)]_0 \,-\, [\mathrm{Ker} (G^\ast)]_0 \;\in\; K_0(\mathcal B)
\;.
$$
This element of $K_0(\Bb)$ is independent of the choice of  $G$ and is determined entirely by the class of $F$ in the Calkin algebra. The generalized Fredholm index $\Ind:\FM(\Hh_\Bb)\to K_0(\Bb)$ is then defined as
$$
\mathrm{Ind}(F )\; =\; [\mathrm{Ker} (F)]_0\, -\, [\mathrm{Ker}( F^\ast)]_0 
\;,
$$
where it is understood that possibly (irrelevant) compact perturbations are used to define the kernels. If $E_\Bb\cong P\Hh_\Bb$ is an arbitrary countably generated Hilbert $\Bb$-module, the definition of the generalized index map is extended to $\Ind:\FM(E_\Bb)\to K_0(\Bb)$ via
$$
\Ind(F)
\;=\;
\Ind(PFP+1-P)
\;,
\qquad
F\in \FM(E_\Bb)
\;.
$$
\end{definition}

\begin{remark} 
{\rm
The index of generalized Fredholm operators over Hilbert C$^\ast$-modules was developed in the works by Kasparov \cite{KasparovJOT1980yt,KasparovMUI1981ui}, Mi\v{s}\v{c}enco and Fomenko \cite{FomenkoMUI1980js}, and Pimsner, Popa, and Voiculescu \cite{PimsnerPopaVoiculescuJOT1980fj} and assembled in a final form by Mingo \cite{MingoTAMS1987fg}. A pedagogical exposition of all this can be found in \cite{WeggeOlsenBook1993de}. 
} 
\hfill $\diamond$ 
\end{remark}
 
\begin{remark}\label{rem-GradedFredInd} 
{\rm In view of the identification in Proposition~\ref{prop-polar}, the natural definition of the generalized Fredholm index of a supersymmetry
$F\in \mathbb F_{\sym,-}(\widehat{\Hh}_\Bb)$ is that of the generalized Fredholm index of the part of $F$ which takes the even grading sector into the odd one. This convention will be adopted throughout.}
\hfill $\diamond$
\end{remark} 

\begin{theorem}[\cite{MingoTAMS1987fg} and  \cite{BlackadarBook1998}]
\label{MingoTh} 
The generalized Fredholm index implements the first isomorphism of Corollary~\ref{Co-TwoIso}. When paired with Proposition~\ref{prop-polar}, this gives the isomorphism
$$
K_0(\mathcal B)\; \simeq \;[\mathbb F(\mathcal H_\mathcal B)]\;=\;[\mathbb F_{\sym,-}(\widehat{\Hh}_\mathcal B)]\;.
$$
\end{theorem}

\proof Given a class $[F]_1\in K_1(\QM(\Hh_\Bb))$, its index in \eqref{eq-GenInd} is, as usual, calculated using an arbitrary lift. If this lift is chosen to be the $G$ of Lemma~\ref{lem-polar}, one  recovers the  generalized Fredholm index. \qed

\begin{remark}[\cite{MingoTAMS1987fg}] {\rm Note that, among many other useful things, the above characterization says that, as for ordinary Fredholm operators, two generalized Fredholm operators belong to the same connected component of $\FM(\Hh_\Bb)$ precisely when they have the same generalized Fredholm index.}
\hfill $\diamond$ 
\end{remark}

\begin{proposition} 
\label{prop-SymRetract}
$[\mathbb \FM_\sa(\Hh_\Bb)]=[\mathbb F_{\sym,+}(\widehat{\Hh}_\Bb)]\cong K_0(\QM(\Hh_\Bb))$.
\end{proposition}

\proof Let us first note that $\Hh_\Bb \cong \Hh_\Bb \oplus \Hh_\Bb$ so that $\FM_\sa(\Hh_\Bb)\cong\FM_\sa(\Hh_\Bb)\oplus\FM_\sa(\Hh_\Bb)$. Furthermore $[\mathbb F_\sa(\Hh_\Bb)]=[\FM_\sym(\Hh_\Bb)]$ follows from functional calculus in $\BM(\Hh_\Bb)$. As even operators on $\widehat{\Hh}_\Bb$ are diagonal, this then implies the first equality. For the second identification, let the class in $[\mathbb F_{\sym,+}(\widehat{\Hh}_\Bb)]$ be specified by diag$(1-2P,1-2Q)$ with two projections $P, \, Q\in \QM(\Hh_\Bb)$. This class is mapped to $[P]_0-[Q]_0\in K_0(\QM(\Hh_\Bb))$.
\qed

\begin{remark}{\rm Given these identifications, the exponential map \eqref{eq-GenExp} now sends elements from $\FM_\sa(\Hh_\Bb)\cong \FM_\sa(\Hh_\Bb)\oplus \FM_\sa(\Hh_\Bb)$ into $K_1(\mathcal B)$. Let us recall that for any $T\in \mathbb F_\sa(\mathcal H_\mathcal B)$ it is always possible to find a $K\in \mathbb K(\mathcal H_\mathcal B)$ such that $(T+K)^*=T+K$. }
\hfill $\diamond$ 
\end{remark}

\begin{definition}[Generalized exponential map]
\label{Def-GenExp} 
Let $T =T^*\in \mathbb F_\sa(\mathcal H_\mathcal B)$. Then the generalized exponential map $\Exp:\FM_\sa(\Hh_\Bb)\to K_1(\Bb)$ is defined as
$$
\mbox{\rm Exp}(T)\;=\; \big [e^{-2\pi \imath T} \big ]_1 \;\in\; K_1(\mathcal B)
\;.
$$
If $E_\Bb\cong P\Hh_\Bb$ is an arbitrary countably generated Hilbert $\Bb$-module, the definition of the generalized index map is extended to $\Exp:\FM_\sa(E_\Bb)\to K_1(\Bb)$ via
$$
\Exp(T)
\;=\;
\Exp(PTP+1-P)
\;,
\qquad
T\in \FM_\sa(E_\Bb)
\;.
$$
\end{definition}

\begin{theorem}
\label{thm-ExpEquiv} 
The generalized exponential map index implements the isomorphism 
$$
K_1(\mathcal B)\; \simeq \;[\mathbb F_{\sym,+}(\widehat{\Hh}_\mathcal B)]\;=\;[\mathbb \FM_\sa(\Hh_\Bb)]\;.
$$
\end{theorem}

\begin{remark} 
\label{rem-ExpInv}
{\rm
In Theorem~\ref{Th-KK0ToK0} below, an inverse for the generalized index map is written out. Apparently, it is not so easy to find an explicit formula for an inverse to the generalized exponential map.
}
\hfill $\diamond$
\end{remark}

\section{Kasparov's $KK$-theory}
\label{chap-KK}

The purpose of this chapter is to introduce and characterize Kasparov's $KK$-groups, in particular, to establish the link with the $K$-theory of the previous chapter, as well as to introduce Kasparov's product and to illustrate how it works in several particular cases of interest.

\subsection{The $KK$-groups} 
\label{sec-KKGroups}

Throughout this section, $\Aa$ and $\Bb$ denote graded unital C$^\ast$-algebras and the symbol $[\, ,\, ]$ stands for graded commutators (see Definition~\ref{def-GradedComm}).

\begin{definition}[Kasparov cycles, \cite{KasparovMUI1981ui,BlackadarBook1998,JensenThomsenBook1991}]
\label{Def-KasparovM}  
A Kasparov $(\Aa,\Bb)$-cycle is a triple 
$$
\mathcal E 
\;=\;
\Big (E_\Bb,\;\pi: \Aa \rightarrow \BM(E_\Bb),\;F \in \BM(E_\Bb) \Big ) 
\; ,
$$
where $E_\Bb$ is a countably generated graded Hilbert $\Bb$-module and the representation $\pi$ and the operator $F$ are such that:

\begin{enumerate}[\rm (i)]

\item $F=-S_{E_\Bb} F S_{E_\Bb}$ ({\it i.e.} $F$ is homogeneous of graded degree 1),

\item $\pi \circ \gamma_\mathcal A = \mbox{\rm Ad}_{S_{E_\Bb}} \circ \pi$ ({\it i.e.} $\pi$ is a graded $\ast$-homomorphism),

\item $[F, \pi(a)] \in \KM(E_\Bb)$ for all $a\in \Aa$,

\item $(F^2 - 1)\pi(a) \in \KM(E_\Bb)$ for all $a\in \Aa$,

\item $(F^\ast -F) \pi(a) \in \KM(E_\Bb)$ for all $a\in \Aa$.

\end{enumerate}
The set of Kasparov $(\mathcal A,\mathcal B)$-cycles is denoted by $\mathbb E(\mathcal A,\mathcal B)$. A  Kasparov $(\mathcal A,\mathcal B)$-cycle is called degenerate if all operators in {\rm (iii)-(v)} vanish. The set of all degenerate cycles is denoted by $\mathbb D(\mathcal A,\mathcal B)$.
\end{definition}  

\begin{remark}{\rm For a unital algebra $\Aa$, conditions (iv) and (v) simply reduce to $F^2 - 1,F^\ast -F \in \KM(E_\Bb)$. As in Connes' quantized calculus, it will be essential below to work with special Kasparov cycles introduced next.
} \hfill $\diamond$ \end{remark}

\begin{definition}[Normalized Kasparov cycles, \cite{Connes:1994wk}]
\label{Def-KasparovNorm}
A Kasparov cycle is called normalized if $F^\ast = F$ and $F^2 = 1$. 
\end{definition} 

\begin{definition}[Kasparov's $KK^0$-group] Let $\mathcal A$ and $\mathcal B$ be graded C$^\ast$-algebras. The Kasparov group $KK^0(\mathcal A,\mathcal B)$ is defined as the group of homotopy equivalence classes of $\mathbb E(\mathcal A, \mathcal B)$, where the binary operation is the one induced by the direct sum of Hilbert C$^\ast$-modules.
\end{definition} 

\begin{remark} {\rm The homotopy equivalence is discussed in depth in \cite{BlackadarBook1998,JensenThomsenBook1991} and, since it does not play a direct role in the computations below, we will not elaborate on it at all. Let us note, however, that the zero element in $KK^0(\mathcal A,\mathcal B)$ is given by the degenerate cycles. Moreover, the inverse is provided by 
$$
-[(E_\Bb,\pi,F)]
\;=\;
[(-E_\Bb,\pi\circ\gamma_\Aa,-F)] 
\;,
$$ 
where $-E_\Bb$ is $E_\Bb$ with grading $-S_E$. One then finds that $KK^0(\mathcal A,\mathcal B)$ is abelian \cite{BlackadarBook1998,JensenThomsenBook1991}. Another important property of the $KK^0$-groups is the stability in each of the two arguments w.r.t. tensorizing with matrices or compact operators. For example, 
$$
KK^0(\KM \hat{\otimes} \mathcal A,\mathcal B)
\;=\;
KK^0(\CM(N) \otimes \Aa,\mathcal B)
\;=\;
KK^0(\mathcal A,\mathcal B)
\;,
$$ 
where an even grading is used on $\KM$ and $\CM(N)$.
} \hfill $\diamond$ \end{remark}

\begin{remark}\label{Re-SecGrading}{\rm By adding the degenerate $(\Aa,\Bb)$-cycle $
\left(\widehat \Hh_\Bb,0,\begin{pmatrix} 0 & 1 \\ 1 & 0 \end{pmatrix}\right)$
and using Kasparov's stabilization theorem, it is easy to see that any class in $KK^0(\Aa,\Bb)$ can be identified with a cycle of the form $(\widehat \Hh_\Bb,\pi,F)$. The implication of this observation is significant for, if one recalls from Example~\ref{Ex-GradedStandardHM1}  that the grading of $\widehat \Hh_\Bb$ is $\widehat S = S \oplus (-S)$, then a second grading $\widehat Q= \begin{pmatrix} 0 & S \\ S & 0 \end{pmatrix}$ can be defined on $\Hh_\Bb \oplus \Hh_\Bb$ such that $\widehat Q \,\widehat S = -\widehat S \,\widehat Q$. Hence, in general, the classes of $KK^0(\Aa,\Bb)$ can be represented by $(\Aa,\Bb)$-cycles $(E_\Bb,\pi,F)$ such that $E_\Bb$ accepts a second grading $Q_{E_\Bb}$ which anti-commutes with original grading $S_{E_\Bb}$ of $E_\Bb$.} 
\hfill $\diamond$ 
\end{remark}

\begin{remark}\label{Re-FProp}
{\rm Since we will be dealing mainly with unital algebras, further simplifications can be assumed. Indeed \cite[p.~152]{BlackadarBook1998}, then the classes of $KK^0(\Aa,\Bb)$ can always be represented by cycles with $F^\ast = F$, $\|F\| \leq 1$ and $F^2 -1 \in \KM(E_\Bb)$ and the representation $\pi$ can be assumed non-degenerate. If the algebras are ungraded, further simplifications occur as explained next.
} 
\hfill $\diamond$ 
\end{remark}

\begin{proposition}
\label{Pr-KK0Cyclebis} 
If $\Aa$ and $\Bb$ are ungraded unital C$^\ast$-algebras, the classes of $KK^0(\Aa,\Bb)$ can be represented by cycles of the form
\begin{equation}
\label{Eq-KK0Cycle}
\left(\widehat{\Hh}_\Bb=\Hh_\mathcal B \oplus \Hh_\mathcal B , \pi \oplus \pi, F=\begin{pmatrix} 0 & T^* \\ T & 0 \end{pmatrix} 
\right)
\;,
\qquad
S
=\begin{pmatrix} 1 & 0 \\ 0 & -1 \end{pmatrix}
\;,
\end{equation}
where $\pi : \mathcal A \rightarrow \mathbb B(\Hh_\mathcal B)$ is a $\ast$-homomorphism and $T\in\BM(\Hh_\Bb)$ is such that, with $\mbox{\rm ev}$ is as in \eqref{eq-Calkin},

\begin{enumerate}[\rm (i)]

\item $[\pi(a), T] \in\KM(\Hh_\Bb)$ for all $a\in\Aa$,

\item $\mbox{\rm ev}(T^*T)=\mbox{\rm ev}(TT^*)=1$ in $\QM(\Hh_\mathcal B)$.

\end{enumerate}

\end{proposition}

\proof This is just a rewriting of Definition~\ref{Def-KasparovM}. It represents the so-called Fredholm picture of $KK^0(\Aa,\Bb)$, see 17.5 in \cite{BlackadarBook1998}. \qed

\begin{theorem}
\label{Th-KK0ToK0} 
Let $\mathcal B$ be an ungraded unital C$^\ast$-algebra. Then there is the isomorphism
\begin{equation}\label{Eq-K0Group}
KK^0(\CM, \mathcal B) 
\;\simeq\; 
K_0(\mathcal B)
\end{equation} 
given by
\begin{equation}
\label{Eq-KK0ToK0}
KK^0(\CM, \mathcal B) \,\ni\, 
\left [ \left (\widehat \Hh_\Bb, s , \begin{pmatrix}
0 & T^* \\ T & 0
\end{pmatrix} \right ) \right ] 
\;\mapsto\;
{\rm Ind}\big (T\big ) \,\in\, K_0(\mathcal B)
\;.
\end{equation}
Here $s:\CM\to\BM(\widehat{\Hh}_\Bb)$ is the representation by scalar multiplication.
The inverse of this isomorphism is given by
\begin{equation}
\label{Eq-InvKK0ToK0}
\big [E_\Bb \big ]_0 - \big [ E'_\Bb\big ]_0 \,\in\, K_0(\mathcal B) 
\;\mapsto\;
\Big [ \big (E_\Bb \oplus E'_\Bb, s , 0\big ) \Big ] \; \in \; KK^0(\CM,\Bb) \;,
\end{equation}
where $E_\Bb$ and $E'_\Bb$ are finitely generated Hilbert $\Bb$-modules ({\it cf.} Corollary~\ref{Co-K0Modules}) and the grading of their direct sum is $1 \oplus (-1)$.
\end{theorem}

\proof Indeed, according to Proposition~\ref{Pr-KK0Cyclebis}, every $(\CM,\Bb)^0$-cycle can be represented as stated in above with $T\in \BM({\mathcal H}_\mathcal B)$ such that $T^*T=1=TT^*$ up to compact perturbations. In other words, $T\in\FM_\gu(\Hh_\Bb)$ and its generalized Fredholm index is well defined and takes values in $K_0(\Bb)$. Regarding the inverse, let $T \in \FM(\Hh_\Bb)$ such that ${\rm Ind}(T) =[{\rm Ker} ( T)]_0 -[{\rm Ker} (T^\ast)]_0$, where it is understood that the kernels are finitely generated. Then $T$ is equivalent in $[\FM(\Hh_\Bb)]$ with the generalized Fredholm operator $0: {\rm Ker} ( T) \rightarrow {\rm Ker} ( T^\ast)$ \cite[Proposition~3.27]{HigsonPSPM1990} and the affirmation follows. \qed

\begin{remark}\label{Re-K0Class} 
{\rm According to Corollary~\ref{Co-FiniteHM}, $E_\Bb \oplus E'_\Bb$ can always be written as $P\, \widehat \Hh_\Bb$, where $P$ is a compact projection from $\BM(\widehat \Hh_\Bb)$ of graded degree zero. If $P_\pm = \frac{1}{2}(1 \pm S)P$ are the sectors of $P$, then $P_+ \widehat \Hh_\Bb = E_\Bb$ and $P_- \widehat \Hh_\Bb = E'_\Bb$. Hence the elements of $K_0(\Bb)$ can be encoded in a single projection $P \in \KM(\widehat \Hh_\Bb)$ of graded degree zero and the element is explicitly represented by $[P_+]_0 - [P_-]_0$, as in Proposition~\ref{Pro-K0G}. The corresponding Kasparov cycle is simply $(P\, \widehat \Hh_\Bb, s, 0)$. Note that the class $[P_+]_0$ is encoded in the ungraded cycle $(P_+\, \Hh_\Bb, s, 0)$.
}
\hfill $\diamond$
\end{remark}

\begin{definition}[Higher Kasparov groups, \cite{BlackadarBook1998} pp. 149] 
\label{def-HigherKK}
For $k\geq 1$, one defines
$$
KK^{-k}(\Aa,\Bb)\;=\;KK^0(\Aa\, \hat \otimes\, \CM_k,\Bb)\;,
\qquad
KK^k(\Aa,\Bb)\;=\;KK^0(\Aa,\Bb\, \hat \otimes \,\CM_k).
$$
\end{definition} 

\begin{remark}{\rm  Due to the stability of the $KK$-functor and basic properties of the complex Clifford algebras, one immediately realizes that it is sufficient to consider $KK^{-1}(\mathcal A,\mathcal B)$ and $KK^{1}(\mathcal A,\mathcal B)$ which are isomorphic, too.  More generally, $KK^{k}(\mathcal A,\mathcal B)$ and $KK^{k+2}(\mathcal A,\mathcal B)$ are isomorphic for all $k\in\ZM$, a fact that is  often referred to as formal Bott periodicity \cite{BlackadarBook1998}.  As such, we will mainly concentrate on $KK^0$ and $KK^1$ groups. Let us stress that the graded external product is used in Definition~\ref{def-HigherKK}.
} 
\hfill $\diamond$ 
\end{remark}

\begin{proposition}
\label{pr-KK1Cycle1} 
Let $\Aa$ and $\Bb$ be ungraded unital C$^\ast$-algebras. Then the classes of $KK^1(\Aa,\Bb)$ can always be represented by cycles of the form
$$
\big ( E_\Bb \hat \otimes \CM_1, \pi \hat \otimes 1 , \DiracSA \hat \otimes \epsilon \big )\; ,
$$
where $E_\Bb$ is an ungraded Hilbert $\Bb$-module, $\pi : \Aa \rightarrow \BM(E_\Bb)$ is an ungraded C$^\ast$-algebra homomorphism and $\DiracSA=\DiracSA^*\in\BM(E_\Bb)$ is such that $\mbox{\rm ev}(\DiracSA)\in \QM(E_\Bb)$ is a symmetry and $\|H\|\leq 1$. Here $\mbox{\rm ev}$ is as in \eqref{eq-Calkin}. In the light of Example~\ref{ex-CliffMatIso}, the above cycle can also be written as (see also Theorem~7.2 in \cite{KNR})
$$
\big ( E_\Bb \hat \otimes (\CM \oplus \CM), \pi \hat \otimes (1 \oplus 1),  \DiracSA \hat \otimes (1 \oplus(-1))  \big ) 
\;=\;
\big (E_\Bb \oplus E_\Bb , \pi \oplus \pi, \DiracSA \oplus (-\DiracSA)  \big ) \; .
$$
On the r.h.s., $E_\Bb \oplus E_\Bb$ is a Hilbert $\Bb \hat \otimes \CM_1$-module with the standard odd grading.
\end{proposition}


\begin{proposition}[\cite{Connes:1994wk} p.~433]
\label{pr-KK1Cycle2} 
Let $\Aa$ and $\Bb$ be ungraded unital C$^\ast$-algebras. Then the classes of $KK^{-1}(\Aa,\Bb)$ can always be represented by $(\Aa \hat{\otimes} \CM_1,\Bb)$-cycles of the form
\begin{equation}
\label{Eq-KK1Cycle}
\big (E_\Bb \hat \otimes \CM^2, \pi \hat \otimes \sigma, \DiracSA \hat \otimes \sigma_2  \big )\;,
\end{equation}
where $\CM^2$ is with the even grading provided by $\sigma_3$, $E_\Bb$ is an ungraded Hilbert $\Bb$-module, and $\sigma:\CM_1 \rightarrow \CM(2)$ is given by $\sigma(\epsilon) = \sigma_1$. Furthermore,  $\pi : \Aa \rightarrow \BM(E_\Bb)$ is an ungraded C$^\ast$-algebra homomorphism and $\DiracSA=\DiracSA^*\in\BM(E_\Bb)$ is such that $\mbox{\rm ev}(\DiracSA)\in \QM(E_\Bb)$ is a symmetry  and $\|H\|\leq 1$.
\end{proposition}

\begin{theorem}
\label{Th-KK1ToK1} 
Let $\mathcal B$ be an ungraded unital C$^\ast$-algebra. Then there are isomorphisms
\begin{equation}\label{Eq-K1Group}
KK^1(\CM, \mathcal B) 
\;\cong\;
KK^{-1}(\CM, \mathcal B) 
\;\cong\; 
K_1(\mathcal B)
\;,
\end{equation} 
given by 
\begin{equation}
\label{Eq-KK1ToK1} 
\Big [ \Big ( E_\Bb \hat \otimes \CM_1 , s \hat \otimes 1, \DiracSA \hat \otimes \epsilon  \Big ) \Big ] 
\;\cong\; \Big [ \Big ( E_\Bb \hat \otimes \CM^2, s \hat \otimes \sigma, \DiracSA \hat \otimes \sigma_2  \big ) \Big ]
\;\mapsto\;
\big [ e^{-\imath \pi \DiracSA} \big ]_1
\;,
\end{equation}
where $E_\Bb$ is ungraded and $\DiracSA=\DiracSA^\ast\in\mathbb B(E_\mathcal B)$ is the lift of a symmetry $\mbox{\rm ev}(\DiracSA)\in\QM(E_\Bb)$  and $\|H\|\leq 1$, {\it cf.} Propositions~\ref{pr-KK1Cycle1} and \ref{pr-KK1Cycle2}, and the graded isomorphism $\sigma:\CM_1 \to \CM(2)$ is defined by $\sigma(\epsilon)=\sigma_1$. 
\end{theorem}

\proof The first isomorphism in \eqref{Eq-KK1ToK1}  is read off Propositions~\ref{pr-KK1Cycle1} and \ref{pr-KK1Cycle2}. 
As to the second, the symmetries $\mbox{\rm ev}(\DiracSA)\in\QM(E_\Bb)$ are in bijection with the projections in $\QM(E_\Bb)$ via $\mbox{\rm ev}(\DiracSA)\mapsto P = \frac{1}{2}(1+\mbox{\rm ev}(\DiracSA) )$. This establishes an isomorphism between $KK^1(\CM,\mathcal B)$ and $K_0(\mathbb Q(\mathcal H_\mathcal B))$.  Then \eqref{Eq-KK1ToK1} is just the  exponential map of the exact sequence \eqref{eq-Calkin}
$$
\mbox{\rm exp}( [P]_0)
\;=\;
\big[\exp(-2\pi\imath\,\tfrac{1}{2}(1+\DiracSA))\big]_1
\;=\;
\big[ -e^{-\imath \pi \DiracSA} \big]_1
\;=\;
\big[ e^{-\imath \pi \DiracSA} \big]_1
\;,
$$
which establishes the isomorphism between $K_0(\mathbb Q(\mathcal H_\mathcal B))$ and $K_1(\mathcal B)$, {\it cf.} equation \eqref{eq-GenExp}.
\qed

\subsection{Kasparov products}
\label{sec-KKProd}

\subsubsection{The notion of connection}
\label{sec-Connection}

The settings for this section are as follows. Let $E_\mathcal A$ and $E_\mathcal B$ be countably generated graded Hilbert C$^\ast$-modules over the graded unital C$^\ast$-algebras $\Aa$ and $\Bb$, respectively, and let $\pi : \mathcal A \rightarrow \mathbb B(E_\mathcal B)$ be a graded homomorphism so that $(E_\Bb,\pi,F_\Bb)$ is actually a Kasparov $(\Aa,\Bb)$-cycle. Recall from Definition~\ref{Def-KasparovM} that this automatically implies that $[F_\Bb,\pi(a)] \in \KM(E_\Bb)$ for all $a \in \Aa$ and that $\pi$ can always chosen to be non-degenerate. Further let $E_\mathcal A \hat \otimes_\mathcal A E_\mathcal B$ be the graded internal tensor product with its standard grading operator. Under these circumstances: 

\begin{definition}[\cite{JensenThomsenBook1991} p.~65, and \cite{BlackadarBook1998} 18.3]
\label{Def-Connection} 
A homogeneous element $G \in \mathbb B(E_\mathcal A \hat \otimes_\mathcal A E_\mathcal B)$ of graded degree one is called an $F_\mathcal B$-connection for $E_\mathcal A$ if: 

\vspace{0.1cm}

\begin{enumerate}[\rm (i)]

\item $T_{\psi_\mathcal A} F_\mathcal B -(-1)^{{\rm deg}(\psi_\Aa)}G T_{\psi_\mathcal A} \in \mathbb K(E_\mathcal B,E_\mathcal A \hat \otimes_\mathcal A E_\mathcal B)$,

\vspace{0.1cm}

\item $F_\mathcal B T^\ast_{\psi_\Aa} - (-1)^{{\rm deg}(\psi_\Aa)}T^\ast_{\psi_\Aa} G  \in \mathbb K(E_\Aa \hat \otimes_\mathcal A E_\Bb, E_\Bb),$

\end{enumerate}
for all ${\psi_\mathcal A} \in E_\Aa$ and where $T_{\psi_\Aa} : E_\Bb \rightarrow E_\Aa \hat \otimes_\Aa E_\Bb$ and its conjugate are, as in  Example~\ref{Ex-Tensorizing}, given by 
$$
T_{\psi_\Aa} (\psi_\Bb) \;=\; \psi_\Aa \hat \otimes_\Aa \psi_\Bb\;, 
\qquad T_{\psi_\Aa}^\ast (\phi_\Aa \hat \otimes_\Aa \psi_\Bb) \;=\; \pi( \langle \psi_\Aa,\phi_\Aa \rangle ) \psi_\Bb \; .
$$
\end{definition}

In the following we elaborate Proposition~2.2.5 in \cite{JensenThomsenBook1991}, summarized in Lemma~\ref{lem-ConnectionConstruction1} below, which gives a standard procedure to construct an $F_\Bb$-connection.

\begin{proposition} 
\label{pro-ConnectionConstruction1}
Assume that the representation $\pi$ actually satisfies $[\pi(a),F_\Bb]=0$ for all $a\in \Aa$, as in Remark~\ref{rem-picommute}. Then 
$$
G
\;=\;
1 \hat \otimes_\Aa F_\Bb 
\;\in\; 
\BM(E_\Aa \hat \otimes_\Aa E_\Bb)
$$ 
is well-defined and is an $F_\Bb$-connection for $E_\Aa$. 
\end{proposition}

\proof We can verify directly that
\begin{align*}
\big ( T_{\psi_\Aa} F_\Bb -(-1)^{{\rm deg}(\psi_\Aa)}G T_{\psi_\mathcal A} \big ) \psi_\Bb & 
\;=\; 
\psi_\Aa \hat \otimes_\Aa F_\Bb \psi_\Bb - (-1)^{{\rm deg}(\psi_\Aa)} (1 \hat \otimes_\Aa F_\Bb)(\psi_\Aa \hat{\otimes}_\Aa \psi_\Bb) \\
& 
\;=\; 
\psi_\Aa \hat \otimes_\Aa F_\Bb \psi_\Bb - (-1)^{2 {\rm deg}(\psi_\Aa)} \psi_\Aa \hat{\otimes}_\Aa F_\Bb \psi_\Bb 
\\
& 
\;=\; 
0 \; ,
\end{align*}
and
\begin{align*}
\big ( F_\mathcal B T^\ast_{\psi_\Aa} \,-\, & (-1)^{{\rm deg}(\psi_\Aa)}T^\ast_{\psi_\Aa} G \big ) (\phi_\Aa \hat{\otimes}_\Aa \psi_\Bb) \\
& \;=\; 
\big (F_\Bb \circ \pi( \langle \psi_\Aa,\phi_\Aa \rangle ) \big ) \psi_\Bb - (-1)^{{\rm deg}(\psi_\Aa)+{\rm deg}(\phi_\Aa)}T_{\psi_\Aa}^\ast (\phi_\Aa \hat \otimes_\Aa F_\Bb \psi_\Bb) \\
& \;= \;
\big ( (F_\Bb \circ \pi( \langle \psi_\Aa,\phi_\Aa \rangle )) \psi_\Bb - (-1)^{{\rm deg}(\psi_\Aa)+{\rm deg}(\phi_\Aa)} \pi(\langle \psi_\Aa,\phi_\Aa \rangle ) F_\Bb \big )\psi_\Bb \\
&\;=\; [F_\Bb,\pi(\langle \psi_\Aa,\phi_\Aa \rangle )] \psi_\Bb 
\\
&\;=\;0\;,
\end{align*}
which used that 
$$
{\rm deg}\big (\pi(\langle \psi_\Aa,\phi_\Aa \rangle ) \big )
\;=\;
{\rm deg} (\langle \psi_\Aa,\phi_\Aa \rangle ) 
\;=\; 
{\rm deg}(\psi_\Aa) + {\rm deg}(\phi_\Aa)
\;,
$$
where the defining relations of a graded module were used (see Definition.~\ref{Def-GradedHM}).
\qed

\begin{example}{\rm Proposition~\ref{pro-ConnectionConstruction1} finds a direct application to the case $\Aa = \CM$ and $E_\Aa = \widehat \Hh$, in which case one can automatically conclude that $G=1 \hat \otimes_\CM F_\Bb \in \BM(\widehat \Hh \hat \otimes_\CM E_\Bb)$ is an $F_\Bb$-connection for $\widehat \Hh$.
}
\hfill $\diamond$
\end{example}

\begin{proposition}
\label{pro-ConnectionConstruction2}
Let $E_\Aa=\Aa$ and recall the graded isomorphism $\Psi_1: \Aa \, \hat \otimes_\mathcal A \, E_\mathcal B\;\to\; E_\mathcal B$ described in Example~\ref{ex-GradedEquiv1}. Then 
$$
G\;=\;\Psi_1^{-1} \, F_\Bb \,\Psi_1 \;\in\; \BM(\Aa \, \hat \otimes_\Aa E_\Bb)
$$ 
is an $F_\Bb$-connection for $\Aa$.
\end{proposition}

\proof The first condition of Definition~\ref{Def-Connection} follows from
\begin{align*}
\big ( T_a F_\Bb - (-1)^{{\rm deg}(a)} GT_a \big ) \psi_\Bb 
& \;=\; a \hat \otimes_\Aa F_\Bb \psi_\Bb - (-1)^{{\rm deg}(a)} (\Psi_1^{-1} F_\Bb \Psi_1) \, a \hat \otimes_\Aa \psi_\Bb \\
& \;=\; 1 \hat \otimes_\Aa \pi(a) F_\Bb \psi_\Bb - (-1)^{{\rm deg}(a)} 1 \hat \otimes_\Aa F_\Bb \pi(a) \psi_\Bb \\
& \;=\; 1 \hat \otimes_\Aa [\pi(a),F_\Bb] \psi_\Bb \; ,
\end{align*}
combined with (iii) from Definition~\ref{Def-KasparovM}. As
\begin{align*}
\big ( F_\Bb T^\ast_a - (-1)^{{\rm deg}(a)}T^\ast_{a} G \big)(a '\hat \otimes_\Aa \psi_\Bb) 
& \;=\; \big (F_\Bb \pi(a^\ast a') - (-1)^{{\rm deg}(a)} \pi(a^\ast) F_\Bb \pi(a') \big ) \psi_\Bb \\
& \;=\; [F_\Bb,\pi(a^*)]\pi(a') \psi_\Bb 
\\
& \;=\;( [F_\Bb,\pi(a^*)] \circ \Psi_1 ) \, a' \hat \otimes_\Aa \psi_\Bb
\end{align*}
with $\Psi_1$ as in Example~\ref{ex-GradedEquiv1}, the second condition of Definition~\ref{Def-Connection}  is also satisfied. \qed

\begin{corollary}
\label{cor-ConnectionConstruction1}
Let $E_\Aa=\widehat{\Hh}_\Aa$ and recall the graded isomorphism $\Psi_2: \widehat \Hh_\Aa \hat \otimes_\Aa \, E_\Bb \;\to\;\widehat{\Hh}\,\hat{\otimes}_\CM\,E_\mathcal B$ described in Example~\ref{ex-GradedEquiv2}. Then 
$$
G\;=\;\Psi_2^{-1} \, (1 \hat \otimes_\CM F_\Bb)\,\Psi_2 \;\in\; \BM(\widehat \Hh_\Aa \hat \otimes_\Aa \, E_\Bb)
$$
is an $F_\Bb$-connection for $\widehat \Hh_\Aa$, where $1 \hat \otimes_\CM F_\Bb$ is defined as in Proposition~\ref{prop-internalProdOp}.
\end{corollary}

\proof Recall that $\widehat \Hh_\Aa \hat \otimes_\Aa \, E_\Bb = \widehat \Hh \hat \otimes_\CM (\Aa \hat \otimes_\Aa E_\Bb)$. Then Proposition~\ref{pro-ConnectionConstruction2} extends $F_\Bb$ to an $F_\Bb$-connection for $\Aa \hat \otimes_\Aa E_\Bb$ and Proposition~\ref{pro-ConnectionConstruction1} extends this connection to an $F_\Bb$-connection over $\widehat \Hh \hat \otimes_\CM (\Aa \hat \otimes_\Aa E_\Bb)$. The final connection has precisely the expression given in the statement.
\qed

\vspace{.2cm}

\begin{proposition}
\label{pro-ConnectionConstruction3}
Let $P \in \BM(E_\Aa)$ be a projection of graded degree 0. If $F \in \BM(E_\Aa \hat \otimes_\Aa E_\Bb)$ is an $F_\Bb$-connection for $E_\Aa$, then 
$$
G
\;=\;
(P \hat \otimes_\Aa 1) \, F \, (P \hat \otimes_\Aa 1) \;\in\; \BM(PE_\Aa \hat \otimes_\Aa E_\Bb)
$$
is an $F_\Bb$-connection for $PE_\Aa$.
\end{proposition}

\proof Take $\psi_\Aa \in P E_\Aa$. Then
\begin{align*}
\big (T_{\psi_\Aa} F_\mathcal B -(-1)^{{\rm deg}(\psi_\Aa)}G T_{\psi_\mathcal A} \big ) \psi_\Bb 
& \;=\; \psi_\Aa \hat \otimes_\Aa F_\Bb \psi_\Bb -(-1)^{{\rm deg}(\psi_\Aa)} (P \hat \otimes_\Aa 1) F \psi_\Aa \hat \otimes_\Aa \psi_\Bb \\
& \;=\; (P \hat \otimes_\Aa 1) \big ( \psi_\Aa \hat \otimes_\Aa F_\Bb \psi_\Bb -(-1)^{{\rm deg}(\psi_\Aa)} F \psi_\Aa \hat \otimes_\Aa \psi_\Bb \big ) \\
& \;=\; (P \hat \otimes_\Aa 1) \big ( T_{\psi_\Aa}F_\Bb -(-1)^{{\rm deg}(\psi_\Aa)} F T_{\psi_\Aa} \big ) \psi_\Bb\;,
\end{align*}
and the first condition in Definition~\ref{Def-Connection} follows because 
$$
(P \hat \otimes_\Aa 1)\KM(E_\Bb, E_\Aa \hat \otimes_\Aa E_\Bb) \;\subset \;\KM(E_\Bb, PE_\Aa \hat \otimes_\Aa E_\Bb)\;.
$$
Furthermore, if $\psi_\Aa, \phi_\Aa \in PE_\Aa \subset E_\Aa$, then
\begin{align*}
\big (F_\Bb T^\ast_{\psi_\Aa} - & (-1)^{{\rm deg}(\psi_\Aa)}T^\ast_{\psi_\Aa} G \big ) \phi_\Aa \hat \otimes_\Aa \psi_\Bb \\
&\;=\; (F_\Bb T^\ast_{\psi_\Aa})\phi_\Aa \hat \otimes_\Aa \psi_\Bb - (-1)^{{\rm deg}(\psi_\Aa)}T^\ast_{\psi_\Aa} (P \hat \otimes_\Aa 1) F) (P \hat \otimes_\Aa 1)\phi_\Aa \hat \otimes_\Aa \psi_\Bb \\
&\;=\; (F_\Bb T^\ast_{\psi_\Aa})\phi_\Aa \hat \otimes_\Aa \psi_\Bb - (-1)^{{\rm deg}(\psi_\Aa)}T^\ast_{\psi_\Aa} F) \phi_\Aa \hat \otimes_\Aa \psi_\Bb \\
& \;=\; \big ( F_\Bb T^\ast_{\psi_\Aa} (-1)^{{\rm deg}(\psi_\Aa)}T^\ast_{\psi_\Aa} F \big ) \phi_\Aa \hat \otimes_\Aa \psi_\Bb
\;,
\end{align*}
and the second condition in Definition~\ref{Def-Connection} follows, too. \qed

\begin{lemma}
\label{lem-ConnectionConstruction1}
Consider the conditions stated at the beginning of the section. Let $P \in \BM(\widehat{\Hh}_\Aa)$ be the projection of graded degree 0  such that $E_\Aa \cong P\,\widehat{\Hh}_\Aa$ as in Corollary~\ref{Co-FiniteHM}, and recall the graded isomorphism 
$$
\Psi_3 \,:\,E_\mathcal A \, \hat \otimes_\mathcal A \, E_\mathcal B\;\to\;\Pi(P)\big(\widehat{\Hh}\,\hat{\otimes}_\CM\,E_\mathcal B\big)
$$ 
from Example~\ref{ex-GradedEquiv3}. Then
$$
G\;=\;\Psi_3^{-1} \,(1 \hat \otimes_\CM F_\Bb)\,\Psi_3 \;\in\; \BM(E_\Aa \hat \otimes_\Aa E_\Bb)
$$
is an $F_\Bb$-connection for $E_\Aa$. More explicitly, if $E_\Aa=P\,\widehat{\Hh}_\Aa$ and $\Psi_2 : \widehat \Hh_\Aa \hat \otimes_\Aa E_\Bb \to   \widehat \Hh \,\hat \otimes_\CM E_\Bb$ as in  Example~\ref{ex-GradedEquiv2}, the connection can be written  as
$$
G\;=\;\Psi_2^{-1} \, \Pi(P)\,(1 \hat \otimes_\CM F_\Bb)\, \Pi(P) \,\Psi_2 \;\in \;\BM(P\,\widehat{\Hh}_\Aa \hat \otimes_\Aa E_\Bb)\;.
$$
\end{lemma}

\proof Starting from $F_\Bb$, Proposition~\ref{pro-ConnectionConstruction1} provides an $F_\Bb$-connection on $\widehat \Hh \hat \otimes_\CM E_\Bb$ in the form $1 \hat \otimes_\CM F_\Bb$. The isomorphism from Proposition~\ref{pro-ConnectionConstruction2} transfers this connection to $\widehat \Hh \hat \otimes_\CM (\Aa \hat \otimes_\Aa E_\Bb) = \widehat \Hh_\Aa \hat \otimes_\Aa E_\Bb$. Lastly, Proposition~\ref{pro-ConnectionConstruction3} then provides the connection on $(P \hat \otimes_\Aa 1)(\widehat \Hh_\Aa \hat \otimes_\Aa E_\Bb) = P\,\widehat \Hh_\Aa \hat \otimes_\Aa E_\Bb$, which is also isomorphic to $E_\Aa \hat \otimes E_\Bb$. Keeping track of the isomorphisms, the final connection can be seen to take the form appearing in the statement. The second expression can be verified by a direct calculation. \qed

\subsubsection{Definition and main properties of the Kasparov product}
\label{sec-Prod}

\begin{definition}[\cite{JensenThomsenBook1991} p.~69, and \cite{BlackadarBook1998,HR}] 
\label{def-KKprod}
Let 
\begin{equation}
\label{eq-ProdData}
\mathcal E_1 \;=\;(E_\mathcal B,\pi,F_\mathcal B) \;\in \;\mathbb E(\mathcal A,\mathcal B)\;, 
\qquad 
\mathcal E_2 \;=\;(E_\mathcal C,\pi',F_\mathcal C) \;\in\; \mathbb E (\mathcal B,\Cc)
\;.
\end{equation}
Then a triple 
$$
\mathcal E_{1,2} 
\;=\; 
(E_\mathcal B \, \hat \otimes_\mathcal B \, E_\mathcal C, \pi \, \hat \otimes_{\mathcal B} \, {\rm id},F) \;\in\; \mathbb E(\mathcal A,\Cc)
$$ 
is called a Kasparov product of $\mathcal E_1$ by $\mathcal E_2$ if

\begin{enumerate}[\rm (i)]

\item $F$ is an $F_\mathcal C$-connection for $E_\mathcal B$,

\vspace{0.1cm}

\item $(\pi \, \hat \otimes_\mathcal B \, {\rm id})[F_\mathcal B \, \hat \otimes_\mathcal B \, {\rm id}, F](\pi \, \hat \otimes_\mathcal B \, {\rm id}) \geq 0 \ {\rm mod} \ \mathbb K(E_\mathcal B \, \hat \otimes_\mathcal B \, E_\mathcal C)$.

\end{enumerate}

\end{definition}

\begin{theorem}[\cite{KasparovMUI1981ui,BlackadarBook1998,JensenThomsenBook1991,HR}]
Let $\mathcal A$, $\mathcal B$ and $\mathcal C$ be separable graded C$^\ast$-algebras and let $\mathcal E_1$ and $\mathcal E_2$ as in \eqref{eq-ProdData}. 

\begin{enumerate}[\rm (i)]

\item There always exists a Kasparov product $\mathcal E_{1,2} \in \EM (\Aa,\Cc)$ of $\mathcal E_1$ by $\mathcal E_2$. Furthermore, $\mathcal E_{1,2}$ is unique up to operator homotopy.

\item The map
$$
KK^0(\mathcal A,\mathcal B) \,\times\, KK^0(\mathcal B,\mathcal C) \;\rightarrow\; KK^0(\mathcal A,\mathcal C)
$$
given by
$$
([\mathcal E_1],[\mathcal E_2])\; \rightarrow \;[\mathcal E_{1,2}]
$$
is well-defined and it is associative and distributive.
\end{enumerate}

\end{theorem}

Given a $F_\Cc$-connection for $E_\Bb$, the following result provides a construction of an $F_\Cc$-connection satisfying also condition (ii) in Definition~\ref{def-KKprod}, which hence represents the Kasparov product.

\begin{proposition}[Proposition~18.10 in \cite{BlackadarBook1998}]
\label{prop-ConnectionToKKProd}
Let $\mathcal A$, $\mathcal B$ and $\mathcal C$ be separable graded C$^\ast$-algebras and let $\mathcal E_1$ and $\mathcal E_2$ as in \eqref{eq-ProdData}, with $F_\Bb=F_\Bb^*$ and $\|F_\Bb\|\leq 1$. Furthermore let $G$ be a $F_\mathcal C$-connection on $E_\Bb \hat \otimes_
\Bb E_\Cc$ for $E_\Bb$. Set
\begin{equation}
\label{eq-KKpr}
F
\;=\;
F_\Bb\hat \otimes_\Bb 1
\;+\;
\Big(\big(1-(F_\Bb)^2\big)^{\frac{1}{2}}\hat \otimes_\Bb 1\Big)G
\;.
\end{equation}
If $[F,\pi(\Aa) \hat \otimes_\Bb 1 ] \subset \KM(E_\Bb \hat \otimes_\Bb E_\Cc)$, it follows that $(E_\mathcal B \, \hat \otimes_\mathcal B \, E_\mathcal C, \pi \, \hat \otimes_{\mathcal B} \, {\rm id},F) \in \mathbb E(\mathcal A,\Cc)$ is a Kasparov $(\Aa,\Cc)$-cycle representing the Kasparov product $\mathcal E_{1,2}$.
\end{proposition}

\begin{remark} 
{\rm
In \eqref{eq-KKpr}, $(1-(F_\Bb)^2)^{\frac{1}{2}}$ is compact on $E_\Bb$, but this does {\em not} imply that $(1-(F_\Bb)^2)^{\frac{1}{2}}\hat \otimes_\Bb 1$ is compact on $E_\Bb \hat \otimes_\Bb E_\Cc$. Hence it is not possible to discard the second summand in \eqref{eq-KKpr}.
}
\hfill $\diamond$
\end{remark}

\begin{remark} 
{\rm
Given the definition of the higher Kasparov groups and the formal Bott periodicity, the Kasparov product accepts the following immediate extension:
$$
KK^i(\mathcal A,\mathcal B)\, \times \,KK^{j}(\mathcal B,\mathcal C)\; \rightarrow \;KK^{i+j}(\mathcal A,\mathcal C)
\;,
$$
where the indices are considered modulo 2. Several of these Kasparov products are written out explicitly below.
}
\hfill $\diamond$
\end{remark}

\subsection{Kasparov products involving $K$-groups}
\label{Sec-RelExamples}

\subsubsection{The product $K_0(\Aa)\times KK^0(\Aa,\Bb)$}
\label{Sec-0x0}

Here we assume that $\Aa$ and $\Bb$ are ungraded unital C$^\ast$-algebras and the isomorphism $KK^0(\CM, \mathcal B) 
\simeq K_0(\mathcal B)$, made explicit in Theorem~\ref{Th-KK0ToK0} and Remark~\ref{Re-K0Class}, will be freely used. The Kasparov product treated below is perhaps the simplest of all and can be found in many places in the literature (see {\it e.g.} \cite{HigsonPSPM1990}).

\begin{proposition}
\label{prop-K0KK0Prod}
The Kasparov product
$$
K_0(\Aa)\; \times\; KK^0(\mathcal A,\mathcal B)\; \rightarrow\; 
KK^0(\mathbb C,\mathcal B) \simeq K_0(\Bb)
$$
is of the form
\begin{align*}
\big [P \big ]_0\times 
\Big [ \big (E_\Bb,\pi,F_\Bb \big ) \Big ]
\, & =\, 
\Big [ \Big (\Pi(P)(\widehat \Hh \, \hat \otimes_\CM E_\Bb), s, \Pi(P)(1 \hat \otimes_\CM F_\Bb) \Pi(P) \Big )\Big ]  \\
\, & \simeq \,
{\rm Ind}\Big ( \Pi(P)(1 \hat \otimes_\CM F_\Bb) \Pi(P) \Big )
\;,
\end{align*}
where $\Pi : \BM(\widehat \Hh_\Aa) \rightarrow \BM(\widehat \Hh \, \hat \otimes_\CM E_\Bb)$ is the canonical homomorphism from Example~\ref{ex-GradedEquiv3}.
\end{proposition} 

\proof Indeed, the classes of $K_0(\Aa)$ can be represented by a Kasparov $(\CM,\Aa)$-cycle $(P\, \widehat \Hh_\Aa,s,0)$, where $P$ is such that $P_-=0$ (see Remark~\ref{Re-K0Class}). All conditions of Proposition~\ref{prop-ConnectionToKKProd} are satisfied and, since $F_\Aa$ is null, the Kasparov product is equal to $\big [ \big ( P\,\widehat \Hh_\Aa \hat \otimes_\Aa E_\Bb, s, G \big ) \big ]$, where $G \in \BM(P\,\widehat \Hh_\Aa \hat \otimes_\Aa E_\Bb)$ is the $F_\Bb$-connection for $P\,\widehat{\Hh}_\Aa$ from Lemma~\ref{lem-ConnectionConstruction1}. We can apply the isomorphism $\Psi_2$ from Example~\ref{ex-GradedEquiv2} on the above cycle without changing its class and the result is precisely the cycle written in the statement. \qed 


\subsubsection{The product $K_1(\Aa)\times KK^1(\Aa,\Bb)$}
\label{Sec-1x1}


In this section, $\Aa$ and $\Bb$ are ungraded unital C$^\ast$ algebras. The following result can be found in \cite{KNR}, and the proof below follows the same line of arguments.

\begin{proposition}
\label{prop-K1KK1Prod}
Let $\mathcal A$ and $\mathcal B$ be ungraded unital C$^\ast$-algebras. Then the Kasparov product
$$
K_1(\Aa)\, \times\, KK^1(\Aa,\Bb)\; \rightarrow\; 
KK^0(\mathbb C,\mathcal B) 
\;\simeq\; 
K_0(\Bb)
$$
is of the form
\begin{align*}
\big [U \big ]_1  \times \Big [\big (  E_\Bb \hat \otimes \CM_1,  \pi \hat \otimes 1 , \DiracSA_\Bb \hat \otimes \epsilon \big ) \Big ]
\;  =\;  \Big [ \big ( \widehat \Hh \, \hat \otimes_\CM E_\Bb, s, F \big )\Big ] 
\; \simeq \;
{\rm Ind}\big ( F \big )
\;,
\end{align*}
with
$$
F
\;=\;
\begin{pmatrix}0 & (1 \otimes_\CM P)\Pi(U^*)(1 \otimes_\CM P)+1 \otimes_\CM (1-P) 
\\ (1 \otimes_\CM P)\Pi(U) (1 \otimes_\CM P)+1 \otimes_\CM (1-P)  & 0\end{pmatrix}
\;,
$$
where $P=\frac{1}{2}(1+\DiracSA_\Bb)=P^*$ and $\Pi : \BM(\Hh_\Aa) \rightarrow \BM(\Hh \, \otimes_\CM E_\Bb)$ is the ungraded version of the canonical homomorphism from Example~\ref{ex-GradedEquiv3}.  Just as an advertisement to the reader, let us point out that $P$ is not a projection, but only $\ev(P)$ is.
\end{proposition} 

\vspace{.2cm}

\proof We recall from Propositions~\ref{Pro-K1G} and Theorem~\ref{Th-KK1ToK1} that the elements of $K_1(\Aa)$ can be represented by $U \in \UM(\Hh_\Aa)$ with $U-1 \in \KM(\Hh_\Aa)$, or by the class of a Kasparov $(\CM,\Aa \hat \otimes \CM_1)$-cycle 
$$
\Big [ \big ( \Hh_\Aa \hat \otimes \CM_1, s,  \DiracSA_\Aa \hat \otimes \epsilon \big ) \Big ] 
\;=\; 
\Big [ \big ( \Hh_{\Aa_{(1)}}, s , \DiracSA_\Aa \hat{\otimes} \epsilon \big ) \Big ]
\;,
$$
where $\DiracSA_\Aa \in\BM(\Hh_\Aa)$ is a symmetry modulo compact operators with $U = e^{-\imath\pi \DiracSA_\Aa}$ and $\|\DiracSA_\Aa\|\leq 1$. Note that $\Hh_{\Aa_{(1)}}$ is graded by $\gamma_{\Aa_{(1)}}=1\otimes\gamma_{\CM_1}$. Now ${\rm ev}(\DiracSA_\Aa \hat{\otimes} \epsilon)\in\QM(\Hh_{\Aa_{(1)}})$ is a symmetry. Setting $Q=\frac{1}{2}(\DiracSA_\Aa+1)\geq 0$, one thus has 
$$
{\rm ev}(\DiracSA_\Aa) 
\;=\;
2 \, {\rm ev}(Q) -1 
\;=\; 
- \cos\big (\pi \, {\rm ev}(Q)\big )
\;=\;
- {\rm ev} \big ( \cos (\pi Q ) \big )
\;,
$$
namely $\DiracSA_\Aa = -\cos(\pi Q)$ modulo compacts and $Q$ is a projection modulo compacts. Consequently
$$
\Big [ \big ( \Hh_{\Aa_{(1)}}, s ,  \DiracSA_\Aa \hat{\otimes} \epsilon \big ) \Big ]
\;=\;
\Big [ \big ( \Hh_{\Aa_{(1)}}, s , -\cos(\pi Q)\hat\otimes \epsilon \big ) \Big ]
\;.
$$
On the other hand, the generic elements of $KK^1(\Aa,\Bb)$ can be represented as ({\it cf.} Propositions~\ref{pr-KK1Cycle1} and \ref{pr-KK1Cycle2})
$$
\Big [ \big ( E_\Bb \hat \otimes \CM_1 , \pi \hat \otimes 1, \DiracSA_\Bb \hat \otimes \epsilon \big ) \Big ] 
\;=\; 
\Big [ \big ( E_\Bb \hat \otimes \CM^2, \pi \hat \otimes \sigma, \DiracSA_\Bb \hat \otimes \sigma_2 \big ) \Big ]
\;,
$$
where $\pi \hat{\otimes} \sigma$ is a graded representation of $\Aa_{(1)}$ on $E_\Bb \hat \otimes \CM^2$. Hence, we need to compute the Kasparov product in $KK^0(\CM,\Aa_{(1)}) \times KK^0(\Aa_{(1)},\Bb)$ of the form
$$
\Big [ \big ( \Hh_{\Aa_{(1)}}, s, -\cos(\pi Q)\hat\otimes \epsilon \big ) \Big ] 
\;\times \;
\Big [\big ( E_\Bb \hat \otimes \CM^2, \pi \hat \otimes \sigma, \DiracSA_\Bb \hat \otimes \sigma_1 \big ) \Big ]
\;,
$$
which is given by ({\it cf.} Proposition~\ref{prop-ConnectionToKKProd})
\begin{equation}
\label{eq-KK1Prod}
\Big [ \big ( \Hh_{\Aa_{(1)}} \hat \otimes_{\Aa_{(1)}} (E_\Bb \hat \otimes \CM^2) , s \hat \otimes_{\Aa_{(1)}} 1, F' \big ) \Big ] 
\,\in\, KK^0(\CM,\Bb)
\;,
\end{equation}
with, due to $(1- (\cos(\pi Q)\hat\otimes \epsilon)^2)^{\frac{1}{2}}=\sin(\pi Q)\hat\otimes 1$ resulting from the positivity of $Q$,
$$
F' \;=\; 
- (\cos(\pi Q) \hat \otimes \epsilon) \hat \otimes_{\Aa_{(1)}} 1 
\,+\, 
\big ( (\sin(\pi Q) \hat \otimes 1) \hat \otimes_{\Aa_{(1)}} 1 \big ) G
\;
$$
and 
$$
G \;=\; \Psi_2^\ast (1 \hat \otimes_\CM F_\Bb) \Psi_2\; ,
$$ 
as given in Proposition~\ref{cor-ConnectionConstruction1}. The input here is $F_\Bb = \DiracSA_\Bb \hat \otimes \sigma_2$ together with the isomorphism 
$$
\Psi_2 
\;:\; 
\Hh_{\Aa_{(1)}} \hat \otimes_{\Aa_{(1)}} (E_\Bb \hat \otimes \CM^2) 
\;\rightarrow\; 
\Hh \otimes_\CM (E_\Bb \hat \otimes \CM^2) 
\;\cong\; 
(\Hh \otimes_\CM E_\Bb) \hat \otimes \CM^2
$$ 
described in Example~\ref{ex-GradedEquiv2} with $\widehat \Hh$ replaced by $\Hh$. Next we apply the isomorphism $\Psi_2$ on the Kasparov cycle \eqref{eq-KK1Prod} and transform it into $( (\Hh\otimes_\CM E_\Bb) \hat \otimes \CM^2, s, \Psi_2 F' \Psi_2^\ast)$.  By using 
$$
\Psi_2 \big ( (Q \hat \otimes \xi) \hat \otimes_{\Aa_{(1)}} 1 \big ) \Psi_2 ^\ast
\;=\;
\Pi(Q) \hat \otimes \sigma(\xi) 
\;\in\; 
\BM \big ((\Hh \otimes_\CM E_\Bb) \hat \otimes \CM^2 \big )
\;,
$$
and the definition of $\sigma$, together with the notation $\tilde Q = \Pi(Q)$ and $\tilde P = 1 \otimes_\CM P$, one finds
\begin{align*}
\Psi_2 F' \Psi_2^\ast  
& \,=\, 
-  \cos(\pi \tilde Q) \hat \otimes \sigma_1
\,+\,
\big(\sin(\pi \tilde Q) \hat \otimes 1 \big)
\;\big( (2\tilde P-1) \hat \otimes \sigma_2\big)
\\
& 
\,=\, 
\big(- \cos(\pi \tilde Q) \hat \otimes \sigma_1+ \sin(\pi \tilde Q) \hat \otimes \sigma_2\big) \tilde P \hat \otimes 1
\,-\,
\big( \cos(\pi \tilde Q) \hat \otimes \sigma_1 + \sin(\pi \tilde Q) \hat \otimes \sigma_2 \big)(1 - \tilde P) \hat \otimes 1
\\
& 
\,=\,
\big(-\cos(\pi \tilde Q) \hat \otimes \sigma_1+ \sin(\pi \tilde Q) \hat \otimes \sigma_2\big) \tilde P \hat \otimes 1
\,+\,
\big(-\cos(\pi \tilde Q) \hat \otimes 1 +\imath \sin(\pi \tilde Q) \hat \otimes \sigma_3\big) (1- \tilde P) \hat \otimes \sigma_1
\\
& 
\,=\, 
\big(-\cos(\pi \tilde Q) \hat \otimes 1 +\imath \sin(\pi \tilde Q) \hat \otimes \sigma_3\big)
\left[
\cos(2\pi \tilde Q) \tilde P \hat \otimes \sigma_1 - \sin(2\pi \tilde Q) \tilde{P} \hat \otimes \sigma_2 
\,+\,
 (1- \tilde P) \hat \otimes \sigma_1
\right]
\,.
\end{align*}
Now $-\cos(\pi \tilde Q) \hat \otimes 1 +\imath\sin(\pi \tilde Q) \hat \otimes \sigma_3$ is unitary and homotopic to the identity by the path obtained by replacing $\tilde Q$ by $\tilde Q(t) = (1-t)\tilde Q -t$. It can be easily checked that $\Psi_2F'\Psi_2^\ast(t)$ respects all properties in Definition~\ref{Def-KasparovM}. Hence the Kasparov $(\CM,\Bb)$-cycle 
$(E_\Bb \hat \otimes \CM^2, s, F'')$ with  
$$
F''
\;= \;
\big(-\cos(\pi \tilde Q) \hat \otimes 1 +\imath \sin(\pi \tilde Q) \hat \otimes \sigma_3\big)^\ast \Psi_2 F' \Psi_2^\ast 
\;=\;
\cos(2\pi \tilde Q) \tilde P \hat \otimes \sigma_1 - \sin(2\pi \tilde Q)\tilde P \hat \otimes \sigma_2
\,+\,
  (1- \tilde P) \hat \otimes \sigma_1
$$ 
defines the same Kasparov class. Since $U =e^{-\imath 2\pi Q}$, we arrive at
$$
F''
\;=\;
\begin{pmatrix}0 & \Pi(U^*)(1 \otimes_\CM P)+1 \otimes_\CM (1-P) 
\\ \Pi(U) (1 \otimes_\CM P)+1 \otimes_\CM (1-P)  & 0\end{pmatrix}
\;,
$$
Lastly, recall that $1-U \in \KM(\Hh_\Aa)$ and $[\pi(a),P]\in \KM(E_\Bb)$ for all $a \in \Aa$, therefore $[\Pi(U),(1 \otimes_\CM P)]$ is compact and the statement now follows.
\qed

\subsubsection{The product $K_0(\Aa)\times KK^1(\Aa,\Bb)$}


\begin{proposition}
\label{prop-K0KK1Prod}
Let $\mathcal A$ and $\mathcal B$ be ungraded unital C$^\ast$-algebras. Then the Kasparov product
$$
K_0(\mathcal A)\; \times\; KK^1(\mathcal A,\mathcal B)\; \rightarrow\; 
KK^1(\mathbb C,\mathcal B) \simeq K_1(\Bb)
$$
is of the form
\begin{align}
[P]_0 \; \times \;
\Big [\big (  E_\Bb \hat \otimes \CM_1,  \pi \hat \otimes 1 , \DiracSA_\Bb \hat \otimes \epsilon \big ) \Big ]
& \;=\; 
\Big [ \big ((\Pi(P) \widehat{\Hh}\hat{\otimes}_\CM E_\mathcal B)\hat\otimes \CM_1, s\hat\otimes 1, \Pi(P) (1\hat{\otimes}_\CM\DiracSA_\Bb) \Pi(P) \hat\otimes \epsilon \big)\Big ] 
\label{eq-K0KK1Prod}
\\
& \;\simeq \;
\big [ e^{-\imath \pi  \Pi(P) (1\hat{\otimes}_\CM\DiracSA_\Bb) \Pi(P)} \big ]_1
\;,
\nonumber
\end{align}
with $\Pi:\BM(\widehat{\Hh}_\Aa)\to\BM(\widehat{\Hh}\hat{\otimes}_\CM E_\Bb)$ as in  Example~\ref{ex-GradedEquiv3}.
\end{proposition} 

\proof The direct application of Proposition~\ref{prop-K0KK0Prod} shows that the product on the l.h.s. of \eqref{eq-K0KK1Prod} is equal to
$$
\Big [\big(\Pi'(P) \widehat{\Hh}\hat{\otimes}_\CM (E_\Bb\hat{\otimes}\CM_1),s, \Pi'(P) (1\hat{\otimes}_\CM\DiracSA_\Bb \hat \otimes \epsilon)\Pi'(P)\big)\Big]
\;,
$$
where $\Pi':\BM(\widehat{\Hh}_\Aa)\to\BM(\widehat{\Hh}\hat{\otimes}_\CM (E_\Bb\hat{\otimes}\CM_1))$ as in  Example~\ref{ex-GradedEquiv3}. But $\widehat{\Hh}\hat{\otimes}_\CM (E_\Bb\hat{\otimes}\CM_1)\cong(\widehat{\Hh}\hat{\otimes}_\CM E_\Bb)\hat{\otimes}\CM_1$ and under this identification $\Pi'(P)\cong\Pi(P) \hat{\otimes} 1$ from which the first claim follows. The second isomorphism follows from Theorem~\ref{Th-KK1ToK1}.
\qed 


\subsubsection{The product $K_1(\Aa)\times KK^0(\Aa,\Bb)$}
\label{sec-K1KK0}

\begin{proposition}
\label{prop-K1KK0Prod}
Let $\mathcal A$ and $\mathcal B$ be ungraded unital C$^\ast$-algebras. Then the Kasparov product
$$
K_1(\mathcal A)\; \times\; KK^0(\mathcal A,\mathcal B)\; \rightarrow\; 
KK^1(\mathbb C,\mathcal B)\,\cong\,KK^0(\CM_1,\Bb\hat{\otimes}\CM(2))
$$
is of the form
\begin{equation}
\label{eq-K1KK0Prod}
\big [U \big ]_1  \times 
\Big [\big ({\Hh}_\Bb\oplus {\Hh}_\Bb, \pi\oplus\pi, F_\Bb \big)\Big]
\;=\; 
\Big [ \big ((\Hh\otimes {\Hh}_\Bb)\hat{\otimes}\CM(4) , \eta, F\big )\Big ]
\;,
\end{equation}
where $\eta$ is a representation of $\CM_1$, the grading on $\CM(4)$ is $\mbox{\rm Ad}_{\sigma_3\otimes 1}$, and
$$
F\;=\;\begin{pmatrix} 0 & V^* \\ V & 0\end{pmatrix}
\;,
\qquad
V\;=\;
P\begin{pmatrix} \Pi(U) & 0 \\ 0 & \Pi(U)^* \end{pmatrix}P\,+\,(1-P)
\;,
$$
with $P = \frac{1}{2}(1\otimes F_\Bb +1)$ and $\Pi : \BM(\Hh_\Aa) \rightarrow \BM(\Hh \, \otimes_\CM \Hh_\Bb)$ is the ungraded version of the canonical homomorphism from Example~\ref{ex-GradedEquiv3}.
\end{proposition} 

\proof Exactly as in the first part of the proof of Proposition~\ref{prop-K1KK1Prod} we begin by representing the class in $K_1(\Aa)$ specified by $U \in \UM(\Hh_\Aa)$ with $U-1 \in \KM(\Hh_\Aa)$, by the class of a Kasparov $(\CM,\Aa \hat \otimes \CM_1)$-cycle $\Big [ \big ( \Hh_\Aa \hat \otimes \CM_1, s,  \DiracSA_\Aa \hat \otimes \epsilon \big ) \Big ]$ where $\DiracSA_\Aa \in\BM(\Hh_\Aa)$ is a symmetry modulo compact operators with $U = e^{-\imath\pi \DiracSA_\Aa}$ and $\|\DiracSA_\Aa\|\leq 1$. Defining a projection modulo compacts $Q$ by $\DiracSA_\Aa = -\cos(\pi Q)$, we then obtained that $[U]_1\in K_1(\Aa)$ is represented by a Kasparov $(\CM,\Aa_{(1)})$-cycle $\Big [ \big ( \Hh_{\Aa_{(1)}}, s , -\cos(\pi Q)\hat\otimes \epsilon \big ) \Big ]$. For the computation of the Kasparov product, we will use the isomorphism   $KK^0(\Aa,\Bb)\cong KK^0(\Aa_{(1)},\Bb_{(1)})$ given
\begin{equation}
\label{eq-CM1iso}
\big (\widehat{\Hh}_\Bb, \pi\oplus\pi, F_\Bb \big)
\;\mapsto\;
\big (\widehat{\Hh}_\Bb\hat{\otimes}\CM_1, (\pi\oplus\pi)\,\hat{\otimes}\,\mu, F_\Bb\hat{\otimes}1 \big)
\;,
\end{equation}
where $\mu$ is the left multiplication with $\CM_1$. Note that the grading on the module $\widehat{\Hh}_\Bb\hat{\otimes}\CM_1$ is given by $\sigma_3\otimes\gamma_{\CM_1}$. Now the product 
$$
\Big [ \big ( \Hh_{\Aa_{(1)}}, s , -\cos(\pi Q)\hat\otimes \epsilon \big ) \Big ]
\times
\Big[\big (\widehat{\Hh}_\Bb\hat{\otimes}\CM_1, \pi\oplus\pi\,\hat{\otimes}\,\mu, F_\Bb\hat{\otimes}1 \big)
\Big]
$$
can be calculated using Propositions~\ref{prop-ConnectionToKKProd} and \ref{cor-ConnectionConstruction1} and is equal to
$$
\Big[\big (\Hh_{\Aa_{(1)}}\hat{\otimes}_{\Aa_{(1)}}\widehat{\Hh}_\Bb\hat{\otimes}\CM_1, s,F' \big)
\Big]
\;,
$$
where
$$
F'\;=\;
-\cos(\pi Q)\,\hat\otimes\, \epsilon\,\hat{\otimes}_{\Aa_{(1)}} 1
\;+\;
(\sin(\pi Q)\,\hat\otimes 1\,\hat{\otimes}_{\Aa_{(1)}} 1)(
\Psi^{-1}\,
1\,\hat{\otimes}_\CM F_\Bb\hat{\otimes}\,1\,
\Psi)
\;,
$$
with graded isomorphism $\Psi:\Hh_{\Aa_{(1)}}\hat{\otimes}_{\Aa_{(1)}}\widehat{\Hh}_\Bb\hat{\otimes}\CM_1\to
\Hh\hat{\otimes}_\CM\widehat{\Hh}_\Bb \hat{\otimes}\CM_1$ given by
$$
\Psi( \psi\, a\,\hat{\otimes}\, \xi\, \hat{\otimes}_{\Aa_{(1)}}\phi\,\hat{\otimes}\,\eta)
\;=\;
\psi\,\hat{\otimes}_\CM(\pi\oplus\pi)(a)\phi\,\hat{\otimes}\,\mu(\xi)\eta
\;,
$$
where $\psi\in\Hh$, $a\in\Aa$, $\xi,\eta\in\CM_1$, $\phi\in\widehat{\Hh}_\Bb$. Next let us introduce $\tilde{Q}\in\BM(\Hh\hat{\otimes}_\CM\Hh_\Bb)$ by  $\tilde{Q}\hat\otimes 1=\Psi (Q\,\hat\otimes 1\,\hat{\otimes}_{\Aa_{(1)}} 1)\Psi^{-1}=1 \hat\otimes_\CM(\pi\oplus\pi)(Q)$. Then we have
$$
\Psi F'\Psi^{-1}
\;=\;
-\cos(\pi\tilde{Q})\,\hat{\otimes}\,\mu(\epsilon)\;+\; (\sin(\pi\tilde{Q})\,\hat{\otimes} \,1)(1\,\hat{\otimes}_\CM F_\Bb\hat{\otimes}1)
\;.
$$
Hence the product is given by 
$$
\Big[\big (
\Hh\hat{\otimes}_\CM\widehat{\Hh}_\Bb \hat{\otimes}\CM_1,
s,
-\cos(\pi\tilde{Q})\,\hat{\otimes}\,\mu(\epsilon)\;+\; (\sin(\pi\tilde{Q})\,\hat{\otimes}\, 1)(1\,\hat{\otimes}_\CM F_\Bb\,\hat{\otimes}\,1)
\big)
\Big]
\;\in\;
KK^0(\CM,\Bb\hat{\otimes}\CM_1)
\;.
$$
Let us note that the grading is given by $1\hat{\otimes}\sigma_3\hat{\otimes}\gamma_{\CM_1}$. In order to read out the result, let us use again the equivalent of the isomorphism \eqref{eq-CM1iso}. Then the product is of the form
$$
\Big[\big (
\Hh\hat{\otimes}_\CM\widehat{\Hh}_\Bb \hat{\otimes}\CM_1\hat{\otimes}\CM_1,
s\hat{\otimes}\mu,
-\cos(\pi\tilde{Q})\,\hat{\otimes}\,\mu(\epsilon)\,\hat{\otimes}\, 1\;+\; (\sin(\pi\tilde{Q})\,\hat{\otimes}\, 1\,\hat{\otimes} \,1)(1\,\hat{\otimes}_\CM F_\Bb\,\hat{\otimes}\,1\,\hat{\otimes}\, 1)
\big)
\Big]
\;,
$$
%
in $KK^0(\CM_1,\Bb\hat{\otimes}\CM_1\hat{\otimes}\CM_1)$ with grading given by $1\hat\otimes\sigma_3\hat\otimes\gamma_{\CM_1}\hat\otimes\gamma_{\CM_1}$. Let us note that the compact commutator condition for the representation of $\CM_1$ is trivially satisfied, so that we will not take care of it in the following. Finally let us apply the isomorphism $\varphi:\CM_1\hat{\otimes}\CM_1\to\CM(2)$ from Example~\ref{ex-Cliffo}. Let us denote $\eta=s\hat\otimes\varphi(1\hat \otimes\mu)\varphi^{-1}$ the left representation of $\CM_1$ (which can be calculated explicitly, namely $\eta(\epsilon)$ is the right multiplication by $\sigma_1$). Further one can check that $\varphi(\mu(\epsilon)\hat\otimes 1)\varphi^{-1}=\sigma_2$ (left multiplication by $\sigma_2$). Therefore the product is given by
$$
\Big[\big (
\Hh\hat{\otimes}_\CM\widehat{\Hh}_\Bb \hat{\otimes}\CM(2),
\eta,
-\cos(\pi\tilde{Q})\,\hat{\otimes}\,\sigma_2\;+\; (\sin(\pi\tilde{Q})\,\hat{\otimes}\, 1)(1\,\hat{\otimes}_\CM F_\Bb\,\hat{\otimes}\,1)
\big)
\Big]
\;\in\;KK^0(\CM_1,\Bb\hat\otimes \CM(2))
\;,
$$
with grading $1\hat\otimes_\CM\sigma_3\hat\otimes\mbox{\rm Ad}_{\sigma_3}$, {\it cf.} Example~\ref{ex-Cliffo}. Next we can decompose $\Hh\hat{\otimes}_\CM\widehat{\Hh}_\Bb \hat{\otimes}\CM(2)=\big((\Hh{\otimes}_\CM{\Hh}_\Bb)\oplus(\Hh{\otimes}_\CM{\Hh}_\Bb)\big)\hat{\otimes}\CM(2)$ as well as $\tilde{Q}=\hat{Q}\oplus\hat{Q}$ with $\hat{Q}\in\BM(\Hh{\otimes}_\CM{\Hh}_\Bb)$ given as in the statement of the proposition. Now the product with values in $KK^0(\CM_1,\Bb\hat\otimes \CM(2))$ becomes
$$
\Big[\big (
\big((\Hh{\otimes}_\CM{\Hh}_\Bb)\oplus(\Hh{\otimes}_\CM{\Hh}_\Bb)\big)\hat{\otimes}\CM(2),
\eta,
\begin{pmatrix}
-\cos(\pi\hat{Q})\,\hat{\otimes}\,\sigma_2 & (\sin(\pi\hat{Q})\,\hat{\otimes}\, 1)(1\hat{\otimes} T^*\hat{\otimes} 1)
\\
(\sin(\pi\hat{Q})\,\hat{\otimes}\, 1)(1\hat{\otimes} T\hat{\otimes} 1)
& -\cos(\pi\hat{Q})\,\hat{\otimes}\,\sigma_2
\end{pmatrix}
\big)
\Big]
\;,
$$
where $F_\Bb=\binom{0 \;\; T^*}{T\;\;0\;}$.  The grading is here $\mbox{\rm Ad}_{\sigma_3}\hat\otimes\mbox{\rm Ad}_{\sigma_3}=\mbox{\rm Ad}_{\sigma_3\hat\otimes \sigma_3}$. Identifying $\CM(2) \hat{\otimes}\CM(2)$ with $\CM(4)$, the product takes the form
$$
\Big[\Big (
(\Hh\otimes\Hh_\Bb)\hat\otimes\CM(4),
\eta,
\begin{pmatrix}
0 & \imath\cos(\pi\hat{Q}) & \sin(\pi\hat{Q})\, 1{\otimes} T^* & 0
\\
-\imath\cos(\pi\hat{Q}) & 0 & 0 &  \sin(\pi\hat{Q})\, 1{\otimes}T^*
\\
\sin(\pi\hat{Q}) \,1{\otimes}T & 0 & 0 & \imath\cos(\pi\hat{Q})
\\
0 & \sin(\pi\hat{Q}) \,1{\otimes}T & -\imath \cos(\pi\hat{Q}) & 0 
\end{pmatrix}
\Big)
\Big]
\;,
$$
The grading is here implemented by $\diag(1,-1,-1,1)=\sigma_3\otimes\sigma_3$. Exchanging multiple rows and columns (implemented by a cyclic right shift) brings the grading into the desired form given by $\diag(1,1,-1,-1)=\sigma_3\otimes 1$, and the module becomes
$$
\Big[\Big (
(\Hh\otimes_\CM\Hh_\Bb)\hat\otimes\CM(4),
\eta',
\begin{pmatrix}
0 & V'' \\
V' & 0
\end{pmatrix}
\Big)
\Big]
\;,
$$
where $\eta'$ is the modified representation of $\CM_1$,
$$
V'
\;=\;
\begin{pmatrix}
-\imath \cos(\pi\hat{Q}) & \sin(\pi\hat{Q}) \, 1{\otimes}T^* 
\\
\sin(\pi\hat{Q}) \, 1{\otimes}T  & \imath\cos(\pi\hat{Q})
\end{pmatrix}
\;=\;
-\imath \cos(\pi\hat{Q})\otimes\sigma_3 \;+\;(\sin(\pi\hat{Q})\otimes 1)( 1\otimes F_\Bb)
\;,
$$
and similarly $V''=\imath \cos(\pi\hat{Q})\otimes\sigma_3 \;+\;(\sin(\pi\hat{Q})\otimes 1 )(1\otimes F_\Bb)$. Note that at this points no gradings are involved any longer in $V'$ and $V''$. Furthermore, $V''=(V')^*$ modulo compacts. We can now homotopically deform $V'$ and $V''$ modulo compacts. This is done by an argument analogous to the last part of the proof of Proposition~\ref{prop-K1KK1Prod}. Write $1\otimes F_\Bb=P-(1-P)$. Then
\begin{align*}
V'
& \,=\, 
-  \imath\cos(\pi \hat Q) \otimes \sigma_3(P+(1-P))
\,+\,
\sin(\pi \hat Q) \otimes 1 (P-(1-P))
\\
& 
\,=\, 
\big(-  \imath\cos(\pi \hat Q) \otimes \sigma_3+\sin(\pi \hat Q)  \otimes 1\big)P
\;+\;
\big(-  \imath\cos(\pi \hat Q)  \otimes \sigma_3-\sin(\pi \hat Q)  \otimes 1\big)(1-P)
\\
& 
\,=\,
\big(-  \imath\cos(\pi \hat Q)  \otimes \sigma_3-\sin(\pi \hat Q)  \otimes 1\big)
\left[
\big(\cos(2\pi \hat Q)  \otimes 1 +\imath \sin(2\pi \hat Q)\otimes\sigma_3\big)P+(1-P)
\right]
\\
& 
\,=\, 
\big(-  \imath\cos(\pi \hat Q)  \otimes \sigma_3-\sin(\pi \hat Q) \otimes 1\big)
\left[
\begin{pmatrix} U & 0 \\ 0 & U^*\end{pmatrix} P+(1-P)
\right]
\,.
\end{align*}
Hence $V'$ is homotopic to the last bracket, and as $P$ and $U$ commute up to compacts, the result now follows.
\hfill $\Box$


\section{The Generalized Connes-Chern characters}
\label{chap-GCCC}



\subsection{A look back at the classical Connes-Chern characters}


This section briefly discusses the classical Chern characters from Connes' work in non-commutative geometry. This serves several purposes. Firstly, we will show more concretely how the Connes-Chern character fits into Kasparov's $K$-theory. As such, it will be a source of further examples of non-trivial Kasparov cycles and will provide an opportunity to see the Kasparov product at work, even though in a restricted setting. Secondly, it introduces the structures which we plan to generalize together with several properties which are to be expected from these generalizations.


\subsubsection{The even Connes-Chern character}


\begin{definition}\label{Def-EvenFM} Let $\Aa$ be a trivially graded C$^\ast$-algebra. Then an even Fredholm module over $\Aa$ is a normalized graded Kasparov $(\Aa,\CM)$-cycle $(\widehat \Hh,\pi,F)$, cf. Definition~\ref{Def-KasparovNorm}.
\end{definition}

\begin{remark}{\rm As noted in Remark~\ref{Re-FProp}, to represent the classes of $KK^0$-groups one can always choose $F^\ast = F$, but the requirement $F^2 = 1$ is special. The quantized calculus developed by Connes depends on it in an essential way. The following statement, reproduced from Connes' book \cite{Connes:1994wk}, is classic in Atyiah's index theory on manifolds \cite{AtiyahJapan1970vc} and it was already recognized by Kasparov in a broader context \cite{KasparovMUI1981ui}. 
}\hfill $\diamond$\end{remark}

\begin{proposition}[\cite{Connes:1994wk} pp. 289] Let us consider the classical characterization of the $K_0$-group of a unital and ungraded C$^\ast$-algebra $\Aa$, and let $[P]_0 \in K_0(\Aa)$ be such that $P \in \Pp_N(\Aa)$. Let $(\widehat \Hh, \pi, F)$ be an even Fredholm module over $\Aa$ and $(\CM^N \hat \otimes \widehat \Hh, {\rm id} \hat \otimes \pi, 1 \hat \otimes F)$ be its straightforward extension over $\CM(N) \hat \otimes \Aa$, where $\CM^N$ is with the trivial grading. Then $({\rm id} \hat \otimes \pi)(P)(1 \hat \otimes F)( {\rm id} \hat \otimes \pi)(P)$ is a Fredholm operator on $\CM^N \hat \otimes \widehat \Hh$, in the classical sense, and there exists the group morphism from $K_0(\Aa)$ to $\ZM$ given by:
\begin{equation}\label{Eq-FredPair}
\Big \langle \big [P \big ]_0, \big ( \widehat H,\pi,F) \Big \rangle 
\;=\; 
{\rm Ind} \Big (({\rm id} \hat \otimes \pi)(P)(1 \hat \otimes F)( {\rm id} \hat \otimes \pi)(P)\Big)
\;,
\end{equation}
where ${\rm Ind}$ on the r.h.s. is the classical Fredholm index. Furthermore, the map in \eqref{Eq-FredPair} depends only on the class of $(\widehat \Hh,\pi,F)$ in $KK^0(\Aa,\CM)$. 
\end{proposition}

\begin{remark}{\rm It is not difficult to see that in fact the pairing in \eqref{Eq-FredPair} connects to a Kasparov product. Indeed, all projectors from $\KM(\Hh_\Aa)$ are finite rank hence are like $P$ in the above statement. Then the classes of $K_0$-group can be characterized by Kasparov cycles $\big (P \Hh_\Aa,s,0 \big)$, {\it cf.} Remark~\ref{Re-K0Class}, and from Section~\ref{Sec-0x0},
$$
\Big [\big (P \Hh_\Aa,s,0 \big) \Big ] \times \Big [ \big ( \widehat \Hh, \pi, F \big ) \Big ] 
\;=\; 
\Big [ \Big (\Pi(P)(\Hh \hat \otimes_\CM \widehat \Hh \, ), s, \Pi(P)(1 \hat \otimes_\CM F) \Pi(P) \Big )\Big ].
$$
On the righthand side, we can replace $\Hh$ by $\CM^N$ and the Fredholm operator becomes identical with the one in \eqref{Eq-FredPair}.
}
\hfill $\diamond$
\end{remark}

\begin{remark}{\rm  Let us recall that for a homogeneous Fredholm operator of degree 1, only half of the operator enters in the index (see Remark~\ref{rem-GradedFredInd}).  If one now interprets the pairing in \eqref{Eq-FredPair} as an analytic index (as in \cite{KhalkhaliEMS2009fj}), Then the next natural step is to search for the corresponding geometric index, which can be indeed defined at least for the special classes of Fredholm modules introduced next.}
\hfill $\diamond$
\end{remark}

\begin{definition}[Smooth subalgebra, \cite{RennieKTh2003vj}] A $\ast$-subalgebra $\mathscr A \subset \mathcal A$ of a unital C$^\ast$-algebra is called smooth if it is a Frechet sub-algebra which is stable w.r.t. the holomorphic functional calculus.
\end{definition}

\begin{remark}{\rm Smooth sub-algebras have the same $K$-theory as their parent C$^\ast$-algebra, see section 4.3 in \cite{KhalkhaliEMS2009fj}.
}
\hfill $\diamond$
\end{remark}

\begin{definition} An even Fredholm module $(\widehat \Hh,\pi,F)$ is said to be $n$-summable over a smooth sub-algebra $\mathscr A \subset \mathcal A$ if $[F,\pi(a)]^n$ is trace class for all $a \in \mathscr A$.
\end{definition}

\begin{definition}[Classic even Connes-Chern character, \cite{Connes:1994wk}] The even Connes-Chern character ${\rm Ch}_{\rm even} (\widehat \Hh, \pi, F)$ of a $n$-summable even Fredholm module over a smooth sub-algebra $\mathscr A \subset \Aa$ is defined as the class of the following cocycles in the periodic cyclic cohomology $HP_{\rm even}(\mathscr A)$ of $\mathscr A$,
\begin{equation}\label{Eq-ChernChar}
\uptau_{2m}(a_0,a_1,\ldots,a_{2m}) 
\;= \;
\Gamma_{2m}\, {\rm Tr}\big (S F[F,\pi(a_0)] \cdots [F,\pi(a_{2m})] \big )
\;,
\end{equation}
Here $2m+1 \geq n$ and $\Gamma_{2m} = \tfrac{(-1)^m m!}{2}$.
\end{definition}

\begin{remark} 
\label{rem-choiceGamma}
{\rm Note that the summability condition is needed precisely to assure that the operator inside the trace belongs to first Schatten class. Note also that, technically the commutators appearing above are the graded commutators but, since the algebra $\Aa$ is trivially graded, all $\pi(a_i)$ have graded degree 0 and the commutator reduces to the ordinary one.
}
\hfill $\diamond$
\end{remark}

\begin{remark} {\rm A cyclic cocycle $\uptau$ can be extended over the stabilization $\KM \otimes \Aa$ using the cup product ${\rm tr} \# \uptau$ with the standard trace over $\KM$, see \cite{Connes:1994wk}~p.~191. This will be routinely assumed in the following, even though we keep the notation unchanged.
}
\hfill $\diamond$
\end{remark}

\begin{theorem}[\cite{Connes:1994wk} p. 296]\label{Th-EvenIndexTh} The following index theorem holds:
\begin{equation}\label{Eq-IndexTh1}
{\rm Ind} \Big (({\rm id} \hat \otimes \pi)(P)(1 \hat \otimes F)( {\rm id} \hat \otimes \pi)(P)\Big) 
\;=\; 
\Big \langle  {\rm Ch}_{\rm even} (\mathcal H,\pi,F), {\rm Ch}^{\rm even}(P) \big \rangle
\;,
\end{equation}
where ${\rm Ch}^{\rm even}$ is the Connes-Chern character in periodic $K$-homology $HP^{\rm even}(\Aa)$ of $\Aa$. More explicitly,
\begin{equation}\label{IndexTh2}
{\rm Ind} \Big (({\rm id} \hat \otimes \pi)(P)(1 \hat \otimes F)( {\rm id} \hat \otimes \pi)(P)\Big) \;=\; \Lambda_{2m} \, \uptau_{2m}(P,\ldots, P)\;,
\end{equation}
where $\Lambda_{2m} = \frac{1}{m!}$.
\end{theorem}

\begin{remark}{\rm From the point of view of index theory, one interpretes \eqref{Eq-IndexTh1} as the equality between an analytic and a geometric index. Hence, \eqref{Eq-IndexTh1} tells us that geometric index is provided precisely by the pairing of the even Connes-Chern character with $K_0$-theory. 
}
\hfill $\diamond$
\end{remark}


\subsubsection{The odd Connes-Chern character}


\begin{definition}\label{Def-OddFM} Let $\mathcal A$ be a trivially graded C$^\ast$-algebra. An odd Fredholm module over  $\mathcal A$ is a normalized Kasparov $(\Aa,\CM_1)$-cycle. Given Proposition~\ref{pr-KK1Cycle1}, the data for such a module is encoded in an ungraded Kasparov cycles $(\Hh,\pi,\DiracSA)$, where $\Hh$ is the ungraded Hilbert space and $\DiracSA=\DiracSA^*$ is a symmetry, while the actual cycle is $(\Hh \hat \otimes \CM_1,\pi \hat \otimes 1, \DiracSA \hat \otimes \epsilon)$. 
\end{definition}

\begin{proposition}[\cite{Connes:1994wk} pp. 289] Let $(\Hh, \pi, \DiracSA)$ be an odd Fredholm module over a unital and ungraded C$^\ast$-algebra $\Aa$, and let $\big [ U \big ]_1 \in K_1(\Aa)$ be such that $U \in \Uu_N(\Aa)$. Let $(\CM^N \otimes \Hh,{\rm id} \otimes \pi,1 \otimes \DiracSA)$ be the straightforward extension of the module over $\CM(N) \otimes \Aa$ and set $P=\frac{1}{2}(1+\DiracSA)$. Then the operator $(1 \otimes P) ({\rm id}\otimes \pi)(U) (1 \otimes P)$ is Fredholm in the classical sense and there exists the group morphism from $K_1(\Aa)$ to $\ZM$ given by
\begin{equation}\label{Eq-FredPairOdd}
\Big \langle \big [ U \big ]_1,(\Hh,\pi,\DiracSA ) \Big \rangle 
\;=\; 
{\rm Ind} \Big ( (1 \otimes P) ({\rm id}\otimes \pi)(U) (1 \otimes P) \Big)\;.
\end{equation}
Furthermore, the map in \eqref{Eq-FredPairOdd} depends only on the class of $(\mathcal H,\pi,\DiracSA)$ in $KK^1(\Aa,\CM)$. 
\end{proposition}

\begin{remark}{\rm As in the even case, the map \eqref{Eq-FredPairOdd} can be connected to the Kasparov product. Indeed, according to Theorem~\ref{Th-KK1ToK1}, the class of $U$ can be represented by the Kasparov cycle $(\Hh \hat \otimes \CM_1, s, T_U \hat \otimes \epsilon)$, with $U = e^{-i\pi T_U}$ and 
$$
\Big [ (\Hh \hat \otimes \CM_1, s , T_U \hat \otimes \epsilon) \Big ]  \times \Big [(\mathcal H \hat \otimes \CM_1, \pi \hat \otimes 1, \DiracSA \hat \otimes \epsilon) \big ] 
\,=\, 
{\rm Ind} \Big ((1 \otimes_\CM P)\Pi(U)(1 \otimes_\CM P)+1 \otimes_\CM (1-P) \Big )
$$
as already seen in Section~\ref{Sec-1x1}. When $U \in \Uu_N(\Aa)$, the last operator can be easily seen to be identical with the Fredholm operator appearing in \eqref{Eq-FredPairOdd}.
}
\hfill $\diamond$
\end{remark}

\begin{definition} An odd Fredholm module $(\Hh,\pi,\DiracSA)$ is said to be $n$-summable over a smooth sub-algebra $\mathscr A \subset \mathcal A$ if $[\DiracSA,\pi(a)]^n$ is trace class for all $a \in \mathscr A$.
\end{definition}

\begin{definition}[Odd classic Connes-Chern character, \cite{Connes:1994wk}] The odd Connes-Chern character ${\rm Ch}_{\rm odd} (\mathcal H, \pi, \DiracSA)$ of a $n$-summable odd Fredholm module over a smooth sub-algebra $\mathscr A \subset \Aa$ is defined as the class of the following cocycles in the periodic cyclic cohomology $HP_{\rm odd}(\mathscr A)$ of $\mathscr A$, 
\begin{equation}\label{Eq-ChernCharOdd}
\uptau_{2m-1}(a_0,\ldots,a_{2m+1}) 
\;= \;
\Gamma_{2m-1}\; {\rm Tr}\big ( \DiracSA[\DiracSA,\pi(a_0)] \cdots [\DiracSA,\pi(a_{2m-1})] \big )
\;,
\end{equation}
where $2m \geq n$ and $\Gamma_{2m-1} = 2(-1)^m (m-\tfrac{1}{2})!$ with  $(m-\tfrac{1}{2})! = (m-\tfrac{1}{2}) \cdots \tfrac{1}{2}$.
\end{definition}

\begin{theorem}[\cite{Connes:1994wk}, p.~296]
\label{Th-OddIndexTh} 
The following index theorem holds:
$$
{\rm Ind} \Big ( (1 \otimes P) ({\rm id}\otimes \pi)(U) (1 \otimes P) \Big)
\;=\; 
\Big \langle  {\rm Ch}_{\rm odd} (\mathcal H,\pi,\DiracSA), {\rm Ch}^{\rm odd}(U) \Big \rangle
\;,
$$
where ${\rm Ch}^{\rm odd}$ is the odd Connes-Chern character in periodic $K$-homology $HP^{\rm odd}(\Aa)$ of $\Aa$. More explicitly,
$$
{\rm Ind} \Big ( (1 \otimes P) ({\rm id}\otimes \pi)(U) (1 \otimes P) \Big) 
\;=\; 
\Lambda_{2m-1} \, \uptau_{2m-1}(U^\ast-1,U-1,\ldots,U^\ast-1, U-1)
\;,
$$
where $\Lambda_{2m-1}=\frac{1}{2^{m-1} (2m-1)!!}$. 
\end{theorem}


\subsubsection{Final remarks and a look ahead}
\label{sec-finalrem}


After the introduction of the quantized calculus in his book, Connes presents \cite[p. 298]{Connes:1994wk} a diagram which captures very well the core of the computations. We reproduce it here in a slightly modified version, essentially as given by Khalkhali \cite[p. 184]{KhalkhaliEMS2009fj}:
\begin{diagram}
& \ K_\ast(\Aa)& \times &  KK^\ast (\Aa,\CM)  & \rTo{\rm Kasparov \ product} &  K_0(\CM) & \rTo{\ \ {\rm Tr} \ \ } &  \ZM &\\
& \dTo{{\rm Ch}^\ast} & \  &  \dTo{{\rm Ch}_\ast} & \ && \ & \dTo{} & \\
& HP^\ast(\Aa)  & \times & HP_\ast(\Aa)  & \rTo{\quad {\rm Connes' \ pairing}\qquad  } &&& \CM &
\end{diagram}
Here it is understood that only the class of finitely-summable Fredholm modules is considered. When this diagram is paired with the local formula of Connes and Moscovici for the Connes-Chern characters \cite{CONNES:1995rv}, the significance of the above diagram becomes more revealing, especially in the context of Kasparov's $K$-theory. On the one hand, the quantized calculus generates cyclic cocycles which pair integrally with $K$-theory and, on the other hand, it provides a concrete way to evaluate Kasparov products.

\vspace{0.2cm}

One straightforward generalization of the above  diagram is
\begin{equation}\label{Diag-Main}
\begin{diagram}
& KK^\ast(\CM,\Aa) & \times & KK^\ast(\Aa,\Bb)  & \rTo{\rm Kasparov \ product} &  K_0(\Bb) & \rTo{\ \ {\rm Tr} \ \ } &  \CM &\\
& \dTo{{\rm Ch}^\ast} & \  &  \dTo{{\rm Ch}_\ast} & \ && \ & \dTo{} & \\
& HP^\ast(\Aa)  & \times & HP_\ast(\Aa)  & \rTo{\quad {\rm Connes' \ pairing}\quad  } &&& \CM &
\end{diagram}
\end{equation}
where the canonical trace appearing in the upper line will be introduced in the following sections and the map $\mbox{\rm Ch}_*$ depends on $\Bb$. Diagram \eqref{Diag-Main} and the program associated with it becomes particularly helpful if it can be completed for separable C$^\ast$-algebras. Indeed, then the calculus is still quantized in the sense that the final numerical invariant obtained by pairing $K_0(\Bb)$ with the trace lies in a discrete subgroup of the real axis. The main goal of our work is to show that this whole program is possible, at least for crossed product algebras by $\ZM^k$.

\subsection{Generalized Calder\'on-Fedosov principle and formula}
\label{Sec-Fedosov}

The equality between the analytic and geometric indices in Theorems~\ref{Th-EvenIndexTh} and \ref{Th-OddIndexTh} can be established by using the Calder\'on-Fedosov principle and formula \cite{CalderonPNAS1967yt,FedosovTMMS1974fd} for ordinary Fredholm operators. A particularly insightful and detailed account of this aspect can be found in \cite{KhalkhaliEMS2009fj}. The goal of this section is to develop an equivalent principle and formula for generalized Fredholm operators. 

\subsubsection{The canonically induced trace} 

Most of the statements of this section are reproduced from \cite{LacaJFA2004gd}. Let $E_\mathcal B$ be a Hilbert $\mathcal B$-module and assume that the C$^\ast$-algebra $\mathcal B$ is unital and equipped with a faithful, continuous and normalized  trace $\mathcal T_\mathcal B$. Then this trace can be promoted to a unique lower-semicontinuous, densely defined trace ${\rm Tr}$ over $\mathbb B(E_\mathcal B)$ such that
$$
{\rm Tr}\big (\Theta_{\psi,\psi'}) 
\;=\; 
\mathcal T_\mathcal B(\langle \psi',\psi \rangle )
\;, 
\qquad \psi,\psi' \in E_\mathcal B
\;.
$$
This  trace ${\rm Tr}$ appears in the upper line of \eqref{Diag-Main} where elements of $K_0(\Bb)$ are, as in Proposition~\ref{Pro-K0G}, seen as finite rank projections in $\KM(\Hh_\Bb)$. Since ${\rm Tr}$ plays a central role later on, let us state its full characterization.

\begin{theorem}[\cite{LacaJFA2004gd}]\label{Th-TraceExtension1}  Let $\BM_+(E_\Bb)$ be the positive cone of the C$^\ast$-algebra $\BM(E_\Bb)$, {\it i.e.} the set of the positive elements of $\BM(E_\Bb)$. Set
\begin{equation}\label{Eq-HatTrace}
{\rm Tr}(T) 
\;=\; \
\sup\Big \{\sum_{\psi \in I} \mathcal T_\mathcal B \big ( \langle \psi, T\psi \rangle \big ) \ : \ \sum_{\psi \in I} \Theta_{\psi,\psi} \leq 1, \ I \ \mbox{finite subset of} \ E_\mathcal B \Big \}
\;,
\end{equation} 
for any $T \in \BM_+(E_\Bb)$. Then the following statements hold over $\BM_+(E_\Bb)$:

\begin{enumerate}[\rm (i)]

\vspace{0.1cm}

\item ${\rm Tr}$ is a strictly lower semicontinuous map. Moreover, 
$$
\liminf_k \Tt_\Bb\big ( \langle \psi,T_k \psi \rangle \big ) 
\;\geq\; \Tt_\Bb\big ( \langle \psi, T \psi \rangle \big ) \;\;\;  \forall \; \psi \in E_\Bb 
\;\;\;\;\Longrightarrow \;\;\;\;
\liminf_k {\rm Tr}\big ( T_k \big ) 
\;\geq\; {\rm Tr}\big ( T \big )
\;.
$$  

\item If $\{P_N\}_{N\geq 1}$ is an approximate unit in $\mathbb K(E_\mathcal B)$ as in \eqref{Eq-ApproxUnit}, then 
$$
{\rm Tr}(T) \;= \;\lim_{N \rightarrow \infty} {\rm Tr} (P_N T P_N)
\;.
$$

\end{enumerate}

\end{theorem} 

\begin{definition} The following sub-spaces of $\BM(E_\Bb)$ will play a central role in the following:
\begin{align*}
& \SM^+_1(E_\mathcal B) \;= \; \{ T\in \BM_+(E_\Bb) \ : \ {\rm Tr}(T) < \infty \}\;, \\
& \mathbb S_2(E_\mathcal B) \;=\; \{ T \in \mathbb B(E_\mathcal B) \ : \ {\rm Tr}(T^\ast T) < \infty \}\;, \\
& \SM_1(E_\Bb) \;= \;{\rm span} \big(\SM^+_1(E_\Bb)\big) 
\;=\; 
\SM_2(E_\Bb)^\ast \, \SM_2(E_\Bb)
\;.
\end{align*}
\end{definition} 

\begin{remark}{\rm We refrained from calling these operator spaces the trace-class spaces because the latter are rather defined by the completion of the above spaces in the appropriate Schatten norms induced by $\rm Tr$. Unlike the case of ordinary Hilbert spaces, these completions can be vastly larger than $\BM(E_\Bb)$. It is important to note that the computations below require at no point to appeal to these closures. Staying inside $\BM(E_\Bb)$ is one of the major differences between our derivation of the local formula for the generalized Connes-Chern characters and the ones obtained with the lower-semicontinuous spectral triples \cite{CareyAMS2014jf}. 
} 
\hfill $\diamond$
\end{remark}

\begin{theorem}[\cite{LacaJFA2004gd}]\label{Th-TraceExtension2} In the above settings:

\begin{enumerate}[\rm (i)]

\item ${\rm Tr}$ can be extended to a positive linear functional on $\mathbb S_1(E_\mathcal B)$.

\item For every pair $\psi, \psi' \in \mathbb E_\mathcal B$, we have $\Theta_{\psi,\psi'} \in \mathbb S_1(E_\mathcal B)$ and 
$$
{\rm Tr}(\Theta_{\psi,\psi'}) \;=\; 
\mathcal T_\mathcal B \big ( \langle \psi', \psi \rangle \big )
\;.
$$

\item ${\rm Tr}$ is a semifinite trace over the C$^\ast$-algebra $\BM(E_\Bb)$, namely each projection in $\BM(E_\Bb)$ can be approximated by an increasing net of projections with finite trace.

\item $\mathbb S_1(E_\mathcal B)$ and $\mathbb S_2(E_\mathcal B)$ are two-sided ideals in $\mathbb B(E_\mathcal B)$. Moreover, $\mathbb S_1(E_\mathcal B)$ is an essential ideal and, if $S,T\in \mathbb S_2(E_\mathcal B)$ or if $S \in \mathbb B(E_\mathcal B)$ and $T\in \mathbb S_1(E_\mathcal B)$, then ${\rm Tr}(ST) = {\rm Tr}(TS)$.

\end{enumerate}

\end{theorem}

\begin{remark}{\rm
When $\mathcal B= \CM$, the inclusions $\mathbb S_1(E_\mathcal B)\subset\mathbb S_2(E_\mathcal B) \subset \mathbb K(E_\mathcal B)$ hold, but this is no longer true in general. We are thankful to Sergey Neshveyev for confirming this statement.
}
\hfill $\diamond$
\end{remark}

\begin{proposition} 
\label{prop-KoTracePairing}
The trace $\Tr$ induces a pairing $\Tr:K_0(\Bb)\to\RM$ via
$$
\Tr([P]-[Q])
\;=\;
\Tr(P)\,-\,\Tr(Q)
\;.
$$
If $\Bb$ is separable, the pairing takes values in a discrete subgroup of $\RM$.
\end{proposition}

\proof Note that the canonical trace is invariant under isomorphisms of Hilbert C$^\ast$-modules (which are necessarily isometric, see Definition~\ref{def-HMIsomorphism}). Hence, we can fix the Hilbert C$^\ast$-module to be $\Hh_\Bb$. Recall from Proposition~\ref{Pr-FiniteRank} that both $P$ and $Q$ are finite rank so that $\Tr(P)$ and $\Tr(Q)$ are finite. Moreover, if $P\sim P'$ are homotopic in $\BM(\Hh_\Bb)$, then there exists a unitary $U\in\KM(\Hh_\Bb)^\sim$ such that $P'=UPU^*$ and hence $\Tr(P')=\Tr(P)$ by Theorem~\ref{Th-TraceExtension2}(iv). As a similar statement holds for $Q$, the first claim follows. The second is a consequence of the countability of $K_0(\Bb)$ for separable $\Bb$. 
\qed

\subsubsection{Generalized Calder\'on-Fedosov formula} 

The first formulation is for operators on the standard Hilbert module $\Hh_\Bb$.

\begin{theorem}[Generic Calder\'on-Fedosov formula]
\label{Th-GenericFedosov} 
Let $\mathcal B$ be an ungraded unital C$^\ast$-algebra equipped with a continuous normalized trace and let $T \in \mathbb B(\mathcal H_\mathcal B)$ with $\|T\| \leq 1$. If there exists a natural number $n$ such that
\begin{equation}\label{Eq-Fedosov1}
(1 - T^\ast T)^n\, , \ (1 - T T^\ast)^n \;\in\; \SM_1(\Hh_\Bb) \cap \KM(\Hh_\Bb)\;,
\end{equation}
then $T \in \mathbb F(\mathcal H_\mathcal B)$ and the trace of its generalized Fredholm index can be computed by the formula
\begin{equation}
\label{Eq-Fedosov2}
{\rm Tr}\Big ({\rm Ind}(T) \Big )  
\;= \;
{\rm Tr} \Big ((1 - T^\ast T )^n\Big) \, -\, {\rm Tr} \Big ((1 - T T^\ast )^n\Big)\;.
\end{equation}
\end{theorem}

\proof ({\rm i}) Assume that $(1 - T^\ast T)^n$ belongs to $\mathbb K(\mathcal H_\mathcal B)$. Taking $m>n$ of the form $m=2^k$, then $(1 - T^\ast T)^{m}$ is again in $\mathbb K(\mathcal H_\mathcal B)$. One can then apply successively the square root until reaching $1 - T^\ast T$. The latter is necessarily compact due to the following generic argument. Suppose $G \geq 0$ and $G \in \mathbb K(\mathcal H_\mathcal B)$. Then its square root can be defined by the ordinary continuous functional calculus and, moreover,
$$
\sqrt{G}
\;=\; 
\lim_{\epsilon \downarrow 0} \,\big (\epsilon + \sqrt{G}\big )^{-1}G\;,
$$
where the limit is in $\mathbb B(\mathcal H_\mathcal B)$. Since $\mathbb K(\mathcal H_\mathcal B)$ is an ideal, the operators inside the limit belong to $\mathbb K(\mathcal H_\mathcal B)$ and, since $\mathbb K(\mathcal H_\mathcal B)$ is closed, the limit is in $\mathbb K(\mathcal H_\mathcal B)$. A $k$-fold application of this fact shows that $1 - T T^\ast \in \mathbb K(\mathcal H_\mathcal B)$ and similarly for $1 - T^\ast T$. Then ${\rm ev} (T)$ is invertible in $\mathbb Q(\mathcal H_\mathcal B)$, with inverse given by $\ev(T^*)$, hence $T$ is a generalized Fredholm operator.

\vspace{.1cm}

({\rm ii}) Let us assume for the beginning that $T$ accepts a polar decomposition $T=W|T|$. The reader is encouraged to read through Chapters~15 and 17 of \cite{WeggeOlsenBook1993de}, which explain well why not all adjointable operators have  a polar decomposition. Recall that $|T| =\sqrt{T^\ast T}$ and
$$
W^\ast W |T| \;=\; |T| W^\ast W \;=\; |T|\;,
\qquad
W \;=\; WW^\ast W\;, 
\qquad 
W^\ast\;=\;W^\ast W W^\ast
\;.
$$
Then
\begin{align*}
& 1 - T T ^\ast\;=\; 1 - WW^\ast +W(1 - |T|^2 ) W^\ast \;,
\\
& (1 - WW^\ast)\big  (W(1 - |T|^2) W^\ast \big )
\;=\;0\;,
\\
& \big (W(1 - |T|^2)W^\ast \big )  (1 - WW^\ast)\;=\;0\;.
\end{align*}
Writing out the power, one deduces
\begin{equation}
\label{Eq-In1}
(1 - T T^\ast)^n 
\;=\; 
1 - WW^\ast\,+\,W(1 - |T|^2 )^n W^\ast\;.
\end{equation}
Using again $W^\ast W |T| = |T|$, it also follows that
$$
W^\ast W(1 - |T|^2 )^n
\;=\;
W^\ast W - 1 \,+\, (1 - |T|^2 )^n 
\;=\; 
-(1-W^\ast W ) \,+\, (1 - T^\ast T )^n
\;.
$$
As both summands are in $\SM_1(\Hh_\Bb)$, this together with the cyclicity of the trace and \eqref{Eq-In1} gives:
\begin{align*}
{\rm Tr} \big ((1 - T^\ast T )^n\big) \, -\, {\rm Tr} \big ((1 - T T^\ast )^n\big)
& \;=\; 
{\rm Tr}(1 - W^\ast W) 
\,+\,
{\rm Tr} \big (W^*W(1 - T^\ast T )^n\big) 
\,-\, 
{\rm Tr} \big ((1 - T T^\ast )^n\big)
\\
& \;=\; 
{\rm Tr}(1 - W^\ast W) 
\,+\,
{\rm Tr} \big (W(1 - T^\ast T )^nW^*\big) 
\,-\, 
{\rm Tr} \big ((1 - T T^\ast )^n\big)
\\
& \;=\;
{\rm Tr}(1 - W^\ast W) - {\rm Tr}(1 - W W^\ast).
\end{align*}
The affirmation is thus proved for the particular case when $T$ accepts a polar decomposition. 

\vspace{.1cm}

For the generic case, note that we already established that $T$ is unitary modulo $\mathbb K (\mathcal H_\mathcal B)$. Then Lemma~7.4 of \cite{PimsnerPopaVoiculescuJOT1980fj} applies which, as noted in \cite{MingoTAMS1987fg}, extends to the present context. It assures that $T(1-P_N)$ accepts the polar decomposition provided $N$ is large enough, where $P_N$ is the approximation of unity defined in \eqref{Eq-ApproxUnit}. Then $T(1-P_N)$ is a compact perturbation of $T$ with a polar decomposition, hence \eqref{Eq-Fedosov2} is already known to hold for $T(1-P_N)$. By definition, the two have the same generalized Fredholm index. Hence
$$
{\rm Tr} \Big ( {\rm Ind}(T)\Big ) 
\;=\; 
{\rm Tr}\big ((1 - (1-P_N)T^\ast T(1-P_N) \big ) 
\;-\; 
{\rm Tr}\big ((1 - T (1-P_N) T^\ast )^n\big)
\;.
$$
Using the cyclic property of the trace, one can check directly that all the terms containing $P_N$ cancel out between the two summands so that the  last equation reduces to \eqref{Eq-Fedosov2}.
\qed

\vspace{0.2cm}

\begin{remark}{\rm Theorem~\ref{Th-GenericFedosov} will next be formulated over generic ungraded and graded Hilbert C$^\ast$-modules. This is useful because it provides a great flexibility for various applications.
} 
\hfill $\diamond$
\end{remark}

 \begin{theorem}[Ungraded Calder\'on-Fedosov formula]\label{Th-UngradedFedosov} Let $\Bb$ be a trivially graded unital C$^\ast$-algebra equipped with a continuous normalized trace and let $E_\Bb$ be a countably generated and trivially graded Hilbert $\Bb$-module. Then Theorem~\ref{Th-GenericFedosov} holds if $\Hh_\Bb$ is replaced by $E_\Bb$.
\end{theorem}

\proof Using the Kasparov stabilization theorem for ungraded Hilbert modules, it follows that $E_\Bb \cong P\,\Hh_\Bb$, with $P$ a projection from $\BM(\Hh_\Bb)$. One can follow the arguments in Theorem~\ref{Th-GenericFedosov}, this time taking 
$$
\widetilde T \;=\; (1-P) + PTP \in \mathbb B(\mathcal H_\mathcal B)\;.
$$
Then $(1 - T^\ast T)^n, \  (1 - T T^\ast )^n\in \mathbb S_1(E_\mathcal B) \cap \mathbb K(E_\mathcal B)$ is equivalent to
$$
(1 - \widetilde T^\ast \widetilde T)^n \;,\ (1 - \widetilde T \widetilde T^\ast)^n 
\;\in\; \mathbb S_1(\mathcal H_\mathcal B) \cap \mathbb K(\mathcal H_\mathcal B)
\;,
$$
and the application of the general Calder\'on-Fedosov formula gives the desired result.\qed

\begin{theorem}[Graded Calder\'on-Fedosov formula]\label{Th-GradedFedosov} 
Let $\mathcal B$ be a graded unital C$^\ast$-algebra equip\-ped with a continuous normalized trace and let $E_\mathcal B$ be a countably generated Hilbert $\mathcal B$-module graded by $S$. Without loss of generality, we can assume the existence of a second grading $Q$ over $E_\Bb$ which anti-commutes with $S$ (see Remark~\ref{Re-SecGrading}). Let $T \in \mathbb B(E_\mathcal B)$, $S_\pm = \frac{1}{2}(1\pm S)$ and 
$$
T_{-+}
\;=\;
S_- T S_+\;, 
\qquad 
T_{+-}
\;=\;
S_+ T^\ast S_-\;, 
\qquad T_{-+}^\ast \;=\; T_{+-}\;.
$$ 
Then:

\begin{enumerate}[\rm (i)]

\item If there exists $n$ such that
\begin{equation}
\label{Eq-Fedosov3}
(S_+ - T_{+-} T_{-+})^n\;, \ (S_- - T_{-+} T_{+-} )^n 
\;\in\; 
\mathbb S_1(E_\mathcal B) \cap \mathbb K(E_\mathcal B)
\;,
\end{equation}
then $T_{-+} \in \mathbb F(E_\mathcal B)$.

\item Assume  that \eqref{Eq-Fedosov3} holds and let ${\rm Tr}$ be the canonical trace over $E_\mathcal B$. Then 
\begin{equation}
\label{Eq-Fedosov4}
{\rm Tr}\Big ( {\rm Ind}(T_{-+})\Big )  
\;=\; 
{\rm Tr} \big ((S_+ - T_{+-}T_{-+} )^n\big) 
\;-\; 
{\rm Tr} \big ( (S_- - T_{-+} T_{+-} )^n\big)
\;.
\end{equation}

\end{enumerate}
\end{theorem}

\proof Let $E_\Bb = E_\Bb^0 \oplus E_\Bb^1$, where $E_\Bb^0 \cong E_\Bb^1$ with $Q$ establishing the isomorphism. Let us define 
$$
\widetilde T 
\;=\; 
1-S_- \;+\; T_{-+} Q 
\;\in\; 
\BM (E_\Bb^1)
\;.
$$
Then
\begin{align*}
1 \;-\; \widetilde T^\ast \widetilde T & =1- \big ( 1-S_- + Q \, T_{+-} \big ) \big (1-S_-  + T_{-+} Q \, \big ) 
\\
& \;=\; S_- - Q T_{+-} T_{-+} Q 
\;=\; 
Q(S_+  -  T_{+-}T_{-+})Q  
\;,
\end{align*}
and
\begin{align*}
1 \;-\; \widetilde T \widetilde T^\ast & =1-  \big (1-S_-  + T_{-+} Q \big )\big ( 1-S_- + Q \, T_{+-} \big ) 
\\
& \;=\; S_- -  T_{-+} T_{+-}  
\;.
\end{align*}
%
%
The conclusion is that \eqref{Eq-Fedosov3} is equivalent to 
$$
(1 - \widetilde T^\ast \widetilde T)^n \;,\  (1 - \widetilde T \widetilde T^\ast)^n 
\;\in\; 
\mathbb S_1(E_\Bb^1) \cap \mathbb K(E_\Bb^1)
\;.
$$
Now, recall that ${\rm Ind}\big  (T\big ) = {\rm Ind}\big  (T_{-+}\big )$ with the later viewed as an operator from $E_\Bb^0$ to $E_\Bb^1$. Then ${\rm Ind}\big  (T\big ) = {\rm Ind}\big  (T_{-+}Q\big )$ and by applying  the Calder\'on-Fedosov formula for $\widetilde T$ from Theorem~\ref{Th-UngradedFedosov},
$$
{\rm Tr}\Big ({\rm Ind} \big (T \big ) \Big ) 
\;=\; 
{\rm Tr}\Big ( \big (S_- - T_{+-} T_{+-} \big )^n \Big )
\,-\, {\rm Tr} \Big (\big ( S_+ - T_{+-} T_{+-} \big )^n\Big )
\;,
$$
where the cyclic property of the trace was used to cancel out $Q$. 
\qed


\subsection{Generalized even Connes-Chern character}
\label{Sec-EvenGChern}


Throughout this section, $\Aa$ is an ungraded unital C$^\ast$-algebra and $\Bb$ is a ungraded unital C$^\ast$-algebra equipped with a continuous normalized trace. Further, ${\rm Tr}$ denotes the canonical trace over Hilbert $\Bb$-modules from Proposition~\ref{Th-TraceExtension1}. Also, all the constants appearing in the section are as previously defined. The characterization of the analytic index defined by the upper line of the Diagram~\ref{Diag-Main} is a direct consequence of the Kasparov product worked out in Proposition~\ref{prop-K0KK0Prod} combined with Proposition~\ref{prop-KoTracePairing}:

\begin{proposition}  
There exists the following group-pairing  between $K_0(\Aa)$ and $KK^0(\Aa,\Bb)$ with values in a discrete subgroup of $\RM$:
$$
\Big \langle \big [P\big ]_0,\Big [ \big (E_\Bb,\pi,F\big ) \Big ] \Big \rangle 
\;=\; 
{\rm Tr}\Big ( \big [P \big ]_0 \times \Big [ \big (E_\Bb,\pi,F\big ) \Big ] \Big )
\;=\; 
{\rm Tr}\Big ({\rm Ind} \big (\Pi(P) (1 \hat \otimes_\CM F) \Pi(P) \big ) \Big) \; ,
$$
where $P \in \KM(\Hh_\Aa)$ and and $\Pi: \BM(\Hh_\Aa) \rightarrow \BM(\Hh \hat \otimes E_\Bb)$ is the canonical homomorphism from Example~\ref{ex-GradedEquiv3} with $\widehat \Hh$ replaced by its ungraded version $\Hh$.
\end{proposition}

\begin{definition}[Finitely-summable Kasparov cycles]\label{Def-FSKasparovM} 
A Kasparov cycle $(E_\mathcal B, \pi, F) \in \EM(\Aa,\Bb)$ is said to be $n$-summable over a smooth sub-algebra $\mathscr A \in \mathcal A$ if, for all $a_i \in \mathscr A$ and $m \geq n$, 
$$
[F,\pi(a_1)] \cdots [F,\pi(a_m)] \;\in \;\SM_1(E_\Bb)\;.
$$ 
\end{definition}

\begin{remark} {\rm The commutators in Definition~\ref{Def-FSKasparovM} are all compact, {\it cf.} Definition~\ref{Def-KasparovM}, hence the product actually lands in $\SM_1(E_\Bb) \cap \KM(E_\Bb)$ and this ultimately allows to apply the Claderon-Fedosov principle. Note also that, since $\Aa$ is ungraded, $\pi(a)$ is of graded degree 0, for all $a \in \Aa$, hence the graded commutators appearing above and throughout this section reduce to the ordinary commutators.}
\hfill $\diamond$
\end{remark} 

\begin{definition}[Even generalized Connes-Chern character] Let $(E_\mathcal B, \pi, F)$ be a normalized graded Kasparov $(\mathcal A,\mathcal B)$-cycle which is $n$-summable over the smooth sub-algebra $\mathscr A \in \Aa$ with $n$ odd, and let $S$ be the grading operator of $E_\mathcal B$. Then the even generalized Connes-Chern character ${\rm Ch}_{\rm even}(E_\mathcal B, \pi, F)$ is defined as the class in the periodic cyclic cohomology of the following even cocycles
\begin{equation}
\label{Eq-EvenChernCh}
\uptau_{2m}(a_0,a_1,\ldots,a_{2m}) 
\;=\; 
\Gamma_{2m} {\rm Tr}\Big (S F[F,\pi(a_0)][F,\pi(a_1)] \cdots [F,\pi(a_{2m})] \Big )
\;,
\end{equation}
where $m$ is such that $2m+1 \geq n$.
\end{definition}

\begin{remark}{\rm One can verify that the r.h.s. of \eqref{Eq-EvenChernCh} is indeed a periodic cyclic cocycle by repeating the steps in the classical setting, {\it e.g.} by following \cite{KhalkhaliEMS2009fj} pp. 169-170. 
}
\hfill $\diamond$
\end{remark}

\begin{theorem}\label{Th-GenEvenIndTh} Let $(E_\mathcal B, \pi, F)$ be a normalized graded Kasparov $(\mathcal A,\mathcal B)$-cycle which is $n$-summable over a smooth sub-algebra $\mathscr A \subset \Aa$, with $n$ odd. Then the following index theorem holds,
$$
{\rm Tr}\Big ({\rm Ind} \Big (\pi(p) F \pi(p) \Big ) \Big)
\;=\;
\Lambda_{2m} \, \tau_{2m}(p,p,\ldots,p)   \;, \qquad 2m+1 \geq n \; .
$$
\end{theorem}

\proof Let $P = \pi(p)$ which, by Remark~\ref{rem-ProjDeg} can be assumed to be even, namely $PS = SP$. Hence $E'_\Bb=PE_\Bb$ is a countably generated graded Hilbert $\Bb$-module with grading operator $S' =PSP$ (as in Example~\ref{Ex-ProjectiveHM}).  Stabilizing as in Remark~\ref{Re-SecGrading} if necessary, one can assume the existence of a second grading $Q'$ on $E'_\Bb$ which anti-commutes with $S'$. We now follow closely the notation from Theorem~\ref{Th-GradedFedosov} and introduce
$$
S'_\pm \;=\; 
\tfrac{1}{2}(P\pm S')\;, 
\qquad 
S' \;=\; S'_+ - S'_-\;, 
\qquad P \;=\; S'_+ +S'_-\;.
$$
Also, let us set
$$
F' \;=\; PFP\;, 
\qquad 
S' F' \;=\; -F' S'\;, 
\qquad 
F'_{+-} \;=\; S'_+ F' S'_-\;, 
\qquad F'_{-+} \;=\; S'_- F' S'_+
\;.
$$ 
Then $S'_\pm F' S'_\pm =0$. Furthermore, from the following list of useful identities inside $\BM(E_\mathcal B)$
\begin{equation}\label{Eq-Id07}
F \big [ F,P \big ]\;=\; -\big [ F,P \big ]F\;, 
\qquad 
P \big [ F,P \big ]\;= \;\big [ F,P \big ](1-P)\;, 
\qquad
P \big [ F,P \big ] ^2 \;=\; P F P F P -P\;,
\end{equation}
one can derive
\begin{align}\label{Eq-KK32}
& (S'_+ - F'_{+-} F'_{-+})^{m+1} \;=\; (-1)^{m+1} S'_+ \big [F,P \big ]^{2m+2} S'_+
\;,
\\
\label{Eq-KK33} 
& (S_- - F'_{-+} F'_{+-})^{m+1} \;=\; (-1)^{m+1} S'_- \big [ F,P \big ]^{2m+2}S'_-
\;.
\end{align}
Given that $(E_\mathcal B, \pi,F)$ is $n$-summable, $[F,P]^{2m+2} \in \SM_1(E'_\Bb)$ provided that $n$ satisfies $2m+2 > n$.  In this case, the Calder\'on-Fedosov principle in the formulation \eqref{Eq-Fedosov3} applies to $F' \in \BM(E'_\mathcal B)$ and we can use the graded Calder\'on-Fedosov formula \eqref{Eq-Fedosov4}. Since the canonical trace on $E'_\mathcal B$ is simply ${\rm Tr}'(.) = {\rm Tr}(P(.)P)$, it follows from  \eqref{Eq-KK32}, \eqref{Eq-KK33} and $S'=S'_++S'_-$ that
$$
{\rm Tr} \big ( {\rm Ind}(PFP) \big) 
\;=\;
(-1)^{m+1}\,{\rm Tr} \Big (S P \big [ F,P \big ]^{2m+2} \Big )
\;.
$$
We can open one commutator without leaving the domain of the trace, so that the r.h.s. becomes
$$
(-1)^{m+1}\,{\rm Tr} \Big ( S P F P\big [ F,P \big ]^{2m+1} \Big ) 
\;+\; 
(-1)^m\,{\rm Tr} \Big ( S PF  \big [ F,P \big ]^{2m+1} \Big )
\;.
$$
The first term vanishes due to the cyclic property of the trace and the second identity in \eqref{Eq-Id07}. The second term can be rewritten as
$$
(-1)^m\,{\rm Tr} \Big ( S F P \big [ F,P \big ]^{2m+1}\Big )
$$
by using the first identity in \eqref{Eq-Id07}, the cyclic property of the trace as well as $SF = -SF$. It can further be rewritten as
$$
(-1)^m\,{\rm Tr} \Big ( S F  \big [ F,P \big ]^{2m+1} P\Big )
$$
by using $SP=PS$ and cyclic property of the trace. Thus
\begin{align*}
{\rm Tr} \big ( {\rm Ind}(PFP) \big) 
& 
\;= \;
\tfrac{(-1)^m}{2} \,{\rm Tr} \Big ( S F P \big [ F,P \big ]^{2m+1} 
\;+\; S F \big [ F,P \big ]^{2m+1} P \Big )
\\
&
\;=\;
\tfrac{(-1)^m}{2}\,{\rm Tr} \Big ( SF \big [  F,P \big ]^{2m+1} \Big ) 
\\
&
\;=\; 
\tfrac{1}{m!}\; \tau_{2m}(p,\ldots,p)
\;,
\end{align*}
where in the second equality we used the identity
$$
P \big [  F,P \big ]^{2m+1}  \;+\;\big [  F,P \big ]^{2m+1}P \;=\;
\big [  F,P \big ]^{2m+1}
\;.
$$
The statement now follows.\qed


\subsection{Generalized odd Connes-Chern character}
\label{Sec-OddGChern}


Throughout, $\Aa$ and $\Bb$ are assumed to be un-graded unital C$^\ast$-algebras and the latter to be equipped with a continuous and normalized trace. Let us start again by characterizing the analytic index associated with the upper line of the diagram~\eqref{Diag-Main} as it follows directly from the theory developed so far.  

\begin{proposition}  There exists the following group-pairing  between $K_1(\Aa)$ and $KK^1(\Aa,\Bb)$ with values in a discrete subgroup of $\RM$:
\begin{align*}
\Big \langle \big [U \big ]_1,\Big [ \big (E_\Bb \hat \otimes \CM_1, \pi \hat \otimes 1, \DiracSA \hat \otimes \epsilon \big ) \Big ] \Big \rangle 
\; &=\; 
{\rm Tr}\Big ( \big [U \big ]_1 \times \Big [ \big (E_\Bb \hat \otimes \CM_1, \pi \hat \otimes 1, \DiracSA \hat \otimes \epsilon \big ) \Big ] \Big ) \\
&=\; 
 {\rm Tr}\Big ({\rm Ind} \big ((1 \otimes_\CM P) \Pi(U) (1 \otimes_\CM P) \big ) \Big) \; ,
\end{align*}
where $U \in \UM(\Hh_\Aa)$ such that $U-1 \in \KM(\Hh_\Aa)$, $P = \tfrac{1}{2}(1+\DiracSA)$ and $\Pi: \BM(\Hh_\Aa) \rightarrow \BM(\Hh \otimes E_\Bb)$ is the ungraded version of the canonical homomorphism from Example~\ref{ex-GradedEquiv3}. 
\end{proposition}

\begin{definition}[Odd generalized Connes-Chern character] Let $\big (E_\mathcal B, \pi, \DiracSA \big )$ be a normalized ungraded Kasparov $(\mathcal A,\mathcal B)$-cycle which is $n$-summable over a smooth sub-algebra $\mathscr A$ of $\Aa$ with $n$ even. Then the odd generalized Connes-Chern character ${\rm Ch}_{\rm odd}(E_\mathcal B, \pi, \DiracSA)$ is defined as the class in the periodic cyclic cohomology of the following odd cocycles
\begin{equation}\label{Eq-EvenChernCh2}
\uptau_{2m-1}(a_0,a_1,\ldots,a_{2m-1}) 
\;=\; 
\Gamma_{2m-1}\, {\rm Tr}\Big (\DiracSA[\DiracSA,\pi(a_0)][\DiracSA,\pi(a_1)] \cdots [\DiracSA,\pi(a_{2m-1})] \Big )
\;,
\end{equation}
where $m$ is such that $2m \geq n$.
\end{definition}

\begin{remark}{\rm One can once again verify that the r.h.s. of \eqref{Eq-EvenChernCh2} is indeed a cyclic cocycle by repeating the steps in the classical setting, {\it e.g.} by following \cite{KhalkhaliEMS2009fj} pp. 166-168. 
}
\hfill $\diamond$
\end{remark}

\begin{theorem}\label{Th-GenOddIndTh} Let $(E_\mathcal B, \pi, \DiracSA)$ be a normalized un-graded Kasparov $(\mathcal A,\mathcal B)$-cycle which is $n$-summable over a smooth sub-algebra $\mathscr A \subset \Aa$, with $n$ even. Then the following index theorem holds,
\begin{equation}\label{Eq-GenOddIndexTh2}
{\rm Tr}\Big ({\rm Ind} \big (P \pi(u) P \big ) \Big)
\;=\;
\Lambda_{2m-1} \, \tau_{2m-1}(u^\ast-1,u-1,\cdots,u^\ast-1,u-1)   \;, \qquad 2m \geq n \; ,
\end{equation}
where $P = \tfrac{1}{2}(1+\DiracSA)$.
\end{theorem}

\proof Let $U = \pi(u)$ and let us start with the list of useful identities
\begin{equation}\label{Eq-Id17}
U^\ast \big [ U,P \big ]\;=\; -\big [ U^\ast,P \big ]U
\;,
\qquad 
P \big [ U,P \big ]\;=\; \big [ U,P \big ](1-P)\;, 
\qquad 
P \big [ U^\ast,P \big ]\;=\; \big [ U^\ast,P \big ](1-P) \; ,
\end{equation}
as well as
$$
P \big [ U^\ast,P \big ] \big [U,P \big ] 
\;=\; 
P(U^\ast P - P U^\ast)(UP -PU)P
\;=\; 
P U^\ast P U P -P\; ,    
$$
from which one can derive 
\begin{align*}
& \big(P - PU^\ast PUP \big )^m 
\;=\; 
(-1)^m P \Big ( \big [ U^\ast,P \big ] \big [U,P \big ] \Big )^m 
\;=\; 
\frac{(-1)^m}{2^{2m}} \,P \Big ( \big [ U^\ast,\DiracSA \big ] \big [U,\DiracSA \big ] \Big )^m\;,
\\
& 
\big(P - PUPU^\ast P \big )^m 
\;=\; 
(-1)^m P \Big ( \big [ U,P \big ] \big [U^\ast,P \big ] \Big )^m 
\;=\; 
\frac{(-1)^m}{2^{2m}} \,P \Big ( \big [ U,\DiracSA \big ] \big [U^\ast,\DiracSA \big ] \Big )^m
\;.
\end{align*}
Provided that $2m \geq n$, the Calder\'on-Fedosov principle of Theorem~\ref{Th-UngradedFedosov} now applies and shows 
$$
{\rm Tr}\Big ( {\rm Ind}\big (PUP) \Big ) 
\;=\; 
\frac{(-1)^m}{2^{2m}}\Big [{\rm Tr}\Big (P \big ( \big [ U^\ast,\DiracSA \big ] \big [U,\DiracSA \big ] \big )^m \Big ) - {\rm Tr}\Big (P \big ( \big [ U,\DiracSA \big ] \big [U^\ast,\DiracSA \big ] \big )^m \Big )\Big ]\; .
$$
At this point, the $P$ in front of the traces can be written as $\frac{1}{2}(1+\DiracSA)$. Using  
$$
{\rm Tr}\Big (\big ( \big [ U^\ast,\DiracSA \big ] \big [U,\DiracSA \big ] \big )^m \Big ) 
\;=\; 
{\rm Tr}\Big (U^\ast\big ( \big [ U,\DiracSA \big ] \big [U^\ast,\DiracSA \big ] \big )^m U\Big ) 
\;=\; {\rm Tr}\Big (\big ( \big [ U,\DiracSA \big ] \big [U^\ast,\DiracSA \big ] \big )^m \Big )
$$
as well as 
$$
{\rm Tr}\Big (\big (\DiracSA  \big [ U,\DiracSA \big ] \big [U^\ast,\DiracSA \big ] \big )^m \Big ) 
\;=\; 
- \,{\rm Tr}\Big (\DiracSA \big ( \big [ U^\ast,\DiracSA \big ] \big [U,\DiracSA \big ] \big )^m \Big )
\;,
$$
one can now conclude that
$$
{\rm Tr}\Big ( {\rm Ind}\big (PUP) \Big ) 
\;=\; 
\frac{(-1)^m}{2^{2m}}\;{\rm Tr}\Big (\DiracSA \big ( \big [ U^\ast,\DiracSA \big ] \big [U,\DiracSA \big ] \big )^m \Big ) \; .
$$
This is precisely the statement of the theorem.
\qed

\section{Implementation for twisted crossed products}
\label{Sec-CrossedProd}

\subsection{Twisted crossed product algebras by $\mathbb Z^k$}

Twisted C$^\ast$-dynamical systems and their covariant representations were introduced and studied in \cite{BusbyTAMS1970vf}, and the associated twisted crossed product algebras were further studied in \cite{PackerMPCPS1989fs}. Most facts and arguments are direct generalizations from crossed products without a twist, as treated in \cite{Ped}. This work will be dealing exclusively with twisted crossed products by $\ZM^k$, in which case the definitions simplify somewhat.

\begin{definition}[Twisted crossed product, \cite{NeshveyevBook2006jf} pp. 121]\label{Def-TwistedCP} Let $\mathcal B$ be a C$^\ast$-algebra, $\theta$ a $k \times k$ anti-symmetric matrix with real entries and $\xi$ a $\ast$-action of $\ZM^k$  on  $\mathcal B$,
$$
\xi\; :\; \ZM^k \rightarrow {\rm Aut}(\mathcal B)
\;, 
\qquad \xi_x \xi_y \;=\; \xi_{x+y}
\;, 
\qquad \xi_0 \;=\; {\rm id}
\;.
$$
The full twisted crossed product algebra $\mathcal A = \mathcal B \rtimes_\xi^\theta \ZM^k$ is the universal unital C$^\ast$-algebra generated by a copy of $\mathcal B$ and by the unitary elements $u_x$, $x \in \ZM^k$, satisfying the commutation relations:
\begin{equation}\label{Eq-CommRel}
u_x u_y \;=\; e^{\imath (x,\theta y)} u_{x+y}
\;, 
\qquad  
b u_x \;=\; u_x \xi_{-x}(b)
\;, 
\qquad 
b \in \mathcal B
\;, 
\qquad x,y \in \ZM^k
\;.
\end{equation} 
Throughout, $(\,\cdot\,,\,\cdot\,)$ represents the Euclidean scalar product.
\end{definition}

\begin{remark}\label{Re-UnitG} {\rm Since any element of a C$^\ast$-algebra can be written as the sum of four unitary elements \cite{BlackadarBook2006vh}, the crossed product algebra is generated by the group of unitary elements of $\mathcal B$, together with the unitaries $u_{e_j}$, where $e_j$, $j = 1, \ldots, k$, generate $\ZM^k$.}
\hfill $\diamond$
\end{remark}

\begin{proposition}[Alternative characterization, \cite{DavidsonBook1996bv} pp.~216]
\label{prop-AltChar}
Let $\mathcal A_0$ be the $\ast$-algebra of polynomials
$$
a\;=\;\sum_{q\in \mathbb Z^k} u_q \, b_q 
\; , 
\qquad \mbox{\rm with }\;b_q=0 \;\; {\rm if} \;\; |q| >R  
\;,
$$
for some $R<\infty$, which satisfy the algebraic operations
%
\begin{equation}\label{Eq-AlgOper}
a' a 
\;=\; 
\sum_{q \in \ZM^k} u_q \Big [\sum_{p \in \ZM^k} e^{\imath (q,\theta p)} \xi_{-p} (b'_{q-p}) b_{p}  \Big ]
\;, 
\qquad
a^\ast \;=\; \sum_{q \in \mathbb Z^k} u_q \, \xi_q (b_{-q}^\ast) 
\;.
\end{equation}
Then the twisted crossed product algebra $\mathcal A = \mathcal B \rtimes_\xi^\theta \ZM^k$ is obtained by closing $\mathcal A_0$ in the universal norm
$$
\|a\| 
\;=\; 
\sup_\pi \| \pi(a)\|_{\mathbb B}
\;,
$$
where the supremum is taken over all $\ast$-representations $\pi$ of $\mathcal A_0$ on separable Hilbert spaces. 
\end{proposition}

\begin{remark}{\rm Opposite to some of the standard textbooks, we chose here to place the unitary elements $u_q$ in front of the coefficients $b_q$, because the Hilbert C$^\ast$-modules on which $\Aa$ will be represented are right $\Bb$-modules.}
\hfill $\diamond$
\end{remark}

\begin{definition}[\cite{DavidsonBook1996bv}, pp.~222] 
The Fourier calculus on $\Aa =\Bb \rtimes_\xi^\theta \ZM^k$ is associated to the continuous group of automorphisms $\big\{\rho_\lambda\big \}_{\lambda \in \TM^k}$ on $\Aa$, often called the torus action, given by 
\begin{equation}\label{Eq-Rho}
\rho_ \lambda(u_q) 
\;=\; 
\lambda^q u_q, \quad \rho_\lambda(b) =b
\;, 
\qquad 
\lambda^q\;=\;
\lambda_1^{q_1} \cdots \lambda_k^{q_k}
\;.
\end{equation}
The Fourier coefficients of $a \in \mathcal B  \rtimes_\xi^\theta \mathbb Z^k$ are defined by the Riemann integral
\begin{equation}\label{Eq-FourierCoeff}
\Phi_q(a) 
\;=\; 
\int_{\TM^k}  \mu(d\lambda) \; \rho_\lambda (u_q^{-1} \, a) 
\;=\;
u_q^{-1} \int_{\TM^k} \mu(d\lambda) \; \lambda^{-q}\rho_\lambda (a)  
\;\in\; \mathcal B
\;,
\end{equation}
where $\mu$ is the Haar measure on the $k$-dimensional torus $\TM^k$. In particular, $\Phi_0$ provides a continuous, positive and faithful expectation from $\mathcal A$ to $\mathcal B$.
\end{definition}

\begin{proposition}[\cite{DavidsonBook1996bv}, pg.~223]\label{Feher} 
For any $a\in \mathcal A = \mathcal B \rtimes_\xi^\theta \ZM^k$, the Ces\`{a}ro sums:
\begin{equation}\label{Eq-Feher}
a_N
\;=\;
\sum_{q_1=-N}^N \ldots \sum_{q_k=-N}^N \prod_{j=1}^k \left (1-\frac{|q_j|}{N+1} \right ) u_ q\, \Phi_q(a)
\end{equation}
converge in norm to $a \in \Aa$ as $N \rightarrow \infty$. As a consequence, all the elements of $\Aa$ (and not just those from $\mathcal A_0$) can be represented as a Fourier series
$$
a
\;=\;
\sum_{q \in \mathbb Z^k} u_q \, b_q \; ,
$$
where $b_q \in \mathcal B$ should be interpreted as the Fourier coefficients of $a$ and the infinite sum as the limit of \eqref{Eq-Feher}.
\end{proposition}

\begin{definition}[Canonical trace, \cite{DavidsonBook1996bv}]\label{Def-CanTrace} 
Assume that there is a continuous $\xi$-invariant trace $\Tt_\Bb$ on $\Bb$. Then the Fourier calculus generates a canonical trace on $\mathcal A = \mathcal B \rtimes_\xi^\theta \ZM^k$. Indeed, the trace $\Tt_\Bb$ can be canonically promoted to a continuous trace on $\Aa$ by
$$
\Tt ( a ) 
\;=\; 
\Tt_\Bb \big (\Phi_0(a) \big ) 
\;=\; 
\Tt_\Bb (b_0)
\;.
$$
\end{definition}

\begin{proposition}
\label{Pro-TBProp} 
The following relations hold true for the trace $\Tt_\Bb$:
\begin{enumerate}[\rm (i)]

\item $\Tt_\Bb(b^\ast) = \Tt_\Bb(b)^\ast \, , \quad \Tt_\Bb(|b|) = \Tt_\Bb(|b^\ast|)$ and $ \Tt_\Bb(|\xi_x(b)|) = \Tt_\Bb(|b|)$

\item $\Tt_\Bb(|b'b|) \leq \|b'\| \, \Tt_\Bb(|b|) $ and $\big | \Tt_\Bb(b) \big | \leq \Tt_\Bb(|b|) $

\item $\Tt_\Bb(|bb'|^s)^\frac{1}{s} \leq \Tt_\Bb(|b|^{r})^\frac{1}{r} \, \Tt_\Bb(|b|^{p})^\frac{1}{p}$ for $\frac{1}{s} =\frac{1}{r} + \frac{1}{p}$ with $s,r,p\geq 1 $

\item $\Tt_\Bb(|b_1 b_2\cdots b_k|) \leq \prod_{i=1}^k \Tt_\Bb(|b_i|^{k})^\frac{1}{k} $

\item $\Tt_\Bb(|b +b'|^s)^\frac{1}{s} \leq \Tt_\Bb(|b|^s)^\frac{1}{s} + \Tt_\Bb(|b'|^s)^\frac{1}{s} $ for $ s \in [1,\infty)$

\end{enumerate}

\end{proposition} 

\proof (i) The first identity follows once $b$ is decomposed into its real and imaginary parts and the fact is used that $\Tt_\Bb$ is a real number when applied on self-adjoint operators. For the second identity, recall that $|b| = \sqrt{b^* b}$ and $|b^\ast|=\sqrt{bb^\ast }$. The continuous functional calculus can be approximated in norm by polynomial functional calculus. Then the statement follows because $\Tt_\Bb\big ( (bb^\ast)^n \big ) = \Tt_\Bb\big ( (b^\ast b)^n \big )$ due to the cyclic property of the trace. As for the third identity, note that the automorphisms commute with the continuous functional calculus and, as a consequence, $|\xi_x(b)| = \xi_x(|b|)$. The statement then follows from the invariance of $\Tt_\Bb$ w.r.t. to the action $\xi$. (ii) For the first inequality, note that $|b'b|^2 = b^\ast |b'|^2 b \leq \|b'\|^2 |b|^2$, which automatically implies $|b'b| \leq \|b'\| |b|$ for any $b,b' \in \Bb$ ({\it cf.} \cite{BlackadarBook2006vh} pp. 66). Since $\Tt_\Bb$ is a positive map, indeed $\Tt_\Bb(|b'b|) \leq \|b'\|\, \Tt_\Bb(|b|)$. For the second inequality, let $e^{\imath \phi}$ be the complex phase of $\Tt_\Bb(b)$. Then $\big |\Tt_\Bb(b)\big |=\Tt_\Bb(e^{\imath \phi}b)$ and the statement reduces to $\Tt_\Bb(b_\phi ) \leq \Tt_\Bb(|b_\phi|)$, with $b_\phi = e^{-\imath \phi}b$. If $b_\phi = b^{(r)}_\phi +\imath  b^{(i)}_\phi$ is the standard decomposition in the real and imaginary parts, then $\Tt_\Bb(b_\phi ) = \Tt_\Bb \big (b^{(r)}_\phi \big )$ and the statement reduced to $\Tt_\Bb \big (b^{(r)}_\phi \big ) \leq \Tt_\Bb(|b_\phi|)$. The latter follows from 
$$
\Tt_\Bb \big (b^{(r)}_\phi \big )
\;=\; 
\tfrac{1}{2}\,
\Tt_\Bb(b_\phi + b_\phi^\ast) 
\;\leq\; 
\tfrac{1}{2}\,
\Tt_\Bb(|b_\phi + b_\phi^\ast|) 
\;\leq\; 
\tfrac{1}{2}\, \big (\Tt_\Bb(|b_\phi|) + \Tt(|b_\phi^\ast|) \big ) 
\;=\; 
\Tt_\Bb(|b_\phi|)
\;,
$$
where we have made use of point (v) and (i), see below. Items (iii)-(v)  By imbedding the C$^\ast$-algebra into its weak von Neumann closure, the statements become standard and the proofs can be found, for example, in \cite{SegalAM1953bv}. 
\qed

\begin{definition}[Canonical derivations, \cite{DavidsonBook1996bv}]\label{Def-CanDer} The Fourier calculus generates a set of un-bounded derivations $\partial=(\partial_1, \ldots, \partial_k)$ on $\mathcal A = \mathcal B \rtimes_\xi^\theta \ZM^k$ as the generators of the $k$-parameter group of automorphisms $\lambda \in \TM^k\mapsto \rho_\lambda$. Explicitly, let $\mathcal C^1(\mathcal A)$ be the set of elements in $\Aa = \mathcal B \rtimes_\xi^\theta \ZM^k$ for which $\rho_\lambda(a)$ is first order differentiable in $\lambda$. On this space, the derivations act as
$$
\partial_i a 
\;=\;-
 \imath \sum_{q \in \mathbb Z^k} q_i u_q\, b_q \;, 
 \qquad i =1,\ldots,k
 \;.
$$
\end{definition}

\begin{definition}[Smooth sub-algebra, \cite{RennieKTh2003vj}]\label{Def-SmoothSubAlg} Let $\mathscr A$ be the set of elements from $\mathcal A$ for which $\rho_\lambda(a)$ are infinitely differentiable w.r.t. $\lambda$. This set is closed under the algebraic operations hence it is a $\ast$-subalgebra of $\mathcal A$. Furthermore, when endowed with the topology induced by the seminorms 
$$
\|a \|_\alpha \;=\; \|\partial^\alpha a \|\;, 
\qquad 
\partial^\alpha = \partial_1^{\alpha_1} \cdots \partial_k^{\alpha_k}
\;,\;\;
\alpha \;=\; 
(\alpha_1, \ldots \alpha_k)
\;,
$$ 
it becomes a dense Fr\'{e}chet sub-algebra of $\Aa$ which is stable under holomorphic calculus. As before, this sub-algebra will be called the smooth sub-algebra and will be denoted by $\mathscr A$.
\end{definition}

\subsection{Representation on the standard Hilbert $C^\ast$-module}

Let $\ell^2(\ZM^k)$ be the Hilbert space of square summable sequences over the lattice $\ZM^k$. Let us denote its standard basis by $(\delta_x)_{x \in \ZM^k}$, and introduce the shifts $S_y \delta_x = \delta_{x+y}$, as well as the magnetic shifts 
$$
U_y \delta_x
\; =\; 
e^{\imath (y,\theta x)}\delta_{x+y}\;, 
\qquad U_y U_x \;=\; e^{\imath 2(y,\theta x)}\,U_x U_y 
\;=\; e^{\imath (y,\theta x)}\,U_{x+y}
\;.
$$
Lastly let us use the usual bra-ket notation for the rank-one operators on  $\ell^2(\ZM^k)$:
$$
\Theta_{\delta_x,\delta_y} 
\;=\; 
|x\rangle \langle y|
\;.
$$

\begin{proposition}
\label{Pr-StandardVec1}
Any $\psi \in \mathcal H_\mathcal B$ can be uniquely represented as
$$
\psi 
\;=\; 
\sum_{x \in \ZM^k} \delta_x \otimes_\CM b_x
\;, 
\qquad 
\lim_{R \rightarrow \infty}\,  
\sum_{|x|<R} b_x^\ast b_x \in \mathcal B
\;,
$$
with coefficients given by $b_x = \langle \delta_x \otimes_\CM 1, \psi \rangle$.  In particular, the Hilbert $\Bb$-module $\mathcal H_\mathcal B$ is countably generated and $(\delta_x \otimes_\CM )_{x \in \ZM^k}$ is a set of generators. 
\end{proposition}

\proof Using a bijection from $\NM$ to $\ZM^k$, an isomorphism between the standard presentation of $\mathcal H_\mathcal B$ of Example~\ref{Ex-StandardHM1} and the presentation from above is readily established. 
\qed

\begin{remark}{\rm If $\psi = (b_x)_{x \in \ZM^k} \in \mathcal H_\mathcal B$ (in the presentation of Example~\ref{Ex-StandardHM1}), then
$$
\langle \psi,\psi \rangle 
\;=\; 
\sum_{x \in \ZM^k} b_x^\ast b_x 
\;=\; 
\sum_{x \in \ZM^k} \langle \psi,\delta_x \otimes_\CM 1 \rangle \langle \delta_x \otimes_\CM 1,\psi \rangle,
$$
with the sums converging in $\Bb$. This shows that $(\delta_x \otimes_\CM 1)_{x \in \ZM^k}$ is a standard normalized tight frame for $\mathcal H_\mathcal B$ and, according to \cite{FrankJOT2002re}, every element of $\Hh_\mathcal B$ accepts a decomposition
$$
\psi 
\;=\; 
\sum_{x \in \ZM^k} (\delta_x \otimes_\CM 1 ) \langle \delta_x \otimes_\CM 1 , \psi \rangle 
\;=\; 
\sum_{x \in \ZM^k} \delta_x \otimes_\CM b_x
\;,
$$
where the sums converge in the norm topology of $\mathcal H_\mathcal B$. This provides an alternative proof of the above statement.
}
\hfill $\diamond$
\end{remark}

\begin{remark} {\rm The following representation of the crossed-product algebra is the heart of our application. It is certainly related to the stabilizing representation considered in \cite{PackerMPCPS1989fs}, as well as the representation used in \cite{BKR}.
}
\hfill $\diamond$
\end{remark}

\begin{proposition}[\cite{Pro3}]\label{Pro-GenRep} 
Let $\Hh_\Bb = \Hh \otimes_\CM \Bb$ be the standard ungraded Hilbert $\Bb$-module and choose specifically $\Hh = \ell^2(\ZM^k)$ as the separable Hilbert space. Then the map
\begin{equation}
\label{Eq-OurRep}
\Bb \rtimes_\xi^\theta \ZM^k \;\ni\; 
a\;=\;
\sum_{q \in \ZM^k} u_q \, b_q 
\;\mapsto\; 
\eta (a) 
\;=\; 
\sum_{q,x \in \ZM^k} U_q |x \rangle \langle x| \otimes_\CM \xi_{-x} (b_q) 
\;\in\; 
\BM(\Hh_\Bb)
\end{equation}
is well-defined and it provides a faithful morphism of C$^\ast$-algebras. Above, $\xi_{-x} (b_q)$ are seen as operators acting by left multiplications on $\Bb$.
\end{proposition}

\proof Let us consider a generic element from $\Hh_\Bb$
$$
\psi 
\;=\; 
\sum_{x \in \ZM^k} \delta_x \otimes_\CM c_x\;, 
\qquad 
\lim_{R \rightarrow \infty}\, \sum_{|x| < R} c_x^\ast c_x \;\in\; \Bb\;,
$$
and let us write out the action of $\eta(a)$ on $\psi$ and its inner product:
\begin{equation}
\label{Eq-Action1}
\eta(a) \psi 
\;= \;
\sum_{q,x \in \ZM^k} U_q \,\delta_x \otimes_\CM \xi_{-x}(b_q) c_x
\;, 
\qquad 
a \;= \;\sum_{q  \in \ZM^k} b_q u_q \;\in\; \Aa
\;,
\end{equation}
and 
\begin{equation}
\label{Eq-NormAction1}
\big \langle \eta(a) \psi , \eta(a) \psi \big \rangle 
\;=\; 
\sum_{x  \in \ZM^k} c_x^\ast \xi_{-x}\Big ( \sum_{q  \in \ZM^k} b_q^\ast b_q \Big ) c_x 
\;=\; 
\sum_x c_x^\ast \xi_{-x}\Big ( \Phi_0(a^\ast a ) \Big ) c_x
\;.
\end{equation}
Since $\xi$ is a $\ast$-action, all operators $\xi_{-x}\Big ( \Phi_0(a^\ast a ) \Big )$ are  positive elements in $\Bb$ and thus
$$
\big \langle \eta(a) \psi , \eta(a) \psi \big \rangle 
\;\leq\; 
\big \|\Phi_0(a^\ast a ) \big \| \sum_{x \in \ZM^k} c_x^\ast  c_x 
\;=\;
 \big \|\Phi_0(a^\ast a )  \big \| \big \langle \psi , \psi \big \rangle
 \;.
$$
Hence, $\eta(a)\psi \in \Hh_\Bb$ for all $a \in \Aa$ and $\psi \in \Hh_\Bb$ and one can check explicitly that $\eta(a)^\ast = \eta(a^\ast)$. The conclusion is that $\eta(a)$ is indeed an adjointable operator for all $a \in \Aa$. Next let us focus on the generators of $\Aa$:
$$
\eta(u_{e_j}) 
\;=\; 
U_{e_j} \otimes 1\;, 
\qquad 
\eta(b) \;=\; 
\sum_{x \in \ZM^k} |x \rangle \langle x| \otimes \xi_{-x}(b)
\;,
$$
with $b$ unitary, {\it cf.} Remark~\ref{Re-UnitG}. From \eqref{Eq-Action1} and \eqref{Eq-NormAction1} one can see that the representations of the generators are indeed unitary operators from $\BM(\mathcal H_\mathcal B)$. Furthermore
$$
\eta(u_{e_i}) \eta(u_{e_j}) 
\;=\; 
U_{e_i} U_{e_j} \otimes 1
\;=\;
 e^{\imath \theta_{i,j}} U_{e_i+e_j} \otimes 1 
 \;=\; 
 e^{\imath \theta_{i,j}} \eta(u_{e_i + e_j})
\;,
$$
and
$$
\eta(b) \eta(u_{e_j}) 
\;=\; 
\sum_{x \in \ZM^k} |x \rangle \langle x|U_{e_j} \otimes \xi_{-x}(b) 
\;=\; 
\sum_{x \in \ZM^k} U_{e_j}|x - e_j \rangle \langle x -e_j| \otimes \xi_{-x}(b)
\;.
$$
After a change of summation variable $x \rightarrow x+e_j$, one obtains
$$
\eta(b) \eta(u_{e_j}) 
\;=\; 
(U_{e_j}\otimes 1) \sum_{x \in \ZM^k} |x  \rangle \langle x | \otimes \xi_{-x-e_j}(b) 
\;=\; 
\eta(u_{e_j}) \eta \big ( \xi_{-e_j}(b) \big )
\;.
$$
The commutation relations \eqref{Eq-CommRel} are thus indeed satisfied.
\qed

\subsection{A generalized Connes-Chern character: the even case}\label{Sec-EvenConnesChern}

Throughout this section, $k$ is considered even, and, as above, $\Bb$ is unital, separable and trivially graded. All this then also holds for the crossed product $\Aa = \Bb \rtimes_\xi^\theta \ZM^k$.  Let us next consider the Clifford algebra $\CM_k$ with the grading described in Example~\ref{ex-CliffMatIso}. We will freely use the isomorphism $\CM_k \cong \CM(2^\frac{k}{2})$ and treat the elements of the Clifford algebras as matrices with complex entries. This provides a trace on $\CM_k$, given by the ordinary trace over $\CM(2^\frac{k}{2})$ and denoted by ${\rm tr}$. The relevant Hilbert $\Bb$-module for our calculation is 
$$
E_\Bb 
\;=\; 
\CM^{2^\frac{k}{2}} \hat \otimes_\CM \Hh_\Bb
\;,
$$
with the grading operator
$$
S \;=\; \gamma_0 \hat{\otimes}_\CM {\rm id}
\;.
$$
Note that this is a graded Hilbert C$^\ast$-module over the ungraded C$^\ast$-algebra $\Bb$. 

\begin{proposition}\label{Pr-TraceClass} Let $\delta_\alpha$, $\alpha =1,\ldots,2^\frac{k}{2}$, be a basis for $\CM^{2^\frac{k}{2}}$ and let us introduce the notation $\Theta_{\delta_\alpha,\delta_\beta} = |\alpha \rangle \langle \beta|$. Then
\begin{enumerate}[\rm (i)]

\item Any $\psi \in E_\mathcal B$ can be uniquely represented as
$$
\psi 
\;=\; 
\sum_{\alpha} \sum_{x \in \ZM^k} \delta_\alpha \hat \otimes_\CM \delta_x \otimes_\CM b_{\alpha, x}
\;, 
\qquad 
\lim_{R \rightarrow \infty}  \sum_{\alpha} \sum_{|x|<R} b_{\alpha, x}^\ast b_{\alpha, x} \in \Bb\;.
$$
Here $b_{\alpha,x} = \langle \delta_\alpha \hat \otimes_\CM \delta_x \otimes_\CM 1, \psi \rangle$. 

\item Any $T \in \KM(E_\mathcal B)$ can be uniquely represented as
$$
T 
\;=\; 
\sum_{\alpha,\beta} \sum_{x, y \in \ZM^k} |\alpha \rangle \langle \beta | \hat \otimes_\CM |x\rangle \langle y | \otimes_\CM b_{\alpha, x;\beta, y}
\;, 
$$
where the coefficients are $b_{\alpha,x;\beta,y} = \langle \delta_\alpha \hat \otimes_\CM \delta_x \otimes_\CM 1 , T(\delta_\beta \hat \otimes_\CM \delta_y \otimes_\CM 1) \rangle$ and satsify
$$
\lim_{R \rightarrow \infty} \sum_{\alpha,\beta} \sum_{|x|, |y| < R} |\alpha \rangle \langle \beta | \hat \otimes_\CM |x\rangle \langle y | \otimes_\CM b_{\alpha, x;\beta, y} \;\in\; \BM (E_\mathcal B)
\;,
$$
with the limit in norm topology of $\BM (E_\Bb)$.

\item Any $T \in \BM(E_\mathcal B)$ can be uniquely represented as
\begin{equation}
\label{eq-Tdecomp}
T 
\;=\; 
\sum_{\alpha,\beta} \sum_{x, y \in \ZM^k} |\alpha \rangle \langle \beta | \hat \otimes_\CM |x\rangle \langle y | \otimes_\CM b_{\alpha, x;\beta, y}
\;, 
\end{equation}
where the coefficients are such that
$$
\lim_{R \rightarrow \infty} \sum_{\alpha,\beta} \sum_{|x|, |y| < R} |\alpha \rangle \langle \beta | \hat \otimes_\CM |x\rangle \langle y | \otimes_\CM b_{\alpha, x;\beta, y} 
\;\in\; \BM (E_\Bb),
$$
with limit in the strict topology of $\BM (E_\Bb)$.

\item The unique extension of the trace $\mathcal T_\mathcal B$ on $\Bb$ to $E_\mathcal B$ is given by
$$
{\rm Tr}(T) 
\;=\; 
\sum_\alpha \sum_{x \in \ZM^k} \mathcal T_\mathcal B(b_{\alpha,x; \alpha,x})
\;, 
\qquad T \in \SM_1(E_\Bb)
\;.
$$

\item Let $T \in \BM(E_\mathcal B)$ be decomposed as in \eqref{eq-Tdecomp}. If
$$
\sum_{\alpha,\beta} \sum_{x,y \in \ZM^k} \Tt_\Bb \big ( |b_{\beta,y;\alpha,x}|\big ) 
\;< \;\infty
\;,
$$
then $T \in \SM_1(E_\mathcal B)$. Since $\Tt_\Bb(|b|) \leq \|b\|$ for any $b\in \Bb$, one also has
$$
\sum_{\alpha,\beta} \sum_{x,y \in \ZM^k} \|b_{\beta,y;\alpha,x}\| \;<\; \infty 
\quad \Longrightarrow \quad T \in \SM_1(E_\mathcal B)
\;.
$$

\end{enumerate}

\end{proposition}

\proof (i) follows from Proposition~\ref{Pr-StandardVec1}. (ii) Since any compact operator can be approximated in norm by a finite sum of rank-one operators it is enough to establish the statement just for the latter. Using the representation of the vectors from (i), one can see that any rank-one operator takes the form
\begin{equation}
\label{Eq-Theta1}
\Theta_{\psi',\psi} 
\;=\; 
\sum_{\alpha,\beta}\sum_{x,y \in \ZM^k} |\alpha \rangle \langle \beta | \hat \otimes_\CM |x \rangle \langle y | \otimes_\CM b'_{\alpha,x} b^\ast_{\beta,y}
\;,
\end{equation}
hence they are indeed of the form stated. Clearly, finite linear combinations preserve this form and can be approximated in norm by sums. (iii) It follows from (ii) since the unit ball of $\KM(E_\Bb)$ is dense in the unit ball of $\BM(E_\Bb)$, when the strict topology is used ({\it cf.} \cite{LanceBook1995vc} pp.~11). The coefficients of the expansion are give by the same formula as in (ii). (iv) It is enough to verify the formula on rank-one operators. From its defining property,
$$
{\rm Tr}(\Theta_{\psi',\psi}) 
\;=\; 
\Tt_\Bb (\langle \psi, \psi' \rangle) 
\;=\; 
\sum_\alpha \sum_{x \in \ZM^k}\Tt_\Bb(b^\ast_{\alpha,x} b'_{\alpha,x})
\;,
$$
where we could inter-change the trace and the summation because $\Tt_\Bb$ is continuous and the sum $\sum_\alpha \sum_x b^\ast_{\alpha,x} b'_{\alpha,x}$ converges in $\Bb$. Using \eqref{Eq-Theta1} and the cyclic property of the trace, the affirmation follows.  (v) The difficulty here is to deal with the absolute value because $T$ may not have a polar decomposition. Recall also that $\BM (E_\mathcal B)$ is stable only under the continuous functional calculus. Let us assume for the beginning that $T$ is self-adjoint, in which case the absolute value can be computed by applying the continuous function $|t|=t \, \sgn(t)$ on the operator $T$. Using the approximation $\sgn_\epsilon(t) = \tanh(t/\epsilon)$, one can write $|T| = \lim\limits_{\epsilon \downarrow 0} \sgn_\epsilon (T) \, T$, with the limit taken w.r.t. the C$^\ast$-norm of $\BM (E_\Bb)$. The point of the last expression is that $\sgn_\epsilon(T)$ belongs to $\BM(E_\Bb)$, while $\sgn(T)$ does not, in general. Now, as $\sgn_\epsilon(T) \, T$ is an increasing sequence of positive operators as $\epsilon \downarrow 0$ and ${\rm Tr}$ is lower semi-continuous, 
$$
{\rm Tr}(|T|) 
\;=\; 
\lim_{\epsilon \downarrow 0} \,{\rm Tr}\big ( \sgn_\epsilon (T) \, T \big ) 
\;=\; 
\lim_{\epsilon \downarrow 0} \;\sum_{\alpha}\sum_{x \in \ZM^k} \Tt_\Bb \Big (\big ( \sgn_\epsilon (T) \, T \big )_{\alpha,x;\alpha,x} \Big )
\;.
$$
The representation from (ii) leads to
$$
\big ( \sgn_\epsilon (T) \, T \big )_{\alpha,x;\alpha,x} 
\;=\; 
\sum_\beta \sum_{y \in \ZM^k} \sgn_\epsilon (T)_{\alpha,x;\beta,y}  \, b_{\beta,y;\alpha,x}
$$
with
$$
\sgn_\epsilon (T)_{\alpha,x;\beta,y} 
\;=\; 
\big \langle \delta_\alpha \hat \otimes_\CM \delta_x \otimes_\CM 1 , \sgn_\epsilon (T) (\delta_\beta \hat \otimes_\CM \delta y \otimes_\CM 1) \big \rangle
\; ,
$$
satisfying
$$
\| \sgn_\epsilon (T)_{\alpha,x;\beta,y}\|_\Bb 
\;\leq \;
\|\sgn_\epsilon (T) \|_{\BM(E_\Bb)} 
\;\leq \;
1
\;.
$$
Thus, by applying points (iv) and (v) of Proposition~\ref{Pro-TBProp},
$$
\Big | \Tt_\Bb \Big ( \big ( \sgn_\epsilon (T) \, T \big )_{\alpha,x;\alpha,x} \Big ) \Big | 
\;\leq\; 
\sum_{\beta}\sum_{y \in \ZM^k} \big \|\ \sgn_\epsilon (T)_{\alpha,x;\beta,y} \big \| \, \Tt_\Bb\big ( |b_{\beta,y;\alpha,x}| \big )
\;,
$$
and hence
\begin{equation}\label{Eq-S1}
{\rm Tr}\big (|T| \big ) 
\;\leq \;
\sum_{\alpha,\beta} \sum_{x,y\in \ZM^k} \Tt_\Bb \big ( |b_{\beta,y;\alpha,x}|\big )
\;.
\end{equation}
A generic operator can be decomposed $T= \tfrac{1}{2}(T+T^\ast) +\imath \tfrac{1}{2\imath}(T-T^\ast)$ into a sum of two self-adjoint onece, and for each summand \eqref{Eq-S1} applies:
$$
{\rm Tr}\big (|T\pm T^\ast| \big ) 
\;\leq\; 
\sum_{\alpha,\beta} \sum_{x,y\in \ZM^k} \Tt_\Bb \big ( |b_{\beta,y;\alpha,x} \pm b_{\beta,y;\alpha,x}^\ast |\big )
\;.
$$
Then point (viii) of Proposition~\ref{Pro-TBProp} assures us that, if  the hypothesis of (v) holds, then both the real and imaginary parts of $T$ belong to $\SM_1(E_\Bb)$, hence $T \in \SM_1(E_\Bb)$. 
\qed

\begin{definition} 
\label{def-DiracDef}
Let $X$ be the position operator over $\ell^2(\ZM^k)$, $X \delta_x= x \delta_x$. Then the shifted Dirac operator on $E_\mathcal B$ is defined as 
\begin{equation}\label{Eq-DiracOp}
D_{x_0} 
\;=\; 
\sum_{i=1}^k \gamma_i \hat \otimes_\CM (X_i + x_0) \otimes_\CM 1
\;=\; \gamma \hat \otimes_\CM (X + x_0) \otimes_\CM 1\;, 
\qquad x_0 \in (0,1)^k
\;,
\end{equation}
where $\gamma=(\gamma_1,\ldots,\gamma_k)$ and contraction over the indices is assumed in the second equality. The phase of $D_{x_0}$,
\begin{equation}\label{Eq-DiracPhase}
F_{x_0} 
\;=\;
\frac{D_{x_0}}{|D_{x_0}|} 
\;=\; 
\sum_{i=1}^k \gamma_i \hat \otimes_\CM \frac{(X_i + x_0)}{\sqrt{(X+x_0)^2}} \otimes_\CM 1
\;=\; 
\gamma \hat \otimes_\CM \widehat{X + x_0} \otimes_\CM 1
\;,
\end{equation}
is an element of $\mathbb B(E_\mathcal B)$. Throughout, we will employ the notation $\widehat{x}= x|x|^{-1}$ for invertible operators as well as ordinary vectors.
\end{definition}

\begin{remark} {\rm The shift $x_0$ assures that $D_{x_0}$ is invertible so that  its phase is well-defined. But the considerations leading to the inclusion of the shift $x_0$ go well beyond that. In particular, the average over $x_0$ is absolutely essential for the key identities discovered in \cite{ProdanJPA2013hg,PS} and used later on.} 
\hfill $\diamond$
\end{remark}

\begin{theorem}\label{Th-EvenKKCycle} The triples 
$$
\mathcal E^0_{x_0} 
\;=\; 
\Big (E_\Bb = \CM^{2^\frac{k}{2}} \hat \otimes_\CM \Hh_\Bb, \pi = 1 \hat \otimes_\CM \eta, F_{x_0} = \gamma \hat \otimes_\CM \widehat{X+x_0} \otimes_\CM 1 \Big )
\;, \qquad x_0 \in (0,1)^k,
$$ 
define a field of Kasparov $(\Aa,\Bb)$-cycles.
\end{theorem}

\proof We will cover one by one the properties (i)-(v) listed in Definition~\ref{Def-KasparovM}.

\vspace{0.1cm}

\noindent (i) follows from $\gamma_0\gamma=-\gamma\gamma_0$. (ii) The algebra $\Aa$ is ungraded, hence this property translates into  $S \pi(a) S = \pi(a)$, and this follows from
$$
S \pi(a) S 
\;=\; 
(\gamma_0 \hat \otimes_\CM {\rm id}) \big (1 \hat \otimes \eta(a) \big )  (\gamma_0 \hat \otimes_\CM {\rm id}) 
\;=\;  
1 \hat \otimes_\CM \eta(a) 
\;=\; 
\pi(a)
\;.
$$

\vspace{0.1cm}

\noindent (iii) This property becomes
$$
\big [ \gamma \hat \otimes_\CM \widehat{X + x_0} \otimes_\CM 1, 1 \hat \otimes_\CM \eta(a) \big ] 
\;=\; 
\gamma \hat \otimes_\CM \big [\widehat{X + x_0} \otimes_\CM 1, \eta(a) \big ]  \in \KM (E_\mathcal B)
\;, 
\qquad a \in \Aa
\;.
$$
Since $\KM(E_\mathcal B)$ is a closed double-sided ideal, it is enough to verify the statement just for the generators. We have $\big [\widehat{X + x_0} \otimes_\CM 1, \pi(b)\big ]=0$ for any $b\in \mathcal B$, and
$$
\big [\widehat{X + x_0} \otimes_\CM 1, \eta(u_q)\big ] 
\; =\; 
\big [\widehat{X + x_0}, U_q \big ] \otimes_\CM 1  
\;=\; 
\Big (\widehat{ X + x_0} - \widehat{X - q + x_0} \Big ) U_q\otimes_\CM 1
\;.
$$
Since
\begin{equation}\label{Eq-Asympt}
\widehat{ x + x_0} - \widehat{x - q+x_0} 
\;\sim\; 
|x+x_0|^{-1} \big ( q +(\widehat{x+x_0},q) \, \widehat{x+x_0} \big )
\end{equation}
as $|x| \rightarrow \infty$, it follows that $\widehat{ X + x_0} - \widehat{X - q + x_0}$ is compact as an operator on $\ell^2(\ZM^k)$. As such, the commutator $ \big [\widehat{X + x_0}, \eta(u_q)\big ]$ can be approximated in the norm topology of $\BM(E_\mathcal B)$ by finite rank operators, and is hence a compact operator over $E_\mathcal B$.

\vspace{0.1cm}

\noindent (iv) and (v)  are satisfied because $F_{x_0}^2 = 1$ and $F_{x_0}^\ast = F_{x_0}$. 
\qed

\begin{theorem}
\label{Th-EvenSummable} 
The Kasparov cycles defined in Theorem~\ref{Th-EvenKKCycle} are $n$-summable over the smooth sub-algebra $\mathscr A$ for any $n \geq k+1$, that is,
$$
\prod_{i=1}^n \big [F_{x_0},\pi(a_i) \big ] \;\in\; \SM_1(E_\mathcal B)
$$
for any $a_i \in \mathscr A$, $i=1,\ldots,n$. 
\end{theorem}

\proof 
We will show that
\begin{equation}\label{Eq-Estimate1}
Y
\;=\;
\sum_{\beta}\sum_{y \in \ZM^k} \Big\| \Big ( \prod_{i=1}^n \big [F_{x_0}, \pi(a_i) \big ] \Big )_{\alpha,x;\beta,y}\Big \| 
\;\leq \;
C\;  |\bm x + \bm x_0|^{-n}\;,
\end{equation}
for some constant $C<\infty$. This implies the claim because the criterion (v) of Proposition~\ref{Pr-TraceClass} applies to $\prod_{i=1}^n \big [F_{x_0},\pi(a_i) \big ]$, as long as $n \geq k+1$. To start out, let us write the $a_i \in \mathscr A$ as in Proposition~\ref{Feher}:
$$
a_i \;=\; \sum_{q \in \ZM^k} u_q \, b^{(i)}_{q} \;\in\; \mathscr A\;, \qquad i=1,\ldots n\;.
$$ 
It is also useful to write each commutator explicitly
$$
\big [F_{x_0}, \pi(a)\big ] 
\;=\;  
\sum_{x_1,x_2 \in \mathbb Z^k} e^{\imath (x_1,\theta x_2)} \gamma \hat \otimes_\CM (\widehat{x_1+x_0} -\widehat{x_2+x_0})|x_1 \rangle \langle x_2 | \otimes_\CM \xi_{-x_2} (b_{x_1-x_2})
\;,
$$
where as above the summation over indices in $\gamma$ and the $x_1,x_2$ is suppressed. Then
\begin{align*}
\prod_{i=1}^n \big [F_{x_0},\pi(a_i) \big ]
= 
\!\!\!
\sum_{x_1,\ldots,x_{n+1} \in \mathbb Z^k} \prod_{i=1}^n e^{\imath (x_i,\theta x_{i+1})}
\gamma \hat \otimes_\CM (\widehat{x_i + x_0}  - \widehat{x_{i+1} +x_0}) |x_{i} \rangle \langle x_{i+1}| \otimes_\CM \xi_{-x_{i+1}}( b^{(i)}_{x_i - x_{i+1}})
\,,
\end{align*}
so that
\begin{align}
\Big(\prod_{i=1}^n \big [F_{x_0}, \pi(a_i) \big ] \Big)_{\alpha,x;\beta,y}  
 \;=\; &
\sum_{x_1,\ldots,x_{n+1} \in \mathbb Z^k} \delta_{x_1,x} \, \delta_{x_{n+1},y} \Big (\prod_{i=1}^n e^{\imath (x_i,\theta x_{i+1})} \Big ) 
\nonumber
\\
& \;\;\times 
\Big (\prod_{i=1}^n \big (\gamma, \, \widehat{x_i + x_0} - \widehat{x_{i+1} +x_0} \big) \Big )_{\alpha,\beta} 
\prod_{i=1}^n \xi_{-x_{i+1}}( b^{(i)}_{x_i - x_{i+1}}) \; .
\label{PP32}
\end{align}
Hence
\begin{align*}
Y 
\;\leq\; 
\sum_{x_1,\ldots,x_{n+1} \in \mathbb Z^k} \delta_{x_1,0}  \prod_{i=1}^n 
\,
\big |\widehat{x_i + x + x_0} - \widehat{x_{i+1} + x + x_0} \big |\; \big \| b^{(i)}_{x_i - x_{i+1}} \big \|
\;.
\end{align*} 
Due to the asymptotic behavior in \eqref{Eq-Asympt}, the supremum 
$$
S(y,y')
\;=\;
\sup \big \{|x| \big |\widehat{ y + x } - \widehat{y' + x} \big |\;:\; x \in \mathbb R^k  \big \}
$$
is finite and it has the scaling property $S(y, y') = s^{-1}S(s y, s y')$. By taking $s = (|y|+|y'|)^{-1}$, one obtains the upper bound
$$
S(y,y') 
\;\leq\; 
(|y|+|y'|)\, \sup\{S(x,x')\;:\;  |x|+|x'| = 1\}
\;.
$$
The conclusion is
$$
\big |\widehat{ y + x }\; - \;\widehat{y' + x} \big | 
\;\leq\; 
C\, |x|^{-1}\,(|y|+|y'|)
\;.
$$
Thus 
$$
Y 
\;\leq\; 
C \,
| x \,+ \, x_0|^{-n} \!\sum_{x_1,\ldots,x_{n+1} \in \mathbb Z^k} \delta_{x_1,0}\,  
\prod_{i=1}^n \,(|x_i|+|x_{i+1}|) 
\;\big \| b^{(i)}_{x_i - x_{i+1}} \big \|
\;.
$$
Let us now make the change of variables $y_i = x_{i+1}-x_i $, $i=1,\ldots,k$. Since $ x_1=0$,
$$
 x_{i+1} 
\;=\;
y_1 + \ldots + y_i
\;\;\; \Longrightarrow \;\;\;
|x_{i+1}|\; \leq\; \prod_{j=1}^n \big (1+| y_j| \big )
\;.
$$
This gives
\begin{align*}
Y & \;\leq\; C\, | x + x_0|^{-n} \sum_{y_1,\ldots,y_n \in \mathbb Z^k}  
\prod_{i=1}^n 2^n\big (1+|y_i| \big )^n \big \| b^{(i)}_{y_i} \big \| \\
 &\; = \;C\,2^n\, | x + x_0|^{-n} \prod_{i=1}^n \Big (\sum_{y_i \in \mathbb Z^k}  \big (1+|y_i|\big )^n \big \| b^{(i)}_{y_i} \big \| \Big ) \; .
\end{align*}
Since the elements $a_i$ belong to the smooth sub-algebra, the norms $\big \|b^{(i)}_x \big \|$ decay in $|x|$ faster than any inverse power, so the sums are finite  and \eqref{Eq-Estimate1} is proved, with another constant $C$. 
\qed

\vspace{.2cm}

\begin{remark} {\rm Now all the conditions of the general theory as developed in the previous chapter are established. Hence we can formulate at once the following statement.} 
\hfill $\diamond$
\end{remark}

\begin{corollary}\label{Co-EvenChernConnesCh} 
The multilinear map
\begin{equation}
\label{Eq-EvenChernConnesCh}
\uptau_k(a_0,\ldots, a_k) 
\;=\;  
\Gamma_k \, {\rm Tr}\Big ( S F_{x_0} \prod_{i=0}^k \big [F_{x_0}, \pi(a_i)\big ] \Big )  
\end{equation}
is a well-defined cyclic $(k+1)$-cocycle over the smooth sub-algebra $\mathscr A$. Its class in the periodic cyclic cohomology of $\mathscr A$ defines the generalized Connes-Chern character of the Kasparov cycle $\mathcal E^0_{x_0}$. It fulfills the following index theorem
\begin{equation}\label{Eq-Index0}
\Lambda_k \, \uptau_k(p,\ldots, p)
\;=\; 
 {\rm Tr} \Big ( {\rm Ind}\big (\pi(p) F_{x_0} \pi (p) \big ) \Big )
\;,
\end{equation}
for any projection $p \in \mathscr A$.
\end{corollary} 

\begin{theorem}[Even local formula] 
\label{Th-EvenLocalFormulaTh}
The generalized Connes-Chern characters of the Kasparov cycles $\mathcal E^0_{x_0}$ coincide and can be represented by a cyclic cocycle $\bar \uptau_k$ which accepts the following local formula,
\begin{equation}\label{Eq-LocalEven}
\bar \uptau_k( a_0,\ldots,a_k)
\;=\; 
\Delta_k \sum_{\rho \in S_k} (-1)^\rho \;\mathcal T \Big (a_0  \prod_{i=1}^k \partial_{\rho_i} a_i \Big )
\;, 
\qquad \Delta_k 
\;=\; 
(2\pi \imath)^{\frac{k}{2}}\;,
\end{equation}
for any $a_i \in \mathscr A$, $i=0,\ldots,k$. Here, $S_k$ denotes the permutation group of order $k$, and $(-1)^\rho$ the signum of $\rho\in S_k$.
\end{theorem}

\proof The first part of the affirmation follows from the fact  that $F_{x_0}$ is continuous of $x_0$ as an element of $\BM(E_\Bb)$ so that the cyclic cocycles $\tau_k$ defined for different $x_0$ are homotopic. Indeed, they therefore belong to the same cohomology class  \cite{Connes:1994wk}. As a consequence, for any $x<x'\in (0,1)$ the cyclic cocycle
\begin{equation}\label{Eq-AverageCycle}
\bar \uptau_k(a_0, \ldots, a_k) 
\;=\; 
(\mbox{Vol}_k[x,x']^k)^{-1}\,\int\limits_{[x,x']^k}dx_0 \, \uptau_k(a_0, \ldots, a_k)
\end{equation}
represents the same cohomology class. Further let us note that, since the summability estimates in Theorem~\ref{Th-EvenSummable}  are uniform w.r.t. $x_0$, the integration over $x_0$ can and will be extended over the whole unit cube $[0,1]^k$ which has unit volume.
Then we also write $\bar \uptau_k= \bar \uptau_k(a_0, \ldots, a_k) $ for sake of brevity.

\vspace{0.1cm}

For the local formula, let us start from \eqref{Eq-EvenChernConnesCh} and apply (iv) of Proposition~\ref{Pr-TraceClass}. By opening the first commutator, 
$$
\bar \uptau_k
\;=\;
\Gamma_k \int\limits_{[0,1]^k}dx_0\sum_\alpha \sum_{x \in \mathbb Z^k} \mathcal T_{\mathcal B} \Big ( \Big (S \big (\pi(a_0)  - F_{x_0} \pi(a_0) F_{x_0} \big ) \prod_{i=1}^{k} \big [F_{x_0},\pi(a_i)\big]\Big )_{\alpha,x;\alpha, x} \Big )
\;.
$$
We recall that the sum is absolute convergent due to Theorem~\ref{Th-EvenSummable}. Given \eqref{Eq-Id07}  and $SF_{x_0} = -F_{x_0}S$, as well as the fact that $k$ is even, one can push the two operators $F_{x_0}$ flanking $\pi(a_0)$ to the edges where they annihilate due to $F_{x_0}^2=1$. Thus
\begin{equation}\label{Eq-PartFedosov}
\bar \uptau_k 
\;=\;
2\, \Gamma_k \int\limits_{[0,1]^k} d x_0\sum_\alpha \sum_{x \in \ZM^k} \mathcal T_{\mathcal B} \Big ( \Big (S \pi(a_0) \prod_{i=1}^{k} \big [F_{x_0},\pi(a_i) \big ]\Big )_{\alpha,x;\alpha,x} \Big ) \; .
\end{equation}
Using \eqref{PP32} for $n=k$, 
\begin{align*}
\bar \uptau_k \; = \;
& 2\, \Gamma_k \int\limits_{[0,1]^k} d x_0  \sum_{x,x_1,\ldots,x_{k+1} \in \ZM^k} \delta_{x,x_{k+1}}e^{\imath (x,\theta x_1)} \Big ( \prod_{i=1}^k e^{\imath (x_i,\theta x_{i+1})} \Big ) \\
& \;\; 
\times \, \mathrm{tr} \Big ( \gamma_0 \prod_{i=1}^k \big ( \gamma,  \widehat{x_0 + x_i} - \widehat{x_0 + x_{i+1}} \big ) \Big )   \;
\mathcal T_{\mathcal B} \Big ( \xi_{-x_1}(b^{(0)}_{x-x_1})\prod_{i=1}^k  \xi_{-x_{i+1}} (b^{(i)}_{x_i - x_{i+1}}) \Big )
\;.
\end{align*}
We next re-define the summation variables as $x_i \rightarrow x_i + x$, $i=1,\ldots,k+1$. By using the invariance of $\mathcal T_\mathcal B$ under the automorphisms $\xi_x$, it follows
\begin{align*}
\bar \uptau_k  \;=\;
& 2\, \Gamma_k \int\limits_{[0,1]^k} d x_0  \sum_{x,x_1,\ldots,x_{k+1}  \in \ZM^k} \delta_{x_{k+1},0}\, \Big(\prod_{i=1}^{k-1} e^{\imath (x_i,\theta x_{i+1})} \Big ) \\
&\; \times \, {\rm tr} \Big ( \gamma_0 \prod_{i=1}^k \big ( \gamma,  \widehat{x_0 +x+ x_i} -  \widehat{x_0 + x+ x_{i+1}} \big ) \Big )  
\;\mathcal T_{\mathcal B} \Big ( \xi_{-x_1}(b^{(0)}_{-x_1})\prod_{i=1}^k  \xi_{-x_{i+1}} (b^{(i)}_{x_i - x_{i+1}}) \Big )
\;.
\end{align*}
Note that the variable $x$ is no longer present in the phase factors. We now combine the summation over $x$ and the integration over $x_0$ in one integration over the whole $\RM^k$:
\begin{align*}
\bar \uptau_k  \;= \;
& 2\, \Gamma_k \sum_{x,x_1,\ldots,x_{k+1}  \in \ZM^k} \delta_{x_{k+1},0}  \, \Big(\prod_{i=1}^{k} e^{\imath (x_i,\theta x_{i+1})} \Big ) \\
& \;\;\times \int\limits_{\RM^k} d x \, {\rm tr} \Big ( \gamma_0 \prod_{i=1}^k \big ( \gamma,  \widehat{x+ x_i} - \widehat{x+ x_{i+1}} \big ) \Big ) 
\;\mathcal T_{\mathcal B} \Big ( \xi_{-x_1}(b^{(0)}_{-x_1})\prod_{i=1}^k  \xi_{-x_{i+1}} (b^{(i)}_{x_i - x_{i+1}}) \Big )
\;.
\end{align*}
At this point we appeal to the identity discovered in \cite{ProdanJPA2013hg} (see also \cite{ProdanSpringer2016ds} for a detailed proof),
\begin{align*}
\delta_{x_{k+1},0}\int\limits_{\mathbb R^k} d x \ {\rm tr} 
\Big ( \gamma_0 \prod_{i=1}^k \big (\gamma, \widehat{x+x_i} -\widehat{x + x_{i+1}}\big) \Big ) 
&
\;=\; 
\tilde \Lambda_k \,\sum_{\rho \in S_k} (-1)^\rho \prod_{i=1}^k (x_i)_{\rho_i}
\\
& \;=\;
\delta_{x_{k+1},0}\,\tilde \Lambda_k \,\sum_{\rho \in S_k} (-1)^\rho \prod_{i=1}^k (x_i-x_{i+1})_{\rho_i}
\;,
\end{align*}
where $\tilde \Lambda_k = \frac{(2\pi\imath)^\frac{k}{2}}{(k/2)!}$ and the second equality is simply due to the anti-symmetrizing factor. Now let us note that $2\,\Gamma_k \tilde \Lambda_k = \imath^k \Delta_k$. Thus
\begin{align*}
\bar \tau_k  & 
\;= \;
\imath^k \Delta_k \sum_{\rho \in S_k} (-1)^\rho \sum_{x_1,\ldots, x_{k+1} \in \ZM^k}  \Tt_\Bb \Big ( b^{(0)}_{-x_1} \prod_{i=1}^k e^{\imath (x_i,\theta x_{i+1})}  \xi_{- x_{i+1}}\big ( (x_i-x_{i+1})_{\rho_i}b^{(i)}_{x_i - x_{i+1}}\big )\delta_{x_{k+1},0} \Big ) \\
& 
\;=\; \Delta_k \sum_{\rho \in S_k} (-1)^\rho  
\; \Tt_\Bb \Big (\sum_{x_1,\ldots, x_{k+1} \in \ZM^k} b^{(0)}_{-x_1} \prod_{i=1}^k e^{\imath (x_i,\theta x_{i+1})}  \xi_{- x_{i+1}}\big ( \partial_{\rho_i}b^{(i)}_{x_i - x_{i+1}}\big ) \delta_{x_{k+1},0}
\Big ) \; .
\end{align*}
By examining the rule of multiplication in the crossed product algebra in \eqref{Eq-AlgOper}, one sees that the last expression is indeed the same as \eqref{Eq-LocalEven}.
\qed

\subsection{A generalized Connes-Chern character: the odd case}\label{Sec-OddConnesChern}

In this section $k$ is odd. Hence the Clifford algebra $\CM_k$ has two distinct irreducible representations on $\CM^{2^\frac{k-1}{2}}$. As in \cite{ProdanSpringer2016ds}, we choose the representation $\gamma_1,\ldots,\gamma_k$ satisfying $\gamma_1\cdots\gamma_k=\imath^k\,\mbox{\rm id}$.
%
%
The relevant Hilbert C$^\ast$-module for  the calculation below is 
$$
E_{\mathcal B_{(1)}} 
\;= \;\Big (\CM^{2^\frac{k-1}{2}} \oplus  \CM^{2^\frac{k-1}{2}} \Big ) \hat \otimes_\CM \mathcal H_\mathcal B 
\;\cong\; 
\big (\CM^{2^\frac{k-1}{2}} \otimes_\CM \mathcal H_\mathcal B \big ) \oplus  \big (\CM^{2^\frac{k-1}{2}}\otimes_\CM \mathcal H_\mathcal B \big ) 
\;\cong\; 
E_\mathcal B \oplus E_\mathcal B
\;,
$$
where similar as above $E_\mathcal B \;=\; \CM^{2^\frac{k-1}{2}}\otimes_\CM \mathcal H_\mathcal B$. The grading operator in the last representation is the swap $S(\psi \oplus \psi') = \psi' \oplus \psi$. As the notation suggests, this will be viewed as a graded Hilbert module over the graded C$^\ast$-algebra $\mathcal B_{(1)}=\Bb\hat{\otimes}\CM_1$. 

\vspace{.2cm}

On $E_\mathcal B$ the representation $\pi = 1 \hat \otimes_\CM \eta$ of the crossed product algebra $\mathcal A = \mathcal B \rtimes_\xi^\theta \ZM^k$ is still given, and it trivially extends to a graded representation (of degree $0$) to $E_{\mathcal B_{(1)}} $, which we still denote by $\pi$. On $E_\mathcal B$ act the adjointable operators
$$
\DiracSA_{x_0}
\;=\;
\gamma \hat \otimes_\CM \widehat{X + x_0} \otimes_\CM 1
\;,
\qquad x_0 \in (0,1)^k
\;,
$$
defined as in Definition~\eqref{def-DiracDef} by contraction over the indices. The operators $\DiracSA_{x_0}$ are selfadjoint and satisfy 
$\DiracSA_{x_0}^2=1$. However, there is no natural grading making them odd for odd $k$. However, proceeding as in Proposition~\ref{pr-KK1Cycle1}, one obtains a candidate for a  $(\mathcal A,\Bb_{(1)})$-cycle by considering $\DiracSA_{x_0} \oplus (-\DiracSA_{x_0})$ on $E_{\mathcal B_{(1)}} $. Arguing similarly as in the proof of Theorem~\ref{Th-EvenKKCycle}, one readily deduces the following. 

\begin{theorem}\label{Th-OddKKCycle} Let $x_0 \in (0,1)^k$. The triples 
$$
\mathcal E^1_{x_0} 
\;=\;
\big  ( E_{\mathcal B_{(1)}},\pi\hat{\otimes}1 ,\DiracSA_{x_0}\hat{\otimes} \epsilon\big  )
\;\cong\;
\big  ( E_\mathcal B \oplus E_\mathcal B, \pi \oplus \pi, \DiracSA_{x_0} \oplus (-\DiracSA_{x_0}) \big )\;, 
\;,
$$
define a field of Kasparov $(\mathcal A,\Bb_{(1)})$-cycles.
\end{theorem}

Proposition~\ref{Pr-TraceClass} and the estimates of Theorem~\ref{Th-EvenSummable} apply verbatim, so that:

\begin{theorem}\label{Th-OddSummable} 
The Kasparov cycles defined in Theorem~\ref{Th-OddKKCycle} are $n$-summable over the smooth sub-algebra $\mathscr A$ for any $n \geq k+1$, that is,
$$
\prod_{i=1}^n \big [\DiracSA_{x_0},\pi(a_i) \big ] \;\in\; \SM_1(E_\mathcal B) \; ,
$$
for any $a_i \in \mathscr A$, $i=1,\ldots,n$. 
\end{theorem}

In conclusion, the results of Section~\ref{Sec-OddGChern} can be applied and lead to the following.

\begin{corollary}\label{Co-OddChernConnesCh} 
The multilinear map
\begin{equation}\label{Eq-OddChernConnesCh}
\uptau_k(a_0,\ldots, a_k) 
\;=\; 
\Gamma_k \;{\rm Tr}\Big (\DiracSA_{x_0} \prod_{i=0}^k \big [\DiracSA_{x_0}, \pi(a_i)\big ] \Big )
\end{equation}
is a well-defined cyclic $(k+1)$-cocycle over the smooth sub-algebra $\mathscr A$. Its class in the periodic cyclic cohomology of $\mathscr A$ defines the Connes-Chern character of the Kasparov cycle $\mathcal E^1_{x_0}$. It fulfills the following index theorem
$$
 \Lambda_k \, \uptau_k(u^\ast,u,\ldots,u^\ast, u)
\;=\; {\rm Tr} \Big ( {\rm Ind}\big (P_{x_0}\pi(u) P_{x_0} \big ) \Big )
\;,
$$
for any unitary element $u \in \mathscr A$ and where $P_{x_0} = \tfrac{1}{2}(1+\DiracSA_{x_0})$.
\end{corollary}

\begin{theorem}[Odd local formula] \label{Th-OddLocalFormulaTh}
The generalized Connes-Chern characters of the Kasparov cycles $\mathcal E^1_{x_0}$ coincide and can be represented by a cyclic cocycle which accepts the following local formula:
\begin{equation}\label{Eq-LocalOdd}
\bar \uptau_k( a_0,\ldots,a_k)
\;=\; \Delta_k \sum_{\rho \in S_k} (-1)^\rho \mathcal T \Big (a_0  \prod_{i=1}^k \partial_{\rho_i} a_i \Big )
\;, 
\qquad 
\Delta_k
\;=\;
-4\imath\,(2\imath \pi)^\frac{k-1}{2}\; , 
\end{equation}
for any $a_i \in \mathscr A$, $i=0,\ldots,k$.
\end{theorem}

\proof The first part of the affirmation follows as in Theorem~\ref{Th-EvenLocalFormulaTh} and $\bar \uptau_k$ is defined as in \eqref{Eq-AverageCycle}, and again we set $\bar \uptau_k=\bar \uptau_k( a_0,\ldots,a_k)$. Let us now apply Proposition~\ref{Pr-TraceClass}(iv) on \eqref{Eq-OddChernConnesCh}. By opening the first commutator, 
$$
\bar \uptau_k
\;=\;
\Gamma_k\; \int\limits_{[0,1]^k} d x_0\sum_\alpha \sum_{x \in \mathbb Z^k} \Tt_\Bb \Big ( \Big (\big (\pi(a_0)  - \DiracSA_{x_0} \pi(a_0) \DiracSA_{x_0} \big ) \prod_{i=1}^{k} \big [\DiracSA_{x_0},\pi(a_i)\big]\Big )_{\alpha,x;\alpha, x} \Big ) \; .
$$
As in the even case, using \ref{Eq-Id07} as well as the fact that $k$ is odd, one can push the two $\DiracSA_{x_0}$ operators flanking $\pi(a_0)$ to the edges where they annihilate. Thus
\begin{align*}
\bar \uptau_k & 
\;=\; 2\,\Gamma_k \int\limits_{[0,1]^k} d x_0\sum_{x \in \ZM^k} \Tt_\Bb \Big ( \Big (\pi(a_0) \prod_{i=1}^{k} \big [\DiracSA_{x_0},\pi(a_i) \big ]\Big )_{\alpha,x;\alpha,x} \Big ) \; .
\end{align*}
In the following, ${\rm tr}$ will represent the trace over $\CM(2^\frac{k-1}{2})$. Using \ref{PP32}, 
\begin{align*}
\bar \uptau_k  \;=\;
& 2 \,\Gamma_k \, \int\limits_{[0,1]^k} d x_0 \sum_{x,x_1,\ldots,x_{k+1} \in \ZM^k} \delta_{x,x_{k+1}} e^{\imath (x,\theta x_1)}  \Big (\prod_{i=1}^k e^{\imath (x_i,\theta x_{i+1})} \Big )
\\
& 
\;\;
\times \, {\rm tr} \Big (\prod_{i=1}^k \big ( \gamma,  \widehat{x_0 + x_i} - \widehat{x_0 + x_{i+1}} \big ) \Big ) \Tt_\Bb \Big ( \xi_{-x_1}(b^{(0)}_{x-x_1})\prod_{i=1}^k  \xi_{-x_{i+1}}( b^{(i)}_{x_i - x_{i+1}}) \Big )
\;.
\end{align*}
Re-define the summation variables as $x_i \rightarrow x_i + x$, $i=1,\ldots,k+1$ and using the invariance of $\mathcal T_\mathcal B$ under the automorphisms $\xi_x$, we obtain
\begin{align*}
\bar \uptau_k  
\;=\;
& 
2\, \Gamma_k \, \int\limits_{[0,1]^k} d x_0 \sum_{x,x_1,\ldots,x_{k+1} \in \ZM^k} \delta_{x_{k+1},0} \, \Big(\prod_{i=1}^{k-1} e^{\imath (x_i,\theta x_{i+1})} \Big ) \\
&
\;\;\times \, 
{\rm tr} \Big (\prod_{i=1}^k \big ( \gamma,  \widehat{x_0 +x+ x_i} -\widehat{x_0 + x+ x_{i+1}} \big ) \Big )   \Tt_\Bb \Big ( \xi_{-x_1}(b^{(0)}_{-x_1})\prod_{i=1}^k  \xi_{-x_{i+1}} (b^{(i)}_{x_i - x_{i+1}}) \Big )
\;.
\end{align*}
Now combining the summation over $x$ and the integration over $x_0$ to an integration over $\RM^k$ leads to
\begin{align*}
\bar \uptau_k  
\;= \; & 2\, \Gamma_k \sum_{x,x_i \in \ZM^k} \delta_{x_{k+1},0}\, \Big(\prod_{i=1}^{k} e^{\imath (x_i,\theta x_{i+1})} \Big ) \\
&  \times \int\limits_{\RM^k} d x \, {\rm tr} \Big ( \prod_{i=1}^k \big ( \gamma,  \widehat{x+ x_i} -\widehat{x+ x_{i+1}} \big ) \Big )  \Tt_\Bb \Big ( \xi_{-x_1}(b^{(0)}_{-x_1})\prod_{i=1}^k  \xi_{-x_{i+1}} (b^{(i)}_{x_i - x_{i+1}} )\Big )
\;.
\end{align*}
At this point let us use the identity discovered in \cite{PS,ProdanSpringer2016ds},
\begin{align*}
\delta_{x_{k+1},0}\int\limits_{\mathbb R^k} d x \, {\rm tr} \Big ( \prod_{i=1}^k \big (\gamma, \widehat{x+x_i} -\widehat{x + x_{i+1}}) \Big ) 
& 
\;=\; \tilde \Lambda_k \;\sum_{\rho \in S_k} (-1)^\rho \prod_{i=1}^k (x_i)_{\rho_i}
\\
& \;=\; \delta_{x_{k+1},0}\,\tilde \Lambda_k \;\sum_{\rho \in S_k} (-1)^\rho \prod_{i=1}^k (x_i-x_{i+1})_{\rho_i}
\;,
\end{align*}
with $\tilde \Lambda_k = -\, \frac{2^k(\imath \pi)^\frac{k-1}{2}}{k!!} $. As $2\Gamma_k \tilde \Lambda_k = \imath^k \Delta_k$, one thus finds
\begin{align*}
\bar \tau_k & \;=\; 
\imath^k \Delta_k \,\sum_{\rho \in S_k} (-1)^\rho \sum_{x_i \in \mathbb Z^k}  \Tt_\Bb \Big ( b^{(0)}_{-x_1} \prod_{i=1}^k e^{\imath (x_i,\theta x_{i+1})}  \xi_{- x_{i+1}}\big ( (x_i-x_{i+1})_{\rho_i}b^{(i)}_{x_i - x_{i+1}} \big )\Big ) \\
& \;=\; \Delta_k\, \sum_{\rho \in S_k} (-1)^\rho \sum_{x_i \in \mathbb Z^k}  \Tt_\Bb \Big ( b^{(0)}_{-x_1} \prod_{i=1}^k e^{\imath (x_i,\theta x_{i+1})}  \xi_{- x_{i+1}}\big ( \partial_{\rho_i}b^{(i)}_{x_i - x_{i+1}} \big )\Big ) \; ,
\end{align*}
and the statement follows.\qed

\section{Application to topological insulators}
\label{Sec-ApplicationTI}

This section presents an application of the previous results to solid state systems describing topological insulators. While this was our original motivation, we attempt to condense the discussion as much as possible and refer the reader to our recent work \cite{ProdanSpringer2016ds} for  physical motivation and a detailed mathematical description of the formalism.

\subsection{Disordered crystals}  

\begin{definition}[\cite{BellissardLNP1986jf,ProdanSpringer2016ds}]
Let $\Omega$ be a compact topological space equipped with a $\ZM^d$-action $\tau$ and an invariant and ergodic probability measure $\PM$. This ergodic dynamical system $(\Omega,\tau,\ZM^d,d\PM)$ encodes the disorder or quasicrystaline configurations. A  homogeneous crystal is a family $H =\{H_\omega\}_{\omega \in \Omega}$ of quantum mechanical Hamiltonians given by selfadjoint operators $H_\omega$ on the physical Hilbert space $\ell^{2}(\ZM^{d}) \otimes \CM^{N}$ of the form
\begin{equation}\label{Eq-GenHam}
H_{\omega} \;=\;
\sum_{q\in R} \;\sum_{x \in \ZM^d} \, U_q|x \rangle \langle x| \otimes W_q(\tau_{-x} \omega)\;, 
\qquad W_q \in C\big (\Omega, \CM(N) \big )
\;. 
\end{equation}
The set $R \subset \ZM^d$ is finite and it is called the hopping range. The operators $U_q$ are the dual magnetic translations on $\ell^2(\ZM^d)$ defined by $U_q \delta_x = e^{\imath (q,\upphi x)} \delta_{x+q}$, where the entries in the skew symmetric matrix $\upphi$ are directly related to the magnetic fluxes through the facets of the unit cell of the crystal. 
\end{definition}

\begin{remark} {\rm The magnetic translations are defined by $V_q \delta_x = e^{-\imath (q,\upphi x)} \delta_{x+q}$ and they commute with the dual magnetic translations. The homogeneity property of $H_\omega$ is reflected in the covariant property $V_q H_\omega V_q^\ast = H_{\tau_q \omega}$.
}
\hfill $\diamond$
\end{remark}

\begin{example}\label{Ex-2DElec}{\rm The low energy physics of the 2-dimensional electron gas subjected to a perpendicular electric field, a periodic potential and white disorder is often analyzed using the following lattice model, 
$$
H_{\omega}=\sum_{j=1,2} \;\sum_{x \in \ZM^2} \, \big (1+\lambda_j \omega_j(x) \big )\Big ( e^{\imath \phi \, e_j \wedge x} |x+e_j \rangle \langle x| + e^{-\imath \phi \, e_j \wedge x} |x \rangle \langle x + e_j| \Big ) ,
$$
where $\phi$ is just a number, equal to the magnetic flux through the repeating cell (in appropriate physical units), $q \wedge x = q_1x_2 - q_2 x_1$, and $\omega_j(x)$ are random numbers drawn independently from the interval $\big [-\frac{1}{2},\frac{1}{2}\big ]$. This model is well known to harbor many topological quantum Hall phases, even at relatively large disorder (see {\it e.g.} \cite{ProdanAMRX2013bn} for numerical representations of its phase diagram). Note that in this model, $W_{e_j}(\omega) = 1+\lambda_j \omega_j(0)$, so that $W_{e_j}(\tau_x\omega) = 1+\lambda_j \omega_j(x)$.} 
\hfill $\diamond$
\end{example}

\begin{definition}\label{De-MobGap} A homogeneous crystal is said to posses a spectral gap, if there exists an open interval $\Delta$ such that the spectra of $H_\omega$ satisfy $\sigma(H_\omega) \cap \Delta = \emptyset$ for all $\omega \in \Omega$. If the chemical potential $\mu$ belongs to a spectral gap, the homogeneous crystal is said to be insulating.
\end{definition}

\begin{remark} {\rm The regime described above is known as the weak disorder regime. The strong disorder regime, which is not treated here, refers to situation where the chemical potential is located in the essential spectrum which is nevertheless Anderson localized (see \cite{ProdanSpringer2016ds} for details).
} 
\hfill $\diamond$
\end{remark}

\begin{definition} A homogeneous crystal posses the chiral symmetry if there exists a symmetry $J =  1 \otimes \Gamma$ on $\ell^2(\ZM^d)\otimes \CM^{2N}$, {\it i.e.} $J^\ast = J$ and $J^2=1$, such that
$$
JH_\omega J^{-1} 
\;=\; 
-H_\omega, \quad \omega \in \Omega
\;,
$$
and if the chemical potential $\mu$ is fixed at the origin. Without loss the generality, we can assume that $J={\rm diag}(1_N,-1_N)$. Note that a chiral symmetric $H$ can posses a spectral gap at the origin only if the dimension of the fiber is even. 
\end{definition}

\begin{remark}{\rm The ground state of a homogeneous rystal is encoded in the covariant family of Fermi projections:
\begin{equation}
\label{eq-FermiP}
P_F 
\;=\; 
\{ \chi(H_\omega \leq \mu) \}_{\omega \in \Omega}
\;.
\end{equation}
If the chiral symmetry is present and the system is insulating, then the unitary $1-2P_F$ (called flat band Hamiltonian) necessarily takes the form
\begin{equation}
\label{eq-FermiU}
1-2P_F 
\;=\; 
\begin{pmatrix} 0 & U_F^\ast \\ U_F & 0 \end{pmatrix}
\;, 
\qquad 
P_F 
\;=\; 
\frac{1}{2} \begin{pmatrix} 1 & -U_F^\ast \\ -U_F & 1 \end{pmatrix}
\;.
\end{equation}
As such, the ground state of a chiral symmetric system is encoded in the family of unitary operators $U_F=\{ U_{F,\omega} \}_{\omega \in \Omega}$, which are called the Fermi unitary operators.
}\hfill $\diamond$
\end{remark}

\subsection{Algebra of physical observables}
\label{Sec-PhysAlg}

To ease the notation, we will use $C_N(\Omega)$ instead of $C\big (\Omega, \CM(N) \big )$ for the C$^\ast$-algebra of continuous functions over $\Omega$ with values in the algebra of $N\times N$ complex matrices.

\begin{definition} 
\label{def-rotalg}
The algebra of the bulk observables is defined as the universal C$^\ast$-algebra generated by $C_N(\Omega)$ and unitaries $u_1,\ldots,u_d$,
$$
\mathcal A_d\;=\;C^\ast\big (C_N(\Omega),u_1,\ldots,u_d \big )\;,
$$
with the following commutation relations:
\begin{equation}
\label{BComm1}
u_i u_j\;=\;e^{\imath \phi_{i,j}}u_ju_i\;, 
\qquad  
u_j^\ast u_j \;=\; u_j^\ast u_j \;=\; 1\;, \qquad i,j=1,\ldots,d,
\end{equation}
as well as
\begin{equation}
\label{BComm2}
f\, u_j\; =\; u_j\,(f\circ \tau_j)\,, 
\qquad \ f \in C_N(\Omega)\,, \;\; j=1,\ldots,d\,.
\end{equation}
\end{definition} 

\begin{proposition}[Equivalent presentations]
\label{Pr-Presentations}
 The algebra of bulk observables admits the following equivalent presentation:
\begin{enumerate}[\rm (i)]
\item As twisted crossed product:
$$
\mathcal A_d \;=\; C_N(\Omega) \rtimes_\tau^\upphi \ZM^d
\;,
$$
where $u_x$ in Definition~\ref{Def-TwistedCP} are given by
$$
u_x \;=\; e^{\imath (x,\upphi_+ x)}u_1^{x_1} \cdots u_d^{x_d}
\;, 
\qquad 
u_x u_y \;=\; e^{\imath(x,\upphi y)} u_{x+y}
\;,
$$
and $\upphi_+$ is the lower triangular part of $\upphi$, and the action of $\ZM^d$ on $C_N (\Omega)$, also denoted by the symbol $\tau$, is given by 
$$
\tau_x(f) 
\;=\; 
u_x f u_x^\ast = f \circ \tau_x^{-1}
\;.
$$ 

\item As decomposed twisted crossed product:
$$
\mathcal A_d \;=\; \Aa_j \rtimes_{\xi}^{\uptheta} \ZM^k \;, 
\qquad \Aa_j \;=\; C_N(\Omega) \rtimes_{\sigma}^{\upvarphi} \ZM^j
\;, 
\qquad \ZM^d \;=\; \ZM^j \times \ZM^k
\;,
$$
where $\sigma = \tau|_{\ZM^j}$ and the matrix $\upvarphi$ is given by
$$
\varphi_{m,n} \;=\; \phi_{m,n}
\;, 
\qquad n,m = 1,\ldots,j
\;,
$$ 
while the matrix $\uptheta$
$$
\theta_{m,n} \;=\; 
\Phi_{j+m,j+n}
\;, 
\qquad n,m = 1,\ldots,k\;,
$$
and
\begin{equation}\label{Eq-Ux}
\xi_x \;=\; {\rm Ad}_{u_x}\in {\rm Aut}(\mathcal A_j)
\;, 
\qquad 
u_x \;=\; e^{\imath (x,\uptheta_+ x)}u_{j+1}^{x_1} \cdots u_{j+k}^{x_k}
\;, 
\qquad x \in \ZM^k
\;.
\end{equation}

\end{enumerate}
\end{proposition}

\proof All statements are direct consequences of the definition of the twisted crossed product given in Section~\ref{Sec-CrossedProd} as well as the standard decomposition of the twisted crossed products in \cite{PackerMPCPS1989fs}. \qed

\begin{remark}{\rm In the following, $\ZM^j$ and $\ZM^k$ will be viewed as the subgroups of $\ZM^d$ generated by $\{e_n\}_{n=1\ldots,,j}$ and  $\{e_n\}_{n=j+1,\ldots,d}$, respectively. As such, for $x \in \ZM^k$ the operator $u_x$  is automatically understood to be the one defined in \eqref{Eq-Ux}. } \hfill $\diamond$ \end{remark}

\begin{proposition}[Canonical Representation, {\it e.g.} \cite{ProdanSpringer2016ds}]\label{Pro-CanRep} The following relations 
$$
\uppi_\omega(u_j) \;=\; U_{e_j} \otimes 1\;, 
\qquad 
\uppi_\omega(f) 
\;=\; 
\sum_{x \in \ZM^d} |x \rangle \langle x| \otimes f(\tau_{-x} \omega)
\;, 
\qquad 
j=1,\ldots,d\;,\;\; f \in C_N(\Omega)\;,
$$
define a field of covariant and faithful $\ast$-representations 
$$
\uppi_\omega : \mathcal A_d \rightarrow \BM\big (\ell^2(\ZM^d) \otimes \CM^N \big )
\;, 
\qquad 
V_x \circ \uppi_\omega \circ V_x^{-1} \;=\; \uppi_{\tau_x \omega}
\;, 
\qquad \omega \in \Omega
\;.
$$
Any homogeneous crystals can be generated from this representation.  
\end{proposition}

\begin{example}\label{Ex-SmoothProj} {\rm The generic model in \eqref{Eq-GenHam} is generated by the element of $\mathcal A_d$
$$
h 
\;=\; 
\sum_{q \in R\subset \ZM^d}  u_q \, W_q\;, 
\qquad W_q \,\in\, C_N (\Omega)
\;.
$$ 
Note that $h$ is actually an element of the smooth algebra. In the weak disorder regime, this automatically implies that the Fermi projector $p_F = \chi(h \leq \mu)$ is also an element of the smooth algebra, because $p_F$ can be computed using the holomorphic functional calculus (as a Riesz projection). When the chiral symmetry is present, it also automatically follows that the element $u_F$ which generates the family of Fermi unitary operators is an element of the smooth algebra. \hfill $\diamond$ }
\end{example}

\subsection{Weak topological invariants in the weak disorder regime}


In this section, it is shown how the mathematical results of the paper can be applied to homogeneous crystals and allow to write out formulas for their weak topological invariants, as in \cite{ProdanSpringer2016ds}. Our new results then provide index theorems and local formulas for these invariants. In order to apply the results of Section~\ref{Sec-CrossedProd}, let first $C_N(\Omega)$ play the role of the C$^\ast$-algebra $\mathcal B$ and $\int_\Omega d\mathbb P (\omega){\rm tr}_N(.)$ the role of trace $\mathcal T_\mathcal B$. Also, $(\partial, \mathcal T)$ will represent the non-commutative calculus on $C_N(\Omega) \rtimes_\tau^\upphi \ZM^d$ and $\mathscr A_d$ the smooth sub-algebra w.r.t. this non-commutative calculus. To start out, let us write out two existing results in order to place our discussion in a proper context. The first assumes $\Omega$ to be contractible. This holds for typical disordered models with continuous distribution of the random variables, but {\it not} for quasicrystals. One of the main points below is that we are able to lift this hypothesis.

\begin{proposition}\label{Pro-CyclicCo} Assume $\Omega$ to be a contractible topological set. Then the $K$-theory of $\Aa_d$ coincides with the $K$-theory of the $d$-dimensional non-commutative torus. Furthermore, the periodic cyclic cohomology of the smooth algebra $\mathscr A_d$ is generated by the cyclic cocycles
\begin{equation}\label{Eq-ChernCoCy}
\zeta_I(a_0,a_1,\ldots,a_{|I|}) 
\;=\; 
\Delta_{|I|} \sum_{\rho \in S_{|I|}} (-1)^\rho \, \mathcal T \Big (a_0 \prod_{i=1}^{|I|} \partial_{\rho_i} a_i \Big )
\;, 
\qquad 
I \subset \{1,\ldots,d\}
\;,
\end{equation}
where the bijections $\rho\in S_{|I|}$ are viewed as maps from $\{1,\ldots,|I|\}$ to the ordered set $I$.
\end{proposition} 

\begin{remark}{\rm The first statement is standard (see \cite{WeggeOlsenBook1993de}). The second statement is due to Nest \cite{NestJFA1988fj,NestCJM1988fj}. The next result follows directly from Connes' book \cite{Connes:1994wk}, see also \cite{ProdanSpringer2016ds}. It does not pend on $\Omega$ being contractible.}
\hfill $\diamond$
\end{remark}

\begin{proposition}\label{Pro-NumericalInv} 
Suppose given a homogeneous crystal with spectral gap and let $p_F$ be the Fermi projection given in \eqref{eq-FermiP} and, in case tha the crystal has a chiral symmetry, $u_F$ the Fermi unitary given by \eqref{eq-FermiU}. Then the pairings of $K$-theory and cyclic cohomology applied to their associated classes in $K$-theory,
$$
\big \langle [p_F]_0, [\zeta_{I}] \big \rangle 
\;=\; 
\zeta_{I}(p_F,\ldots,p_F)
\;, 
\qquad |I|\mbox{ even} 
\;,
$$
and 
$$
\big \langle [u_F]_1, [\zeta_{I}] \big \rangle 
\;=\; \zeta_{I}(u_F^\ast,u_F,\ldots,u_F^*,u_F)\;, 
\qquad |I|\mbox{ odd} 
\;,
$$
generate numerical invariants for the homogeneous crystals which are stable at least in the regime of weak disorder. 
\end{proposition}

\begin{remark}{\rm The physical interpretation of these invariants was established in \cite{ProdanSpringer2016ds}. For $|I|=\mbox{even}$, the invariants are directly related to the linear and non-linear magneto-electric transport coefficients of the crystal, while for $|I|=\mbox{odd}$ they are directly related to the chiral electric polarization and its derivatives w.r.t. the magnetic field. \hfill $\diamond$
}
\end{remark}

\begin{remark}{\rm The range of the numerical invariants of Proposition~\ref{Pro-NumericalInv} was computed in \cite{ProdanSpringer2016ds} under the hypothesis that $\Omega$ is contractible. It was also established that, still under this hypothesis, that the invariants pinpoint uniquely the class of $p_F$ in $K_0(\mathcal A_d)$ and that of $u_F$ in $K_1(\mathcal A_d)$. \hfill $\diamond$
}
\end{remark}


The following is our new contribution. By re-labeling the space coordinates, one can always assume that the set $I \subset \{1,\ldots,d\}$ defining the cyclic cocycles in Proposition~\ref{Pro-CyclicCo} is of the form $I=\{d-k+1,\ldots,d\}$, with $k=|I|$. Let us now consider the presentation (ii) of Proposition~\ref{Pr-Presentations} of the algebra of physical observables, namely $\mathcal A_d = \Aa_j \rtimes_{\xi}^{\uptheta} \ZM^k$. Thus let us take $\mathcal B=\mathcal A_j= C_N(\Omega) \rtimes^{\upvarphi}_{\sigma} \ZM^j$ with canonical trace $\mathcal T_\mathcal B$ from Definition~\ref{Def-CanTrace} now induced by the trace $\mathcal T_{C_N(\Omega)}(.) = \int_\Omega d\PM(\omega){\rm tr}_N(.)$ on $C_N(\Omega)$. Following again the standard procedure from Definitions~\ref{Def-CanTrace} and \ref{Def-CanDer}, we define next the non-commutative calculus on $\mathcal A_d =\mathcal B \rtimes_\xi^\uptheta \ZM^k$. The procedure leads to the same trace $\Tt$ on $\mathcal A_d$ and to the derivations $\partial_{1}, \ldots, \partial_k$ which actually coincide with the derivations $\partial_{i_1}, \ldots, \partial_{i_k}$, $\{i_1,\ldots,i_k\}=I$. As a consequence, the cyclic cocycles defined by the local formulas in Theorems~\ref{Th-EvenLocalFormulaTh} and \ref{Th-OddLocalFormulaTh} coincide with the cocycles $\zeta_I$ defined in Proposition~\ref{Pro-CyclicCo}. Then, in the notations from Sections~\ref{Sec-EvenConnesChern} and \ref{Sec-OddConnesChern}, we can state at once the following results, which, as already mentioned above, do not require $\Omega$ to be contractible as in \cite{ProdanSpringer2016ds}.

\begin{corollary}[Index theorems for weak topological invariants]  The class of $\zeta_{I}$ in the periodic cyclic cohomology is equal to the generalized Connes-Chern character canonically associated to $\mathcal A_j \rtimes_\xi^\uptheta \ZM^k$. Consequently, the numerical topological invariants defined in Proposition~\ref{Pro-CyclicCo} accept a generalized index formula:
$$
\Lambda_k \, \zeta_{I}(p_F,\ldots,p_F) 
\;=\; 
{\rm Tr} \Big ( {\rm Ind}\big(\pi(p_F) F_{x_0} \pi(p_F) \big ) \Big )
\;, 
\qquad |I|\mbox{ even}
\;,
$$
and, if the chiral symmetry is present,
$$
\Lambda_k \, \zeta_{I}(u_F^\ast,\ldots,u_F) \;=\; {\rm Tr} \Big ( {\rm Ind}\big (P_{x_0} \pi(u_F) P_{x_0} \big ) \Big )
\;, 
\qquad |I|\mbox{ odd}
\;.
$$
\end{corollary}

\begin{theorem}[Range of weak topological invariants]\label{Th-Range} The topological invariants defined in Proposition~\ref{Pro-NumericalInv} take values in the image of $K_0(\mathcal A_j)$ through the trace ${\rm Tr}$,
$$
{\rm Ran}(\zeta_I)\; \subset \;{\rm Tr}\big ( K_0(\mathcal A_j) \big )
\;,
$$
regardless of the parity of $|I|$.
\end{theorem}

\noindent {\bf Acknowledgements:} The work of E.~P. was partially supported by the U.S. NSF grant DMR-1056168, that of H.~S.-B. by the DFG grant SCHU-1358/6. E.~P. also wants to acknowledge the hospitality of the Department of Mathematics of Friedrich-Alexander Universit\"at, where a large part of this work was completed.



\end{document}